\documentclass{olplainarticle}

\usepackage{amssymb}
\usepackage{amsmath}
\usepackage{placeins}
\usepackage{graphicx}
\usepackage[colorlinks=true, linkcolor=blue, citecolor=blue, urlcolor=blue]{hyperref}
\usepackage{gensymb}
\usepackage{caption,subcaption}
\usepackage{tabularray}
\usepackage[acronym, nonumberlist]{glossaries}
\usepackage{xcolor}
\usepackage{graphicx}
\usepackage[T1]{fontenc}   
\usepackage[utf8]{inputenc} 
\usepackage{titlecaps}

\title{ADJOINT-BASED RECOVERY OF THERMAL FIELDS FROM DISPLACEMENT OR STRAIN MEASUREMENTS}

\author[1]{Talhah Shamshad Ali Ansari}
\author[2]{Rainald Löhner}
\author[1]{Roland Wüchner}
\author[3]{Harbir Antil}
\author[1]{Suneth Warnakulasuriya}
\author[4,5]{Ihar Antonau}
\author[2]{Facundo Airaudo}
\affil[1]{Chair of Structural Analysis, Technical University of Munich, Munich, Germany}
\affil[2]{Center for Computational Fluid Dynamics and Department of Physics and Astronomy, George Mason University, Fairfax, VA 22030, USA}
\affil[3]{Center for Mathematics and Artificial Intelligence (CMAI) and Department of Mathematical Sciences, George Mason University, Fairfax, VA 22030, USA}
\affil[4]{Cluster of Excellence SE$^{2}$A -- Sustainable and Energy-Efficient Aviation, Technische Universit\"at Braunschweig, Germany}
\affil[5]{Institute of Structural Analysis, Technische Universit\"at Braunschweig, Braunschweig, Germany}


\begin{abstract}
A finite-element method dependant adjoint-based procedure to determine the temperature field of structures based on measured displacements/strains and a set of standard loads is developed and tested. Given a series of force and deformation measurements, the temperature field is obtained by minimizing the adequately weighted differences between the measured and computed values. Three numerical examples — a Plate With a Hole, a Bridge, and a Hoover Dam example — each with multiple sensors distributed in different configurations, demonstrate the procedure's capabilities. A target temperature distribution is prescribed in all cases, and the displacement sensor data is recorded. The optimization algorithm (here, steepest descent with Barzilai-Borwein step) uses this data to optimize the temperatures such that the same deformation is obtained at the sensor locations. Vertex Morphing is used as a filter to mitigate the ill-conditioning. Results show that the proposed approach can accurately reconstruct the target thermal distribution, especially when more sensors are used. Additionally, it is observed that the sensors do not need to be positioned in the region of interest; the method remains effective as long as the sensors can detect changes related to that area. A comparison with standard spatial interpolation techniques, namely, k-nearest neighbors and ordinary and universal kriging, is performed using temperature sensors in the same configurations. The proposed approach performs remarkably better than the interpolation techniques with a reduction in the root-mean-squared error of up to 38.4\%, 94\%, and 40\%, for the Plate With a Hole, the Bridge, and the Dam examples, respectively.
\end{abstract}

\begin{document}

\flushbottom
\maketitle
\thispagestyle{empty}
\noindent
\printglossary[type=\acronymtype]

\section{Introduction}
\label{sec:introduction}

There is a great need to accurately model the thermal distribution of structures, especially those exposed to uneven solar radiation. Such typical structures include bridges where the deck and the soffit are exposed to sunlight differently (\cite{kulprapha2012structural}); for long bridges, temperatures may have non-uniformity even in the longitudinal direction (\cite{ma2023statistical}). Large structures like dams are also subjected to non-uniform irradiation on the upstream and downstream surfaces. The non-uniform thermal distribution of the structure is further complicated by shadow and other environmental effects like wind, humidity, cloud cover, time of the day influencing the angle of incidence, season, etc.

\cite{bayraktar2022long} studied the yearly and daily strain-temperature behavior on the deck, pylon, and cables in cable-stayed bridges taking the Kömürhan cable-stayed bridge in Turkey as an example. They observed a strong correlation between the strain in the main steel deck and the air temperature for both yearly and daily monitoring, the mid-section concrete pylon leg strains on all four sides were different, and the cable forces also showed a direct correlation with temperature on yearly variation but the daily variation showed an opposite correlation between cable force and temperature.

Temperature monitoring is not only required during the structure's service life but may also be critical during construction. 
During the construction phase of massive concrete structures, factors such as hydration heat, external environmental conditions, cooling technique, etc., can influence the temperature distribution within the structure and thereby lead to crack development. Hence accurate temperature monitoring and control is required to adequately address these issues (\cite{lin20213d, liu2015precise,zhou2019temperature}). Additionally, temperature at different dam sections may be controlled by different factors. For instance, the ambient air temperature and solar radiation are the dominant factors influencing the temperature of the upstream, crest, and exposed sections of the dam, while the upstream submerged section temperature is mainly influenced by the water temperature (\cite{bui2019evaluation}).

Thermal-induced effects are one of the significant loads under regular operation for bridges (\cite{bayraktar2022long,glashier2024temperature}). The influence on structural response due to environmental (mainly temperature) conditions could even be greater than that caused by live loads and damage (\cite{glashier2024temperature};~\cite{kromanis2014predicting};\ \cite{zhou2014summary}). 
If the thermal loads are not adequately monitored, modeled, and integrated into the analysis and long-term monitoring framework, it may lead to incompetent digital twins and faulty Structural Health Monitoring (\acrshort{SHM}) systems. 

With the increasing demand for accurate digital twins (\cite{HAntil_2024a}) of civil infrastructure usually integrated with sensors enabling long-term monitoring, there is a necessity of separating the thermally induced response changes from the structural aging and weakening to get a realistic status of the structure. With regard to practicality, only a finite number of sensors are placed in/on a structure. With optimized sensor placement using a sparse network of sensors, the exclusion of some sensors due to improper installation or failure (as in the case of MX3D bridge with 10 excluded sensors out of a total of 84 installed sensors (\cite{glashier2024temperature})) could lead serious lack of information in that region. As proposed in this work, there may arise a need to calculate the temperature distribution of the structure using the deformation sensors either due to a lack of adequate or reliable temperature data. Therefore, the deformation sensors could play a secondary role when not detecting weakness-induced deformations, particularly during the healthy early life of the structure (\cite{kromanis2013support}).

This paper is structured as follows: Section \ref{sec:thermal field reconstruction} presents a review of the existing research related to thermal field reconstruction of structures with a hint towards Structural Health Monitoring (\acrshort{SHM}). Design specifications and different techniques available in the literature are discussed. Section \ref{sec:methodology} explains the proposed methodology including the thermal stress consideration, the thermal field determination, the optimization problem, and the Vertex Morphing filtering technique. Section \ref{sec:examples} presents three numerical examples: a Plate with a Hole, a Bridge, and a simplified FE model of the Hoover Dam. All the examples are tested with different numbers of sensors and sensor locations, with and without Vertex Morphing filtering. Section \ref{sec:spatial interpolation} compares the results obtained from the proposed approach with those obtained from k-nearest neighbors (\acrshort{kNN}) interpolation, ordinary kriging (herein abbreviated as \acrshort{OK}), and universal kriging (herein abbreviated as \acrshort{UK}). The spatial interpolation techniques take as input the temperature sensor data from the same locations as the displacement sensors used in the proposed methodology, i.e., the number and location of measurement points remain the same, and only the type of data (displacement or temperature) used, determines the type of sensor. Section \ref{sec:conclusion} presents a summary of the work and a discussion of the results. A brief overview of the limitations and future research is also provided.


\section{Thermal field reconstruction from a Structural Health Monitoring (\acrshort{SHM}) perspective}
\label{sec:thermal field reconstruction}


\begin{figure}[!b]
\begin{minipage}[c][][t]{\textwidth}
  \centering
  \includegraphics[width=0.4\paperwidth]{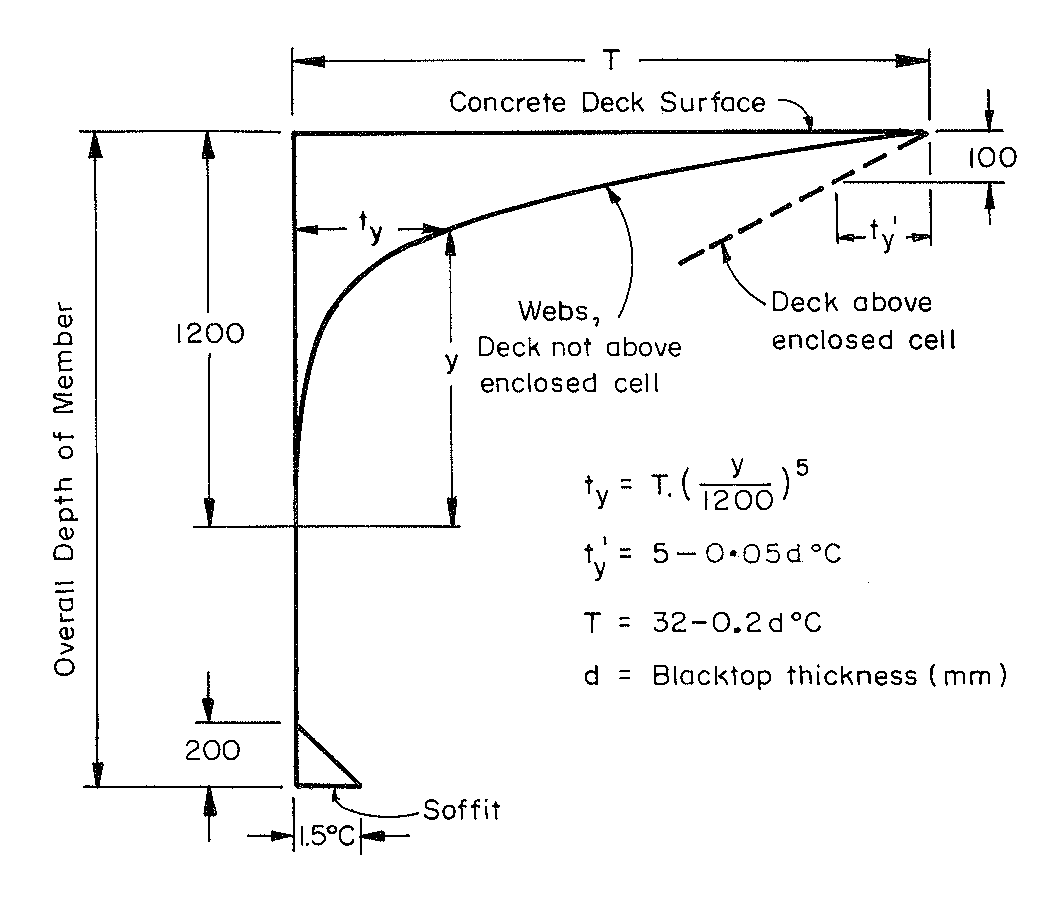}
  \caption{Temperature variation with depth according to the \cite{waka2013bridge}. \textbf{From the Waka Kotahi NZ Transport Agency Bridge Manual SP/M/022, Third Edition, Amendment 4. Used with permission.}}
  \label{fig:Figure_1}
\end{minipage}%
\end{figure}

Extensive research has been carried out in the past few decades to study the thermal distribution and behavior of civil structures. 
For most bridges, the temperature in the longitudinal and transverse directions is considered to be uniform, and variation in the vertical direction is mainly studied.   
For example, the New Zealand Code (The Waka Kotahi NZ Transport Agency Bridge manual SP/M/022 (\cite{waka2013bridge})) considers the temperature varies from the maximum temperature at the deck surface to a minimum over a depth of $1200$ $mm$ of the structure following a fifth-order parabolic curve, where the maximum deck surface temperature depends on the blacktop thickness (asphalt). The soffit temperature is assumed to be $1.5$ $\degree C$ above the minimum and decreases linearly over a depth of $200$ $mm$ into the structure. In case the structure depth is less than $1400$ $mm$, both curves are superimposed. The temperature variation with depth according to the New Zealand Code is shown in Figure \ref{fig:Figure_1}.

\begin{figure}[!t]
\begin{minipage}[c][][t]{.5\textwidth}
  \vspace*{\fill}
  \centering
  \includegraphics[width=0.95\textwidth]{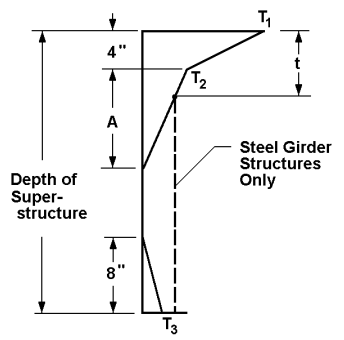}
  \subcaption{Positive vertical temperature gradient in Concrete and Steel Superstructures.}
  \label{fig:Figure_2a}
\end{minipage}%
\begin{minipage}[c][][t]{.5\textwidth}
  \vspace*{\fill}
  \centering
  \includegraphics[width=0.95\textwidth]{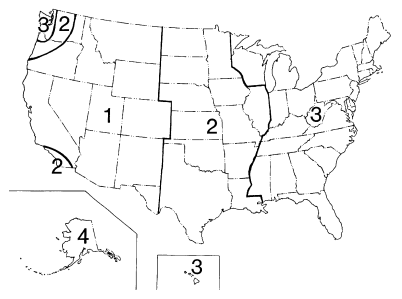}
  \subcaption{Solar radiation zones for the United States}
  \label{fig:Figure_2b}\par
  \vspace{1em}
  \includegraphics[width=0.95\textwidth]{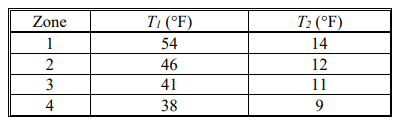}
  \subcaption{Basis for temperature gradients}
  \label{fig:Figure_2c}
\end{minipage}
\caption{Vertical thermal design gradient according to the \cite{aashto2020lrfd}. \textbf{From AASHTO LRFD Bridge Design Specifications, 2020, published by the American Association of State Highway and Transportation Officials, Washington, D.C., USA. Used with permission.}}
\label{fig:Figure_2}
\end{figure}

The \acrshort{AASHTO LRFD} (American Association of State Highway and Transportation Officials Load and Resistance Factor Design) Bridge Design Specification has a relatively simpler vertical temperature gradient design specification considering a bi-linear curve in the top section of the girder and a linear curve in the bottom section. For a positive vertical gradient (top surface warmer than the bottom surface), the temperature decreases linearly for a depth of 4 $inches$ from the top maximum temperature ($T_1$) to an intermediate temperature ($T_2$); a further linear decrease to the minimum temperature happens over a certain depth $A$, where the distance $A$ is taken as 12 $inches$ for concrete superstructures with depth more than 16 $inches$, and 4 $inches$ less than the actual depth if the depth is less than 16 $inches$. The linear variation from the soffit temperature ($T_3$) happens over a depth of 8 $inches$ into the structure. The temperature gradient profile according to the AASHTO LRFD 2020 specification is shown in Figure \ref{fig:Figure_2a}. In 2007, the temperature gradient profile was simplified to a bi-linear curve in the upper section as shown in Figure \ref{fig:Figure_2a} from a previous tri-linear curve. The temperatures $T_1$ and $T_2$ considered for the temperature gradient basis depend on the geographical location. The AASHTO LRFD 2020 specification divides the United States into four zones depending on the solar radiation. The US solar radiation zones and the corresponding temperature gradient basis are shown in Figures \ref{fig:Figure_2b} and \ref{fig:Figure_2c}, respectively. Unless otherwise specified, the bottom temperature $T_3$ is generally taken as $0$ $\degree F$.

\cite{hagedorn2019temperature} presented their experimental findings regarding the temperature gradients in concrete bridge girders subjected to environmental conditions. They concluded that the maximum ambient temperatures did not cause a maximum thermal gradient in the girders, but rather, higher temperature fluctuations during the day caused higher thermal gradients. They also observed a noticeable difference between the measured and AASHTO LRFD predictions for the vertical temperature gradient profile. Along the transverse direction, a temperature difference of up to 20 $\degree C$ was observed. They also compared the measured temperature rise to the theoretical rise (calculated from the strain measurements assuming any expansion was solely due to temperature change) and found reasonable agreement between the two. However, there were some deviations in the lower section of the girder, which they attributed to errors in the two methods. Most importantly, the strain was calculated on the outer surface while the temperature was measured in the middle section.
\cite{hagedorn2019temperature} recommended more detailed thermal vertical and thermal gradient models for the AASHTO Type II and Type IV girders to improve upon the AASHTO LRFD specifications.

Regarding the longitudinal temperature distribution, \cite{ma2023statistical} conducted statistical analyses of the longitudinal temperature distribution for a steel box girder long-span suspension bridge using one-year temperature monitoring data. They observed the probability density functions of the girder surface temperatures at different sections in the longitudinal direction were different. Upon further inspection, they noted that the left and right span of the bridge, while symmetric about the middle section, had different temperatures at different sections. \cite{ma2023statistical} concluded that the girder surface temperature distribution was three-dimensional non-uniform in both space and time.
These observations call for the need for more detailed temperature distribution analysis in structures to avoid making wrong assumptions.

In general, the thermal field of structures can be reconstructed in multiple ways (\cite{lin20213d}): (1) the simplest approach is by using thermometers to measure the internal and/or surface temperatures at limited discrete locations. Infrared cameras can provide continuous temperature distribution on the surface. More recently, optical fibers have gained popularity due to the large amount of measurement data available along the length of the fiber, (2) another approach is by solving the thermal heat equations analytically. Although, suitable for simple geometries, they can run into difficulties for complex geometries and boundary conditions (\cite{santillan2015new}), (3) another popular approach is the use of numerical analysis such as Finite Element Method (\acrshort{FEM}) (\cite{bui2019evaluation, castilho2018fea,chen2019temperature, yang2012fem, zhouc2016optimization}), extended Finite Element Method (XFEM) (\cite{zuo2015extended}), Composite Element Method (CEM) (\cite{zhong2017improved}), etc to simulate the temperature profile. 

The article by \cite{jin2010practical} presented a practical approach to predict the thermal field of the exposed face of arch dams accounting also for solar radiation and shadow effects. It observed non-uniform temperature distribution on the surface due to shadow and irregular irradiation. 
\cite{santillan2015new} introduced a 1-D analytical solution of the heat diffusion equation to compute the temperature field of concrete dams, also accounting for environmental conditions.

\cite{tong2001numerical} developed a numerical modeling approach to analyze the 2-D cross-sectional thermal field in steel bridges using finite-element models taking into account direct and diffuse solar radiation calculation. They conducted a sensitivity study based on the method and reported that the film coefficient (total surface thermal energy lost per unit area per unit temperature difference between surface and surrounding) of the outer surfaces and the absorptivity of the outer surfaces showed a significant effect on the thermal distribution.

\cite{bkaczkiewicz2018experimental} conducted an experimental and numerical study to study the temperature distribution of square hollow section joints from a fire protection perspective. They observed that the temperature distribution depended on the joint configuration and also varied at different sections of the joint based on geometrical properties and brace connection. 

2-D and 3-D spatial interpolation techniques are very commonly employed to estimate quantities at unsampled locations from discrete location measurements. Inverse distance weighting, kriging interpolation, trend surface analysis, natural neighbor interpolation, spline interpolation, etc are some methods for spatial interpolation.

\cite{zheng2019simulation} used kriging interpolation to estimate the 2-D thermal distribution of arc dams based on limited temperature monitoring points. They observed that the temperature fields obtained using different correlation models were quite similar.
\cite{zhou2019temperature} investigated concrete dam thermal field reconstruction employing the Distributed Temperature Sensing technology using optical fiber temperature measurements. They used the kriging difference algorithm to determine the 2-D temperature fields of longitudinal profiles of dam blocks from optical fiber temperature measurements. \cite{lin20213d} presented a 3-D thermal field reconstruction method for concrete dams using the kriging interpolation technique based on optical fiber temperature measurements. 
\cite{peng2020positioning} proposed a positioning method for temperature sensors to monitor the thermal field of dams based on the reconstructed thermal field using the nearest neighbor interpolation technique and cross-validation.

\cite{pan2022novel} presented a downstream surface thermal field reconstruction approach for arch dams based on convolutional neural networks capable of considering solar radiation, air temperature, and climatological data (rainfall, wind speed, cloud cover) as well. They observed that in addition to varying solar radiation causing non-uniform temperature distribution, the weather conditions also had a significant impact on the surface temperature distribution. They advocated for a high-frequency monitoring system that is able to capture the peak temperatures accurately and thus reconstruct more realistic temperature and stress fields, allowing for accurate detection of cracks and damage.

\cite{kulprapha2012structural} investigated the feasibility of ambient thermal load integrated Structural Health Monitoring for multi-span prestressed concrete bridges. Due to uneven thermal loading from the solar radiation, they noted the difference in temperatures and rate of heating and cooling of the different bridge surfaces depending on the time of the day, with maximum difference of around 20 $\degree C$ between the bridge deck and soffit in the early afternoon.
It has been observed in the literature that the temperature profile along the cross-section of bridges is non-linear when exposed to uneven ambient solar radiation (\cite{abid2022temperatures,fan2022efficient,kulprapha2012structural}).

In an experimental study, \cite{nguyen2016static} investigated the change in the deflection line of a statically loaded prestressed concrete bridge with varying levels of damage based on displacement measurements and proposed a temperature compensation algorithm based on the slope of the deflection-temperature curve.
\cite{sun2019bridge} considered temperature variation in their proposed real-time damage identification approach for bridges using inclination and strain measurements. They assumed that the temperature in the longitudinal and transverse direction was constant, while the temperature-induced inclinations and deflections at a point were only a function of time and not related to the material property, sectional property, or spatial temperature distribution. This allowed them to use even simple temperature gradients like the linear or bi-linear thermal distributions along the depth of the bridge.

\cite{abid2022temperatures} studied the temperature and temperature gradients in concrete bridges using experimental data to validate the finite-element model and conduct a parametric study and obtain temperature gradients along the depth of the concrete specimen in different regions of Turkey.
\cite{sheng2022prediction} proposed a prediction equation for concrete box-girder temperature gradients in different regions of China based on a validated finite element simulation model. They compared the long-term simulation results for the vertical temperature gradient with the results from empirical equations proposed by \cite{abid2016experimental}; \cite{lee2012analysis}; \cite{potgieter1983response}; and \cite{roberts2002measurements}.

\cite{fan2022efficient} proposed a dimensionally-reduced vertical discrete numerical model for rapid calculation of the temperature field and thermal load in steel-concrete composite bridges. They concluded that the calibrated one-dimensional numerical model was efficiently able to calculate the thermal field along the height of the bridge. However, simplifications regarding the lateral direction of the bridge limit the applicability of this method to certain bridge geometries.

Researchers have also investigated the thermal response of the structure rather than the actual thermal distribution when focusing on structural health monitoring (\acrshort{SHM}) (\cite{glashier2024iterative}; \cite{glashier2024temperature}; \cite{kromanis2014predicting}). Thermal response prediction is usually done with either model-based methods or data-driven methods (\cite{glashier2024iterative}; \cite{glashier2024temperature}; \cite{jin2016damage}; \cite{jin2015structural}; \cite{kromanis2014predicting}; \cite{peeters2001vibration}).

\cite{kromanis2014predicting} presented a generic regression-based thermal response prediction methodology to predict and isolate the thermal component from the distributed temperature and response measurements. They investigated the performance of four popular machine learning algorithms, namely multiple linear regression (MLR), robust regression (RR), artificial neural network (ANN), and support vector regression (SVR) based on two test bridges and concluded the regression-based models were able to accurately predict the thermal response of the bridges.   

\cite{glashier2024iterative, glashier2024temperature} analyzed the monitoring data obtained from the comprehensive network of sensors installed on the MX3D bridge (the world's first wire and arc additive manufactured metal structure) located in Amsterdam. Physics-based structural assessment of wire and arc additively manufactured structures is difficult due to geometric, material, and manufacturing variability; therefore, data-based techniques are preferred. \cite{glashier2024temperature} observed that the structural response of the bridge was mainly temperature-driven. They applied a data-driven iterative regression-based thermal response prediction method to use the predicted thermal response to correct and eliminate environmental and operational variations (temperature - induced (\cite{glashier2024iterative}) and pedestrian loading (\cite{glashier2024temperature})) component from response signals.

\cite{wang2023causal} studied the lag effect, i.e., the influence of change in air temperature and reservoir water level on the evolution of the dam thermal field. They mentioned that different regions of the dam are subjected to different lag influencing laws since the dominating factor causing the lag changes from the air-exposed region to the dam base. Based on this information, they divided the thermal field into the air temperature influencing zone, air temperature and reservoir water level joint influencing zone (mainly upstream surface), and trend changing zone (mainly dam heel and interior of dam body). To obtain the temperature field prediction of the dam, they trained a nonlinear autoregressive with exogenous inputs (NARX) neural network model for each temperature measuring point.

\cite{zhang2024temperature} presented a temperature analysis and prediction method for assessing the effective temperatures and temperature differences for road-rail steel truss bridges from a structural health monitoring point of view. They observed a strong direct correlation between the atmospheric temperature and the effective temperature of the main girder and a significant negative correlation between the atmospheric humidity and the temperature of the main girder. They also noticed a time lag (approximately 3 hours in winter and 5 hours in summer) between the solar radiation and the effective temperature of the main girder. They employed a Long-Short-Term Memory (LSTM) deep learning network to predict the temperatures of the girder based on the environmental conditions (air temperature and humidity, solar radiation, wind speed, and direction), factoring in also the past environmental conditions. They concluded that the LSTM networks gave better temperature predictions than linear regression models.

\cite{hou2022simulation} used the actual measured temperature data as feedback to invert and analyze the temperature evolution for the dam simulated using the finite-element method. A close agreement between the measured and simulated temperature indicated that appropriate \acrshort{FEM} analysis method and thermal parameters were chosen. These parameters were used as a basis to further simulate the thermal profile of the dam.

From a Structural Health Monitoring and structural damage assessment viewpoint, as evident from above, extensive research has been done to compensate for or eliminate the effect of ambient thermal loads. Response-temperature correlation (\cite{jin2015structural}), Principal Component Analysis based approach incorporating temperature as latent variables (\cite{nguyen2014damage, reynders2014output, wang2020damage}), machine learning based techniques (\cite{corbally2022data, gu2017damage, jin2016damage, kromanis2014predicting, pan2022novel, sharma2021bridge}), time-series analysis (\cite{kostic2017vibration, zhang2019eliminating}), and baseline temperature approach (\cite{soo2020new}) are some of the common methods for compensating temperature effects. These methods have proven to be effective in indicating the presence of damage, but to localize the damage, the consideration of a scalar temperature for the entire structure may lead to false results if the structure is unevenly heated.

Much of the research and design standards (like the \acrshort{AASHTO LRFD} 2020 (\cite{aashto2020lrfd})) to assess the thermal distribution of bridges focuses on the temperature distribution in the vertical direction. While this simplification reduces the fidelity of the problem and allows faster computations for thermal fields, it may be insufficient for use in applications such as Structural Health Monitoring and damage detection and localization where temperature distribution in the entire structure is desired with sufficient accuracy to prevent false alarms.

For instance, the conventional assumption of uniform longitudinal and transverse temperature distribution is found to be untrue (\cite{hagedorn2019temperature, ma2023statistical, zhou2013thermal}).
Additionally, the thermal loads included in design specifications are not an actual real-time reflection of the thermal profile of the structure (\cite{hagedorn2019temperature}). The actual thermal field of the structure changes throughout the day due to ambient temperature, solar radiation, wind speed, rainfall, humidity, etc. Accurately quantifying and taking into account all these factors to calculate thermal profile by analytical or numerical models, especially for several days and months, is cumbersome and prone to error, if not impossible. This leaves, using measurement data to reconstruct the thermal fields as the most practical solution. Direct temperature measurements or indirect calculations using displacement or strain measurements can be used to recover the temperature distribution in the structure.    

While spatial interpolation techniques based on discrete temperature monitoring data can satisfactorily estimate the thermal distribution within the structure, it is essential for these non-physics-based methods that the sensors detect the non-uniform temperatures within the structure to be able to reconstruct a realistic temperature distribution. 
In this work, an adjoint-based high-fidelity approach to obtain/recover the thermal field in structures using displacement and/or strain measurements is proposed. This physics-based approach, while also dependent on the sensor density and location, does not explicitly require the sensors to be placed in locations of non-uniform temperatures in the structure.  
This technique is similar to the approach used for detecting and localizing damage/weaknesses in structures, i.e., obtaining the material strength distribution in the structure (\cite{airaudo2023adjoint, FAiraudo_HAntil_RLoehner_URakhimov_2024a, lohner2024high, antonau2024formulation, RLoehner_Fairaudo_HAntil_RWuechner_FMeister_SWarnakulasuriya_2024b}).
To the knowledge of the authors, this approach has not been used before to reconstruct the structural thermal distribution in the context of high-fidelity Structural Health Monitoring using adjoint-based optimization. The proposed approach is able to identify the non-uniform thermal field even when none of the sensors are located in a different temperature zone as demonstrated in the Plate with a Hole example in Section \ref{sec:Plate With a Hole}. It is applicable as long as the sensors can detect any temperature-induced change in the structure. This approach is not specific to any particular type of civil engineering structure and can conveniently be applied to other areas of engineering as well, such as to structural wind, aerospace, and mechanical engineering.

\section{Methodology}
\label{sec:methodology}

This section describes the methodology used to determine the thermal field of a structure based on displacement and/or strain measurements at limited sensor locations. Section \ref{sec:thermal_stress} explains the consideration of thermal stresses in a mechanical structural FE analysis. Section \ref{sec:thermal_field} expands on the proposed methodology of obtaining the temperature distribution from the displacement/strain measurements. Section \ref{sec:optimization} describes the optimization procedure with a brief overview of the Vertex Morphing method for regularization.

\subsection{Thermal Stresses and Internal Forces}
\label{sec:thermal_stress}

For a one-dimensional rod of length $L$ with thermal expansion coefficient $\alpha$, an increment of temperature $\Delta T$ results in an expansion in length of the rod given by:
\begin{equation}
 \Delta L = \alpha \cdot L \cdot \Delta T \quad.
\end{equation}
This implies that the thermal strain $\epsilon$ is given by:
\begin{equation}
    \epsilon = \alpha \cdot \Delta T \quad.
\end{equation}
Assuming an isotropic thermal strain (as is most often the case), the stress tensor is
given by:
\begin{equation}
    \sigma^{ij} = -\frac{\alpha E}{1-2\nu} \cdot \Delta T \cdot \delta^{ij} \quad,
\label{eq:thermal stress}
\end{equation}
where $E$, $\nu$ denote the Young’s modulus and Poisson coefficient. If we now take the conservation of momentum (forces) equations, given by:
\begin{equation}
    \frac{\partial\sigma^{ij}}{\partial x^{j}} = f \quad,
\end{equation}
and only consider the internal forces $f_{int}^{k}$ due to thermal strain, we have:
\begin{equation}
f^{k}_{int} = - \frac{\partial}{\partial x^{k}}\left(   \frac{\alpha E}{1-2\nu} \Delta T    \right) \quad,
\end{equation}
Within a finite element discretization (and after integration by parts) this translates into
\begin{equation}
    \int N^{i} \cdot f^{k}_{int} d\Omega = \int \frac{\partial N^{i}}{\partial x^{k}} \cdot  \frac{\alpha E}{1-2\nu} \cdot  N^{j} \cdot  \Delta \hat{T}_{j} d\Omega \quad,
\end{equation}
where $N^i$ denotes the shape function of node $i$ and $\Delta \hat{T}_{j}$ is the temperature change of the node $j$. For linear elements, this results in:
\begin{equation}
\hat{f}^{i}_{k} = \sum_{els} \frac{\partial N^{i}}{\partial x^{k}} \cdot \frac{\alpha E}{1-2\nu}\cdot  N^{j} \cdot Vol_{el} \cdot \Delta T^{avg}_{el} \quad, \quad \Delta T^{avg}_{el} = \frac{1}{nn} \sum_{j=1}^{nn} \Delta \hat{T}_j \quad,
\end{equation}
where $\Delta T^{avg}_{el}$ is the average temperature change of the element and $nn$ is the number of nodes associated with that element.

\subsection{Determining the Thermal Field from Measurements}
\label{sec:thermal_field}

Let $\Delta T(\mathbf{x}) $  be the spatial temperature distribution throughout the structure. The $\Delta$ refers to the temperature difference compared to the ambient; therefore, the temperature difference is used instead of absolute temperatures.
For a discretized system, determining the thermal field $\Delta T(\mathbf{x}) $ from measurements may be formulated as an optimization problem as follows: Given $n$ test cases with each recording measurements at $m$ measuring points/locations $\mathbf{x}_j$, $j = 1, m$ of their respective deformations $\mathbf{u}_{ij}^{md}$, $i = 1, n$, $j = 1, m$ or strains $\mathbf{s}_{ij}^{ms}$, $i = 1, n$, $j = 1, m$, obtain the spatial distribution of the temperature $\Delta \mathbf{T}$ that minimizes the cost function:
\begin{equation}
I(\mathbf{u}, \Delta\mathbf{T}) = \frac{1}{2} \sum_{i=1}^{n} \sum_{j=1}^{m} w_{ij}^{md} (\mathbf{u}_{ij}^{md} - \mathbf{I}_{ij}^{d} \cdot \mathbf{u}_i)^ 2 + \frac{1}{2} \sum_{i=1}^{n} \sum_{j=1}^{m} w_{ij}^{ms} (\mathbf{u}_{ij}^{ms} - \mathbf{I}_{ij}^{s} \cdot \mathbf{s}_i)^ 2 \quad,
\label{eq:cost_function}
\end{equation}
with the finite element description given by: 
\begin{subequations}\label{eq:fe_equations}
    \begin{align}
        \mathbf{K} \cdot \mathbf{u} & = \mathbf{f}_{ext} + \mathbf{f}_{\Delta\mathbf{T}} \quad,
        \label{eq:fe_equation_a} \\
        \mathbf{K}_{ij} & = \sum_{e=1}^{N_e} \mathbf{K}_{ij}^{e} \quad ,
        \label{eq:fe_equation_b}
    \end{align}
\end{subequations}
where $w_{ij}^{md}, w_{ij}^{ms}$ are displacement and strain weights, $\mathbf{I}^{d} , \mathbf{I}^{s}$ are the interpolation matrices used to obtain the displacements and strains from the finite element mesh at the measurement locations, $\mathbf{K}$ is the usual stiffness matrix, and $\mathbf{f}_{ext}, \, \mathbf{f}_{\Delta \mathbf{T}}$ are the external forces and the forces due to thermal strains. We emphasize that the optimization problem given by Eqns.(\ref{eq:cost_function}-\ref{eq:fe_equations}) does not assume any specific choice of finite element basis functions, and is therefore widely applicable.

Eqn.\eqref{eq:fe_equation_a} can be reformulated to obtain the residual formulation ($\mathbf{R} = \mathbf{0}$) and the change in residual as:
\begin{subequations}\label{eq:residual_equations}
    \begin{gather}
        \mathbf{R} =  \mathbf{f}_{ext} + \mathbf{f}_{\Delta\mathbf{T}} - \mathbf{K} \cdot \mathbf{u} = \mathbf{0}\quad, \label{eq:residual_equations_a} \\
        \frac{d\mathbf{R}}{d\Delta\mathbf{T}} = \frac{\partial \mathbf{R}}{\partial \Delta\mathbf{T}} + \frac{\partial \mathbf{R}}{\partial \mathbf{u}} \cdot \frac{d \mathbf{u}}{d \Delta \mathbf{T}} = \mathbf{0} \quad \implies \frac{d \mathbf{u}}{d \Delta \mathbf{T}} = - \left[ \frac{\partial \mathbf{R}}{\partial \mathbf{u}} \right]^{-1} \cdot \frac{\partial \mathbf{R}}{\partial \Delta\mathbf{T}}      \quad       . \label{eq:residual_equations_b} 
    \end{gather}
\end{subequations}
The objective function can be extended to the Lagrangian functional (cf.~\cite{HAntil_DPKouri_MDLacasse_DRidzal_2018a})
\begin{equation}
    L(\mathbf{u},\Delta\mathbf{T}, \mathbf{\tilde{u}}^t) = I(\mathbf{u}, \Delta\mathbf{T}) + \mathbf{\tilde{u}}^t \cdot \mathbf{R} = I(\mathbf{u}, \Delta\mathbf{T}) + \mathbf{\tilde{u}}^t \cdot ( \mathbf{f}_{ext} + \mathbf{f}_{\Delta\mathbf{T}} - \mathbf{K} \cdot \mathbf{u} ) \quad,
\label{eq:lagrangian}
\end{equation}
where $\mathbf{\tilde{u}}$ are the Lagrange multipliers (adjoints) and the superscript $t$ refers to the transpose. The variation of the Lagrangian with respect to each of its unknowns then results in:
\begin{subequations}\label{eq:lagrangian_derivative}
    \begin{gather}
        \frac{\partial L}{\partial \Delta \mathbf{T}} = \frac{\partial I}{\partial \Delta \mathbf{T}} + \mathbf{\tilde{u}}^t \cdot \frac{\partial \mathbf{R}}{\partial \Delta \mathbf{T}}  = \mathbf{0} \quad , \label{eq:lagrangian_derivative_a} \\
        \frac{\partial L}{\partial\mathbf{\tilde{u}}} = \mathbf{R}=  \mathbf{f}_{ext} + \mathbf{f}_{\Delta\mathbf{T}} - \mathbf{K} \cdot \mathbf{u} = \mathbf{0} \quad, \label{eq:lagrangian_derivative_b} \\
        \frac{\partial L}{\partial \mathbf{u}} = \frac{\partial I}{\partial \mathbf{u}} + \mathbf{\tilde{u}}^t \cdot \frac{\partial \mathbf{R}}{\partial \mathbf{u}} = \mathbf{0} \quad, 
        \label{eq:lagrangian_derivative_c}
    \end{gather}
\end{subequations}
where $\mathbf{J}_{ij}^{s}$ denotes the relationship between the displacements and strains (i.e. the derivatives of the displacement field on the finite element mesh and the location $\mathbf{x}_j$). Eqn.\eqref{eq:lagrangian_derivative_b} is the usual forward problem.
Eqn.\eqref{eq:lagrangian_derivative_c} can be rewritten to represent the adjoint problem as:
\begin{equation}
    \left[- \frac{\partial \mathbf{R}}{\partial \mathbf{u}} \right]^{t} \cdot \mathbf{\tilde{u}} = \frac{\partial I}{\partial \mathbf{u}} \implies \mathbf{K}^{t} \cdot \mathbf{\tilde{u}} = \frac{\partial I}{\partial \mathbf{u}} \quad,
\label{eq:adjoint_problem}
\end{equation}
where $\mathbf{K}^{t} = \mathbf{K}$ for structural systems, thereby allowing using the same stiffness matrix as the forward problem, for the adjoint solve. The right-hand side component in Eqn.\eqref{eq:adjoint_problem} is understood as the pseudo-force when comparing it to the right-hand side of the forward problem ($\mathbf{K}\cdot\mathbf{u}= \mathbf{f}_{ext}+\mathbf{f}_{\Delta\mathbf{T}}$).

The important part of the methodology comes in the gradient calculation. The gradient of the Lagrangian w.r.t. the temperature distribution $\Delta\mathbf{T}$ can be written as: 
\begin{subequations}\label{eq:gradient_equations}
    \begin{gather}
        \frac{dL}{d\Delta\mathbf{T}} = \frac{dI}{d\Delta\mathbf{T}}  +   \mathbf{\tilde{u}}^t \cdot \frac{d\mathbf{R}}{d\Delta\mathbf{T}}  \quad, \label{eq:gradient_equations_a} \\
        \intertext{Inserting $\frac{d\mathbf{R}}{d\Delta\mathbf{T}}$ from Eqn.\eqref{eq:residual_equations_b}:}
        \frac{dL}{d\Delta\mathbf{T}} = \frac{dI}{d\Delta\mathbf{T}} = \frac{\partial I}{\partial \Delta\mathbf{T}} + \frac{\partial I}{\partial \mathbf{u}} \cdot \frac{d \mathbf{u}}{d \Delta \mathbf{T}}  \quad,         \label{eq:gradient_equations_b} \\
        \intertext{Inserting $\frac{d \mathbf{u}}{d \Delta \mathbf{T}}$ from Eqn.\eqref{eq:residual_equations_b}:}
        \frac{dI}{d\Delta\mathbf{T}} = \frac{\partial I}{\partial \Delta\mathbf{T}} + \frac{\partial I}{\partial \mathbf{u}} \cdot \left(  - \left[ \frac{\partial \mathbf{R}}{\partial \mathbf{u}} \right]^{-1} \cdot \frac{\partial \mathbf{R}}{\partial \Delta\mathbf{T}}  \right) \quad,        \label{eq:gradient_equations_c} \\
        \intertext{Rearranging the brackets and inserting $\mathbf{\tilde{u}}^t$ from Eqn.\eqref{eq:adjoint_problem}:}
        \frac{dI}{d\Delta\mathbf{T}} =  \frac{\partial I}{\partial \Delta\mathbf{T}} +  \left( -\frac{\partial I}{\partial \mathbf{u}} \cdot \left[ \frac{\partial \mathbf{R}}{\partial \mathbf{u}} \right]^{-1}  \right)   \cdot \frac{\partial \mathbf{R}}{\partial \Delta\mathbf{T}}  = \frac{\partial I}{\partial \Delta\mathbf{T}} +  \mathbf{\tilde{u}}^t \cdot \frac{\partial \mathbf{R}}{\partial\Delta\mathbf{T}}  \quad.        \label{eq:gradient_equations_d} 
    \end{gather}
\end{subequations}
In this work, the temperature difference between different regions of the structure is considered to be small, such that the Young's modulus of the material is assumed to be constant. Of course, for large temperature changes, this assumption would not hold true, and Young's modulus-temperature variation needs to be included. Considering this simplified assumption, Eqns. \eqref{eq:residual_equations_a} and \eqref{eq:gradient_equations_d} can be combined to obtain the gradient used in this work as:
\begin{equation}
    \frac{dI}{d\Delta\mathbf{T}} = \nabla I (\Delta \mathbf{T}) =\frac{\partial I}{\partial \Delta\mathbf{T}} +  \mathbf{\tilde{u}}^t \cdot \frac{\partial\mathbf{R}}{\partial\Delta\mathbf{T}} = \frac{\partial I}{\partial \Delta\mathbf{T}} +  \mathbf{\tilde{u}}^t \cdot \frac{\partial\mathbf{f}_{\Delta\mathbf{T}}}{\partial\Delta\mathbf{T}}   \quad.        \label{eq:gradient_final} \\
\end{equation}

\medskip

\subsection{Optimization}
\label{sec:optimization}

An optimization cycle using the adjoint approach is composed of the following steps:
\begin{itemize}
    \item With the current temperature distribution $\Delta \mathbf{T}$: solve the $n$ forward problems (Eqn.\eqref{eq:lagrangian_derivative_b}) corresponding to the $n$ load cases to obtain the displacement fields $\xrightarrow{} \mathbf{u}= \left[ \mathbf{u}_1, \dots, \mathbf{u}_i, \dots, \mathbf{u}_n  \right]$,
    \item With the $n$ current displacement fields $\mathbf{u}$ and sensor measurements $\mathbf{u}_{ij}^{md}, \mathbf{s}_{ij}^{ms}$: solve adjoint problem (Eqn.\eqref{eq:adjoint_problem}) corresponding to each load cases to obtain $n$ sets of adjoint variables $\xrightarrow{} \tilde{\mathbf{u}}= \left[ \tilde{\mathbf{u}}_1, \dots, \tilde{\mathbf{u}}_i, \dots, \tilde{\mathbf{u}}_n  \right] $,
    \item With the $n$ sets of adjoint variables $\tilde{\mathbf{u}}$: calculate the gradients (Eqn.\eqref{eq:gradient_final}) for each load case $\xrightarrow{} \frac{d I}{d\Delta\mathbf{T}} = \left[ \frac{d I^{1}}{d\Delta\mathbf{T}}, \dots, \frac{d I^{i}}{d\Delta\mathbf{T}}, \dots, \frac{d I^{n}}{d\Delta\mathbf{T}}           \right] $ \quad,
    \item Combine and if necessary smoothen the gradients $\xrightarrow{}  \frac{d I^{smooth}}{d\Delta\mathbf{T}} $   \quad,
    \item Update the temperature distribution $\xrightarrow{} \Delta \mathbf{T}_{new} = \Delta \mathbf{T}_{old} - \gamma \cdot \frac{d I^{smooth}}{d\Delta\mathbf{T}} $\quad.
\end{itemize}

Here $\gamma$ is a small stepsize that can be adjusted so as to obtain optimal convergence (e.g. via a line search). In this study, the Barzilai-Borwein (\acrshort{BB}) method for stepsize calculation is used. The stepsize $\gamma$ at iteration $i$ is calculated as:

\begin{equation}
\gamma^{(i)} = \frac{\mathbf{d}^{(i-1),t} \cdot \mathbf{d}^{(i-1)}}{\mathbf{d}^{(i-1),t} \cdot \mathbf{y}^{(i)}} \quad ,
\label{eq:bb}
\end{equation}
where the superscript $(i)$ and $(i-1)$ refer to the iteration count, $\mathbf{y}^{(i)} = \nabla I(\Delta\mathbf{T}^{(i)}) - \nabla I(\Delta\mathbf{T}^{(i-1)})$ is the change of cost function gradients from the previous iteration, and $\mathbf{d}^{(i-1)} = \Delta\mathbf{T}^{(i)} - \Delta\mathbf{T}^{(i-1)} $ is the design update of the previous iteration. Detailed explanation about the method can be found in \cite{antonau2023enhanced, fletcher2005barzilai}.

\subsubsection{Vertex Morphing}
\label{sec:vertex_morphing}

Vertex Morphing is a shape control method commonly used for node-based shape optimization (\cite{hojjat2014vertex}). The approach aims at identifying the design field in the physical space which is controlled by a field in the non-physical control space through a mapping function. 

However, before this physical to control space mapping, the physical design variables (in this case, nodal temperature differences ($\mathbf{\Delta T}$) are projected to a range of $[0,1]$, called as \textit{physical phi field} ($\mathbf{\Phi}$), using the Sigmoid function. This projection has the advantage of smoothly approaching zero gradient at the design bounds, rather than strictly cutting off the gradients when they exceed the bounds.

Vertex Morphing is further applied on this physical $\mathbf{\Phi}$ field. The discretized formulation is written in Eqn.\eqref{eq:vm_forward_filter}, where $\mathbf{\Phi}$ is the vector of sigmoid-projected nodal design variables in the physical space, $\mathbf{\tilde{\Phi}}$ is the vector of corresponding discretized control variables in control space, and $\mathbf{A}$ is the kernel function/operator matrix facilitating the mapping of  $\mathbf{\tilde{\Phi}}$  to $\mathbf{\Phi}$. This step is also known as \textit{forward filtering}. It is considered here that the control field and the physical field have the same discretization, therefore the size of the vectors $\mathbf{\Phi}$ and $\mathbf{\tilde{\Phi}}$ are the same.
\begin{equation}
    \mathbf{\Phi} = \mathbf{A} \cdot \mathbf{\tilde{\Phi}} \qquad or \qquad \Phi_p = A_{pq} \cdot \tilde{\Phi}_q \qquad ,
\label{eq:vm_forward_filter}
\end{equation}
where $A_{pq}$ refers to the filter effect corresponding to interaction between node $p$ and node $q$ based on their positions and Euclidean distance.

Since the algorithm works solely in the control space, the gradient of the cost function w.r.t. the physical phi variables i.e. physical gradients ($\frac{dI}{d\mathbf{\Phi}}$) need to be mapped to the control space gradients ($\frac{dI}{d\mathbf{\tilde{\Phi}}}$) using the chain rule of differentiation:
\begin{subequations}\label{eq:vm_backward_filter} 
    \begin{gather}
         \frac{dI}{d\tilde{\Phi}_p} = \frac{d I}{d \Phi_q} \cdot \frac{d\Phi_q}{d\tilde{\Phi}_p} \quad, 
         \label{eq:vm_backward_filter_a}\\
         \implies \frac{dI}{d\tilde{\Phi}_p} = A_{qp} \cdot \frac{d I}{ d \Phi_q} = A_{qp} \cdot b_{q} \qquad or \qquad  \frac{dI}{d\mathbf{\tilde{\Phi}}} = \mathbf{A}^{t} \cdot \mathbf{b} = \mathbf{A}^{t} \cdot \frac{dI}{d\mathbf{\Phi}} \quad,
         \label{eq:vm_backward_filter_b}
    \end{gather} 
\end{subequations}
where $b_q = \frac{dI}{d\Phi_q}$ is the gradient of the cost function w.r.t. the $q$-th design variable in the physical field. This step is also known as \textit{backward filtering}. The operator matrix $\mathbf{A}$ (which is symmetric for most kernels) is used twice; firstly, for the forward filter which can be seen as a way to suppress high oscillations when generating physical space design using the control field; and secondly, for smoothing the raw oscillating physical gradients to the control gradients.  

Since the shape optimization or the current thermal field identification optimization problems usually have a larger number of design variables than the responses, the design variables behave independently and vary in a haphazard manner, resulting in a rough surface or temperature field with local and sharply changing temperatures. To remedy this and obtain smooth design updates, Vertex Morphing works as a smoothing operator that performs a convolution in a specified radius, thereby resulting in smoother results. An in-depth explanation of Vertex Morphing can be found in \cite{antonau2022latest, ghantasala2021realization, hojjat2014vertex}. Other methods, such as simple point averaging, Laplacian smoothing, Pseudo-Laplacian smoothing, etc, can also be used alternatively to achieve gradient smoothing (\cite{airaudo2023adjoint, lohner2024high}).

\section{Examples}
\label{sec:examples}

The applicability of the proposed approach has been demonstrated using a Plate With a Hole example modeled with triangular shell elements, a Bridge example modeled with linear trusses, and a (simplified) Hoover dam example modeled with tetrahedral solid elements. All the numerical examples were simulated and analyzed using Kratos Multiphysics (\cite{dadvand2010object,vicente_mataix_ferrandiz_2024_6926179}), a general, open-source finite element framework for multiphysics applications supporting many element types, material models, and other options.

For all the examples, a `target' temperature distribution of $\Delta \mathbf{T}(\mathbf{x})$ was given, together with the external loads $\mathbf{f}_{ext}$. The system was then simulated and the displacements $\mathbf{u}(\mathbf{x})$ and strains $\mathbf{s}(\mathbf{x})$ were obtained and recorded at the 'measurement locations' $\mathbf{x}_j$, $j=1,m$. This then yielded the 'measurement pair' $\mathbf{f},\mathbf{u}_j$, $j=1,m$ or $\mathbf{f},\mathbf{s}_j$, $j=1,m$ that was used to determine the spatial temperature distribution $\Delta \mathbf{T} (\mathbf{x})$ throughout the structure.

\subsection{Plate With a Hole}
\label{sec:Plate With a Hole}

\begin{figure}[!t]
\begin{minipage}[c][][t]{.5\textwidth}
  \vspace*{\fill}
  \centering
   \includegraphics[trim= 0 0 0 0, clip, width=0.3\paperwidth]{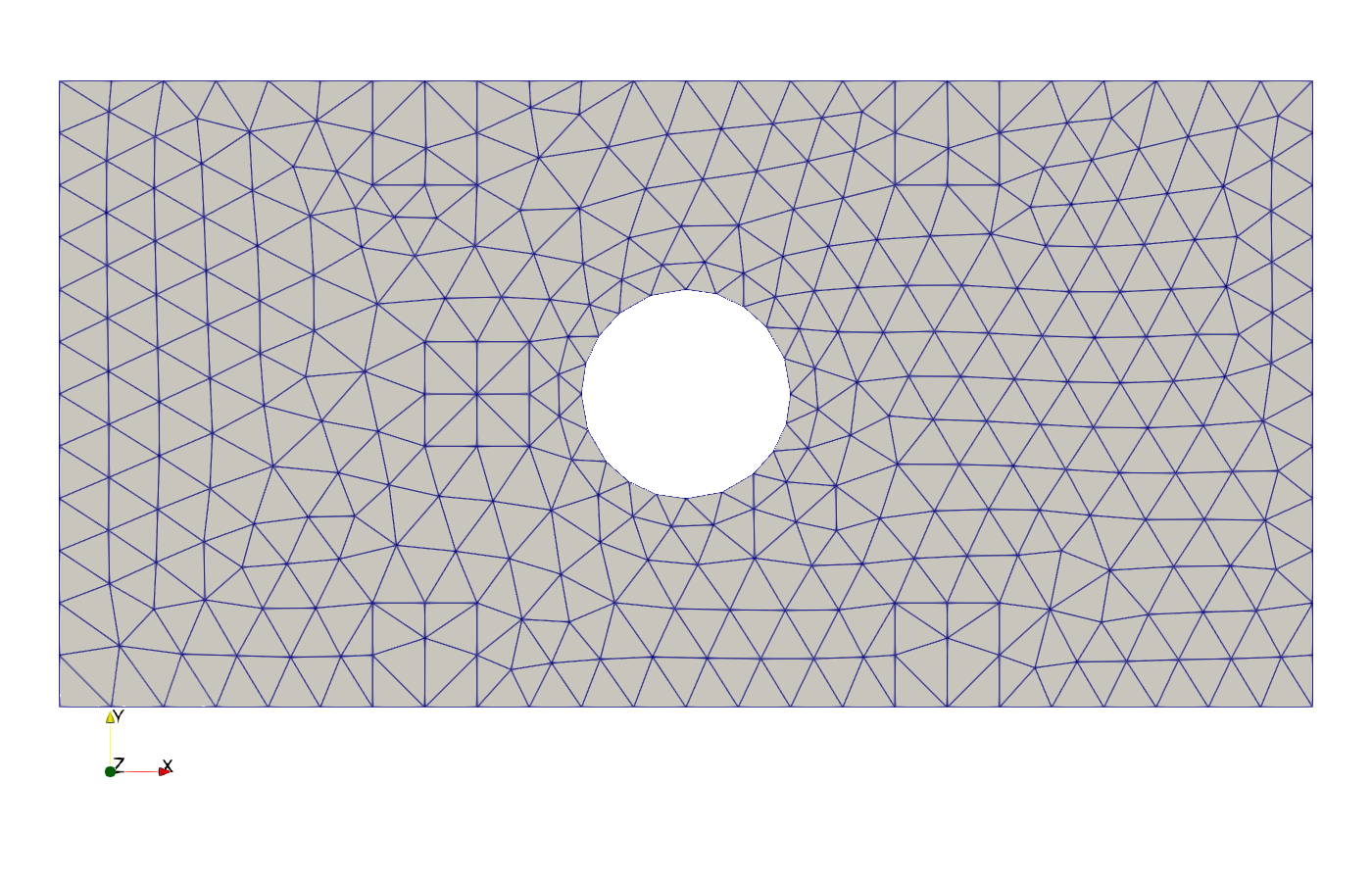}
\end{minipage}
\begin{minipage}[c][][t]{.5\textwidth}
  \vspace*{\fill}
  \centering
  \includegraphics[trim= 0 0 0 0, clip, width=0.3\paperwidth]{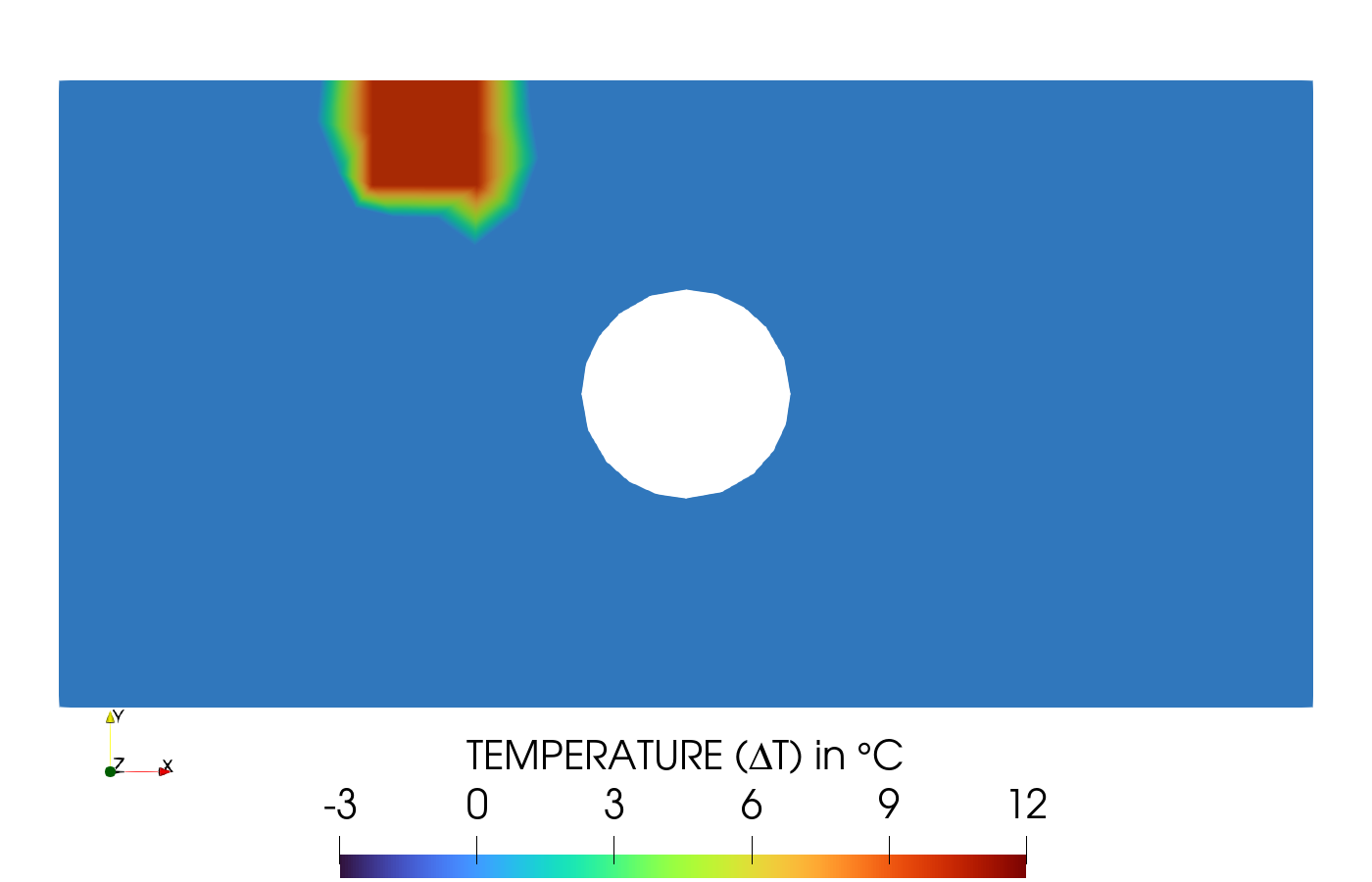}
\end{minipage}
  \caption{Plate With a Hole: Mesh used (left) and the target temperature distribution (right).}
  \label{fig:Figure_3}
\end{figure}

\begin{figure}[!b]
\begin{minipage}[c][][t]{.5\textwidth}
  \vspace*{\fill}
  \centering
   \includegraphics[trim= 0 0 0 0, clip, width=0.3\paperwidth]{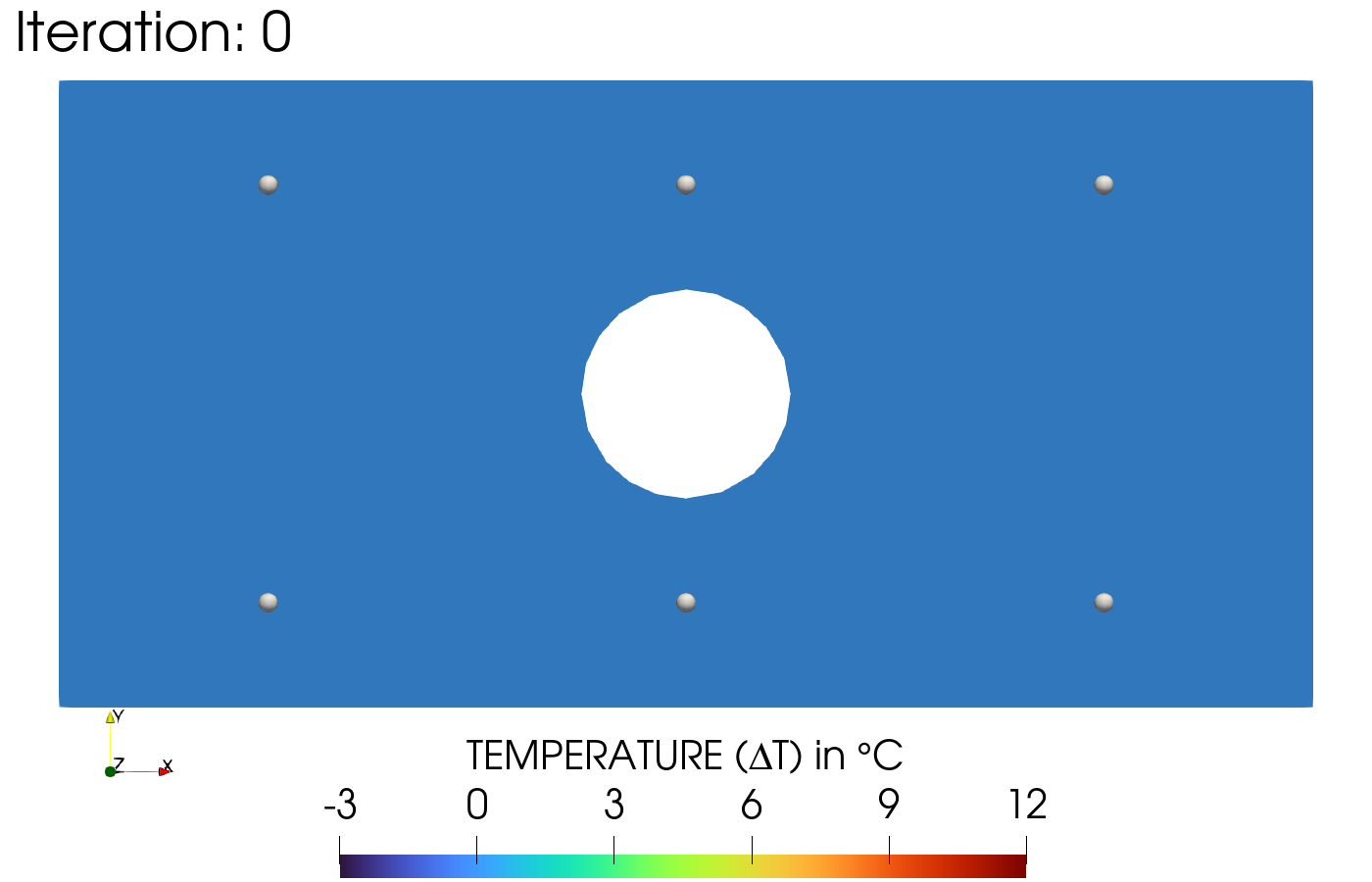}
\end{minipage}
\begin{minipage}[c][][t]{.5\textwidth}
  \vspace*{\fill}
  \centering
  \includegraphics[trim= 0 0 0 0, clip, width=0.3\paperwidth]{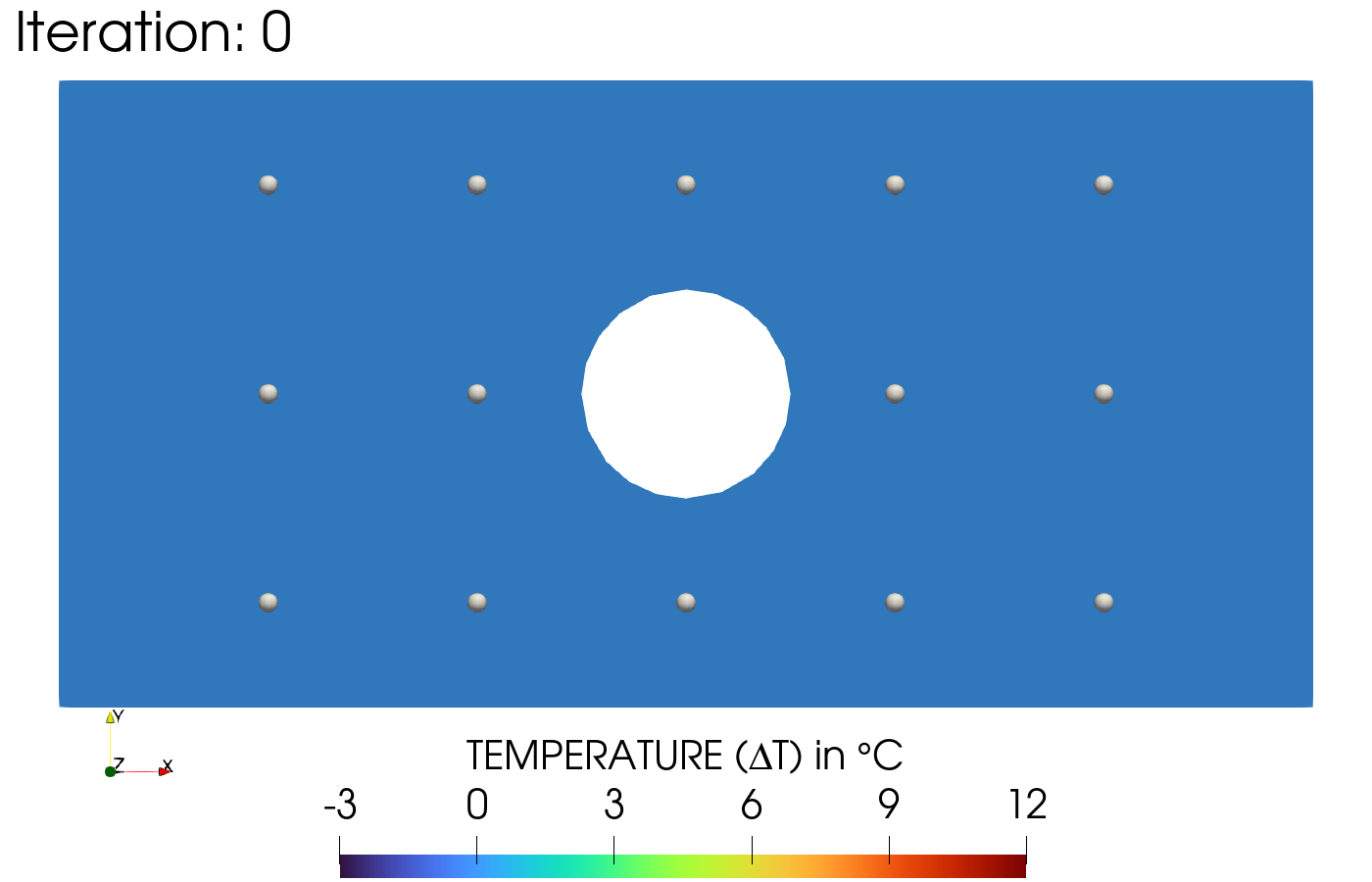}
\end{minipage}
  \caption{Plate With a Hole: 6 (left) and 14 (right) sensors configurations and the temperature distribution at optimization start $\Delta \mathbf{T} = \mathbf{0}$.}
  \label{fig:Figure_4}
\end{figure}

The case is shown in Figures \ref{fig:Figure_3} and considers a plate with a hole. The dimensions are (all units in SI): $0 \leq x \leq 60$, $0 \leq y \leq 30$, $0 \leq z \leq 0.1$. A hole of diameter $d = 10$ is placed in the middle ($x = 30$, $y = 15$). Density, Young's modulus, Poisson ratio, and thermal expansion coefficient were set to $\rho = 7800$, $E = 2\cdot 10^{12}$, $\nu = 0.3$, $\alpha = 1.0 \cdot 10^{-5} /\degree K$ respectively. 646 linear, triangular, plane stress elements were used. The mesh is displayed in Figure \ref{fig:Figure_3} (left). The left boundary of the plate was assumed clamped ($\mathbf{u} = \mathbf{0}$), while a horizontal load of $q = (10^{5}, 0, 0)$ was prescribed at the right end. The temperature difference was prescribed to be $\Delta T  = 0$ $\degree C$ everywhere except in the region $15 \leq x \leq 20$, $26 \leq y \leq 30$, where a temperature difference of $\Delta T = 10$ $\degree C$ was set. This distribution is shown in Figure \ref{fig:Figure_3}.

Two sensor configurations consisting of 6 and 14 displacement sensors distributed over the plate were analyzed. Figure \ref{fig:Figure_4} shows the sensor locations and the temperature distribution at the start of the optimization ($\Delta T = 0$ $\degree C$ everywhere).

For both the 6 and 14 sensors cases, the steepest descent optimization algorithm with the Barzilai-Borwein stepsize method (Eqn.\eqref{eq:bb}) was used. A maximum stepsize of $2.5\cdot 10^{-3}$ was also set to prevent huge updates in each iteration. The optimization starts with the presumption that the $\Delta T = 0$ $\degree C$ everywhere on the plate.
Stopping criteria was set to 5 magnitudes reduction of the initial cost function, i.e., the algorithm stops when the $\text{current cost function} \leq (1\cdot 10^{-5} * \text{initial cost function})$. 
Vertex morphing with radius $r = 5$ was used to smoothen the gradients and subsequently get a smoother field. The radius was chosen according to the rule of thumb such that the radius is roughly equal to 2-3 element side lengths.
For better visualization of the results, the color bar is shown in the range $[-3,12]$ $\degree C$ as the temperature ($\Delta T$) in no case falls below $-3$ $\degree C$, neither does it exceed $12$ $\degree C$.  

\medskip

\begin{figure}[!t]
\begin{minipage}[c][][t]{.5\textwidth}
  \vspace*{\fill}
  \centering
   \includegraphics[trim= 0 0 0 0, clip, width=0.3\paperwidth]{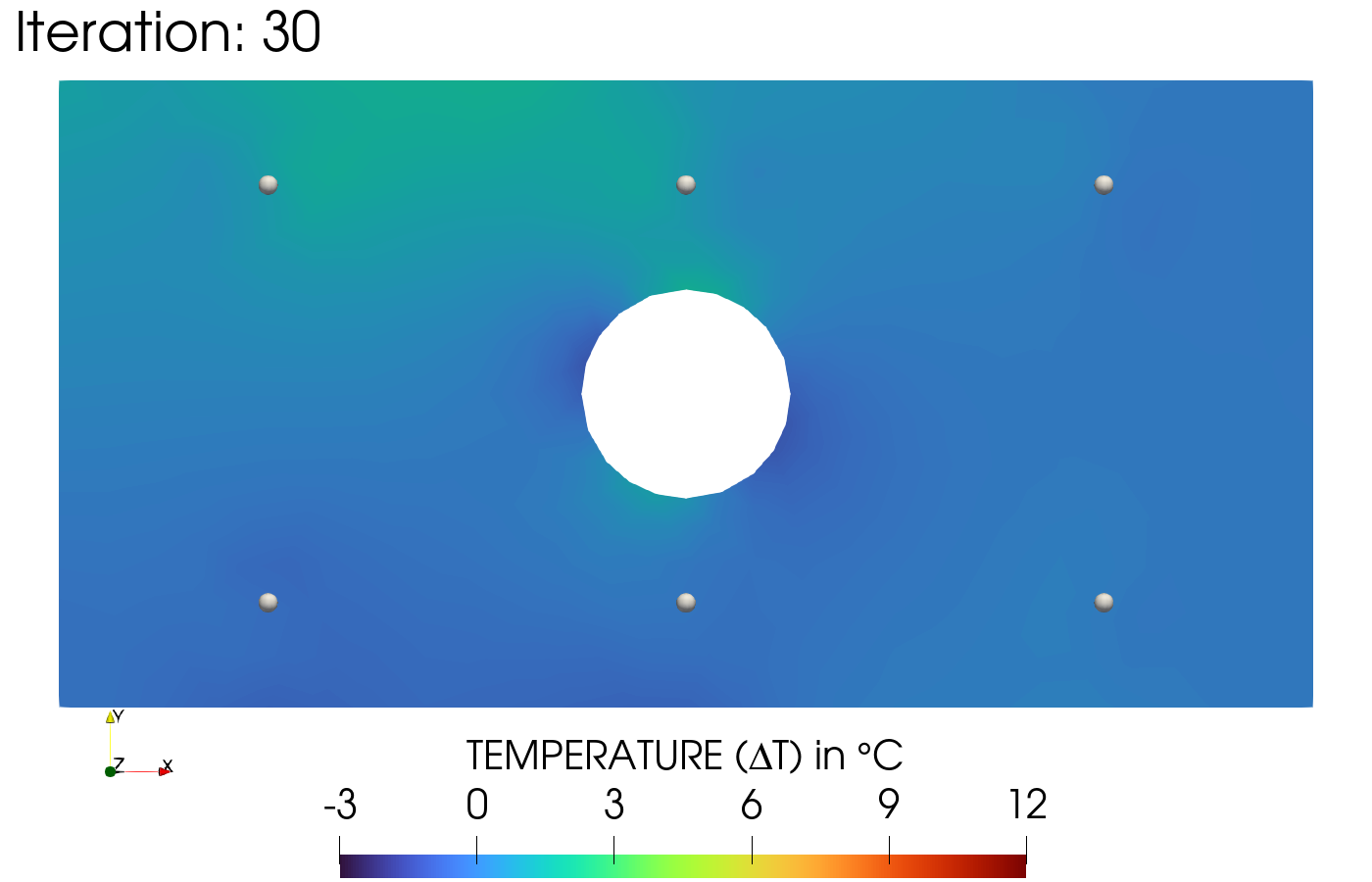}
  \subcaption{Without VM: Iteration 30}
   \label{fig:Figure_5a}
   \vspace{0.5em}
  \includegraphics[trim= 0 0 0 0, clip, width=0.3\paperwidth]{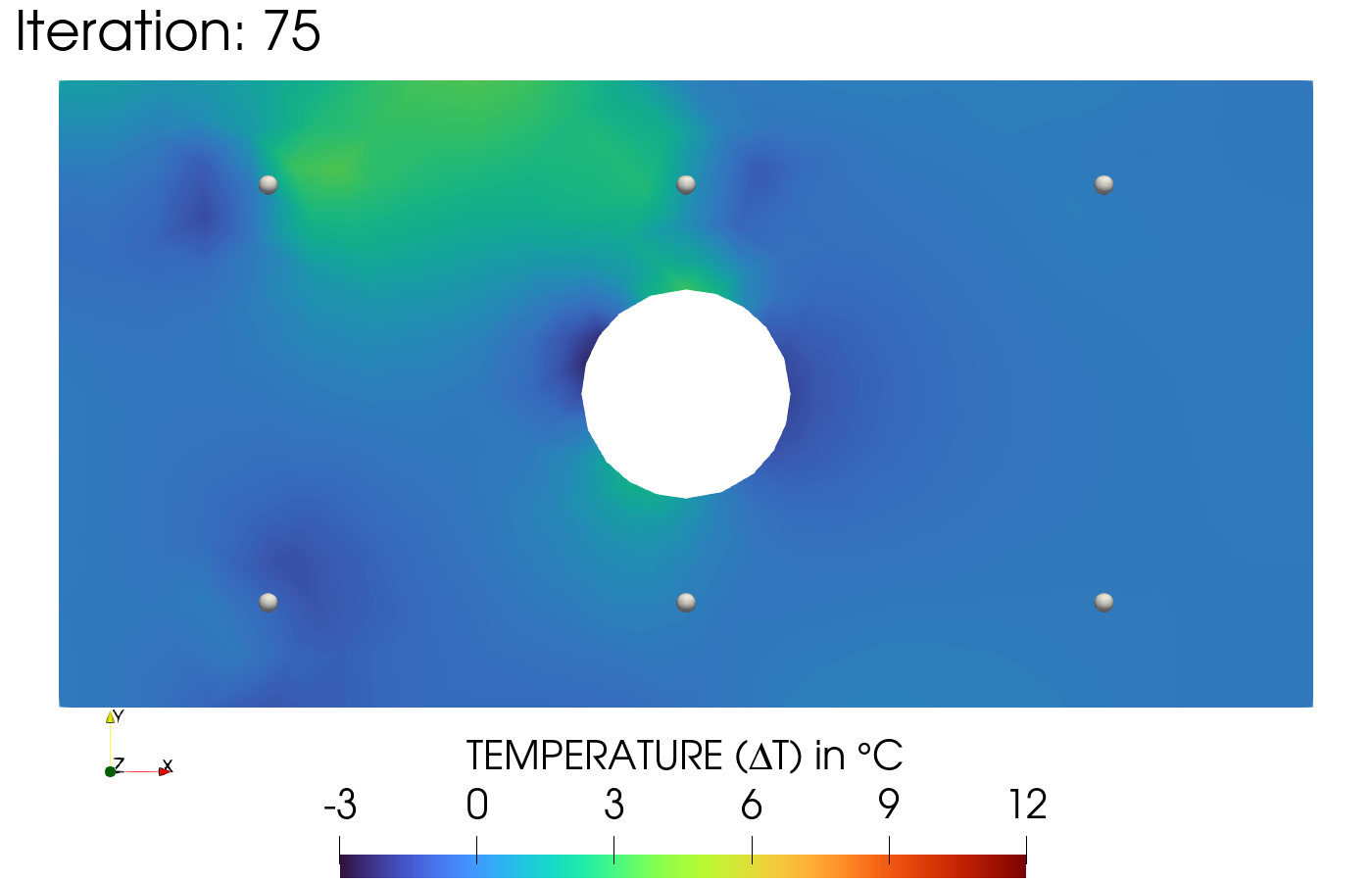}
  \subcaption{Without VM: Iteration 75}
   \label{fig:Figure_5b}
   \vspace{0.5em}
  \includegraphics[trim= 0 0 0 0, clip, width=0.3\paperwidth]{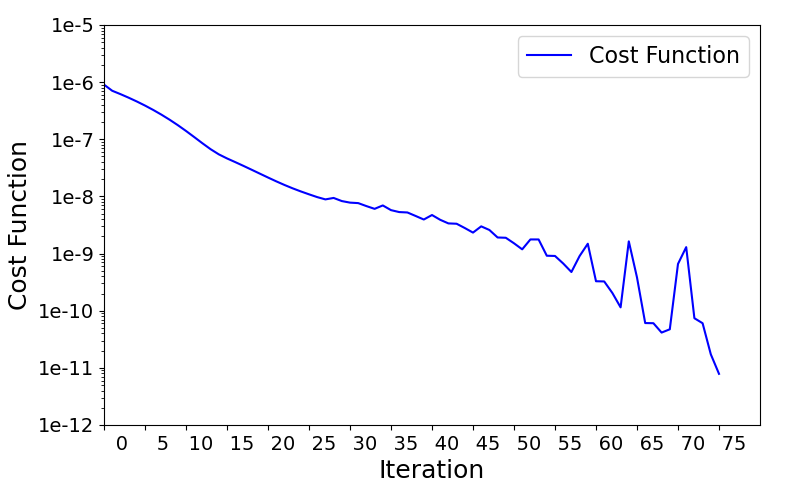}
  \subcaption{Without VM: Cost function convergence}
  \label{fig:Figure_5c}\par
\end{minipage}
\begin{minipage}[c][][t]{.5\textwidth}
  \vspace*{\fill}
  \centering
  \includegraphics[trim= 0 0 0 0, clip, width=0.3\paperwidth]{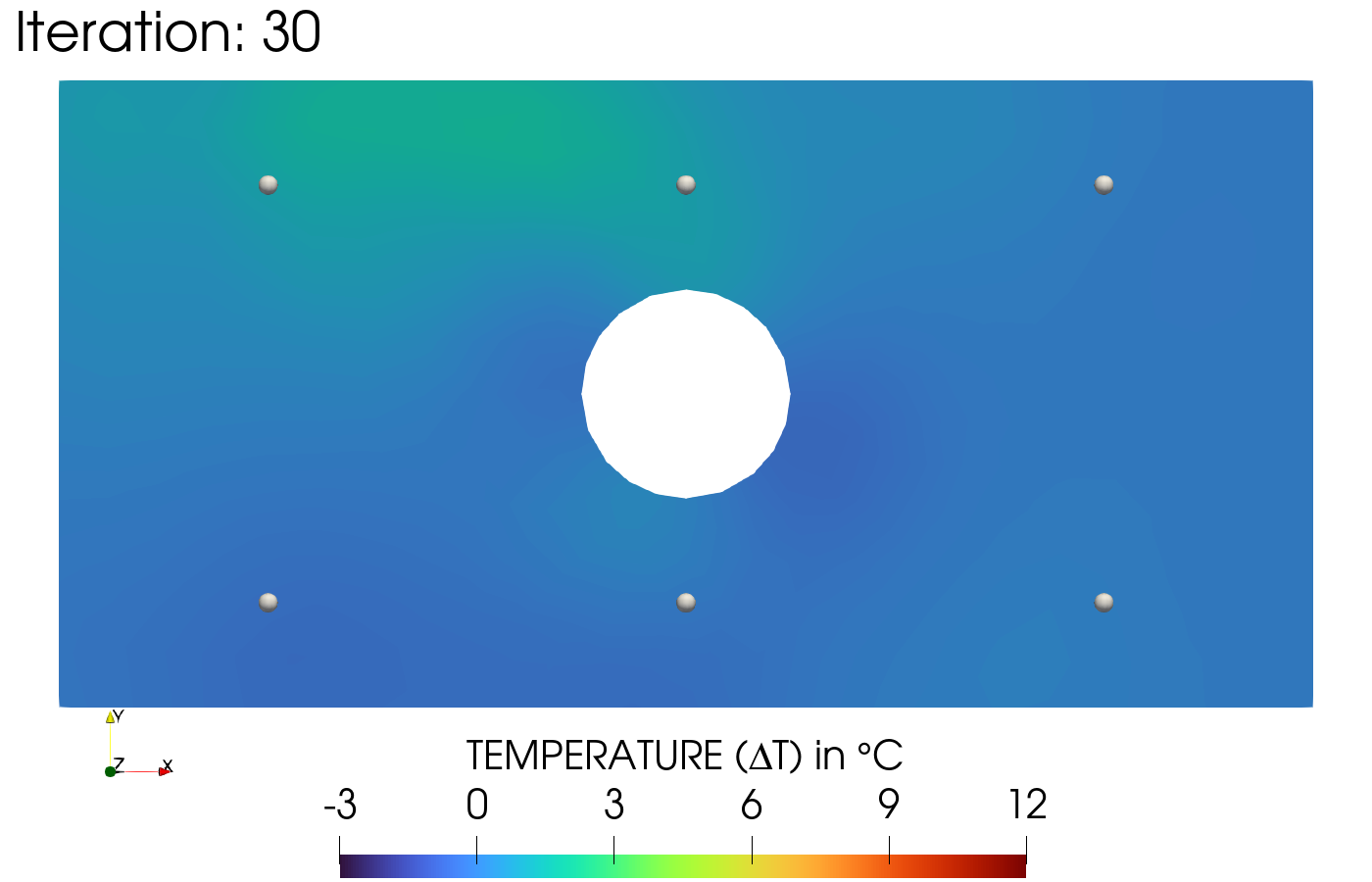}
  \subcaption{With VM ($r=5$): Iteration 30}
   \label{fig:Figure_5d}
   \vspace{0.5em}
  \includegraphics[trim= 0 0 0 0, clip, width=0.3\paperwidth]{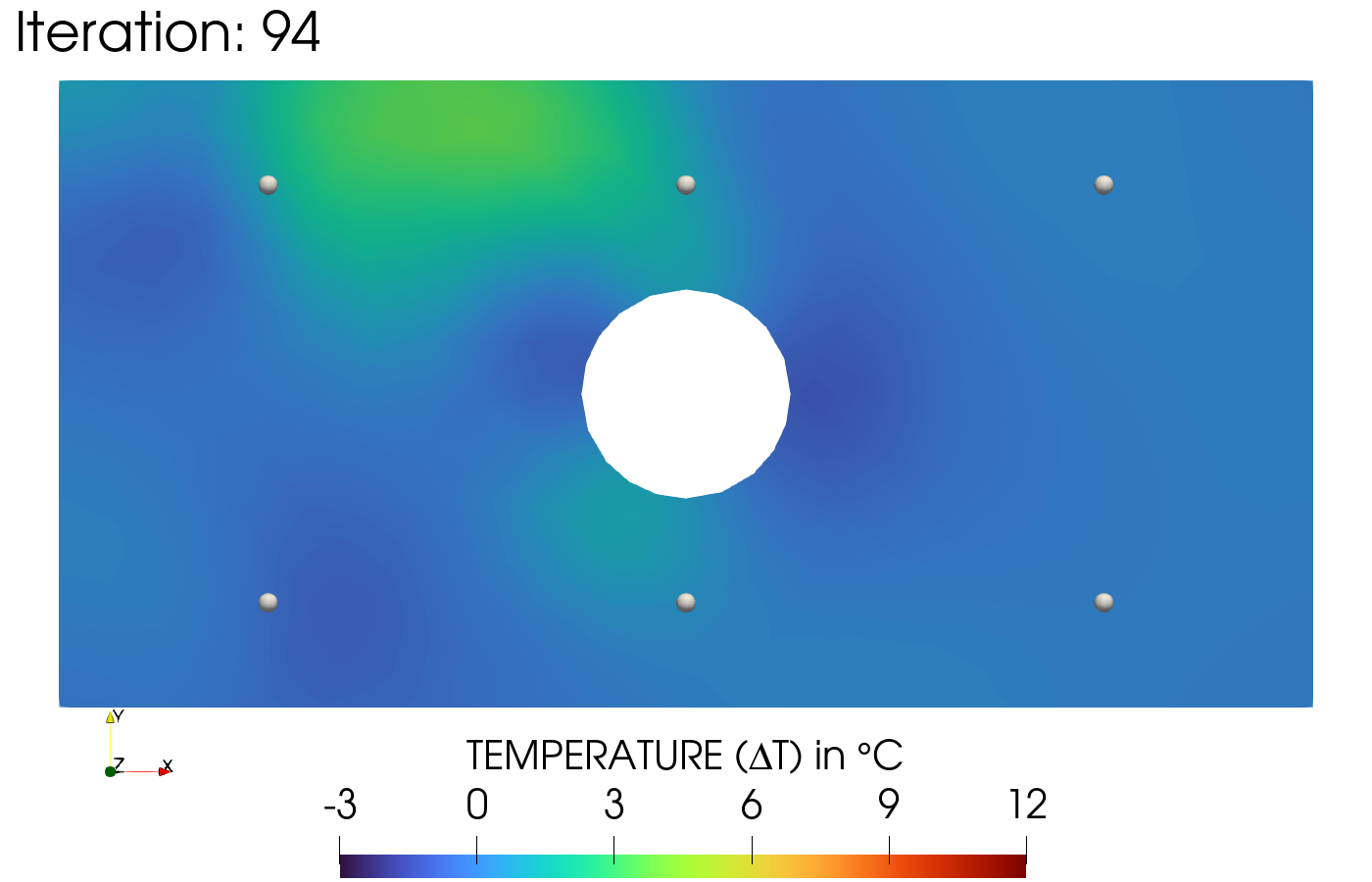}
  \subcaption{With VM ($r=5$): Iteration 94}
   \label{fig:Figure_5e}
  \vspace{0.5em}
  \includegraphics[trim= 0 0 0 0, clip, width=0.3\paperwidth]{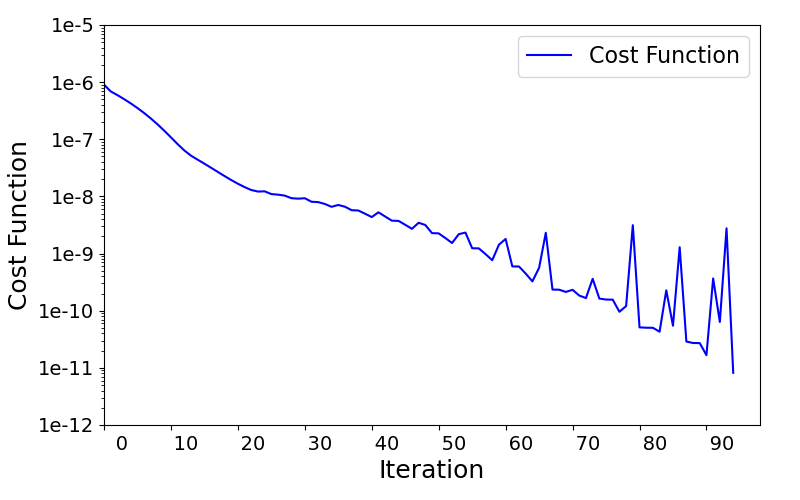}
  \subcaption{With VM ($r=5$): Cost function convergence}
  \label{fig:Figure_5f}\par
\end{minipage}
\caption{Plate With a Hole, 6 sensors configuration: Temperature distribution obtained without and with Vertex Morphing (VM) filtering at different iterations and the cost function convergence plots.}
\label{fig:Figure_5}
\end{figure}

The results obtained at iteration 30 and the final iterations along with the cost function convergence plots for the 6 sensors configuration, 'without' and 'with' Vertex Morphing filtering are shown in Figures \ref{fig:Figure_5a}-\ref{fig:Figure_5c} and \ref{fig:Figure_5d}-\ref{fig:Figure_5f} respectively. 
It can be observed from Figures \ref{fig:Figure_5a} and \ref{fig:Figure_5d} that a rough profile of the temperature distribution is already obtained at 30 iterations. Up to this iteration, the cost function value has also dropped around 2 magnitudes compared to the starting point. Further optimization iterations focus on refining the identified region.

Looking at the converged results in Figures \ref{fig:Figure_5b} and \ref{fig:Figure_5e}, it can be noticed that the maximum $\Delta T$ identified is around $3.6$ $\degree C$, which is a lot less than the prescribed temperature.

The spiking behavior of the cost function convergence (seen in Figures \ref{fig:Figure_5c} and \ref{fig:Figure_5f}) is a well-observed phenomenon when using the Barzilai-Borwein stepsize method (\cite{fletcher2005barzilai}). Since this work does not focus on the optimization algorithm but rather the end result of the optimization, any optimization algorithm (in this case, the Steepest Descent algorithm with the Barzilai-Borwein method) can be used. 

The effect of Vertex Morphing filtering can be seen by comparing Figure \ref{fig:Figure_5b} with Figure \ref{fig:Figure_5e}. Localized sharp temperatures are observed when no filtering is used. Vertex Morphing helps smooth out these localized effects.
The 'with' Vertex Morphing case is noted to require a higher number of iterations to converge when compared to the 'without' Vertex Morphing case. 

\medskip

\begin{figure}[!t]
\begin{minipage}[c][][t]{.5\textwidth}
  \vspace*{\fill}
  \centering
   \includegraphics[trim= 0 0 0 0, clip, width=0.3\paperwidth]{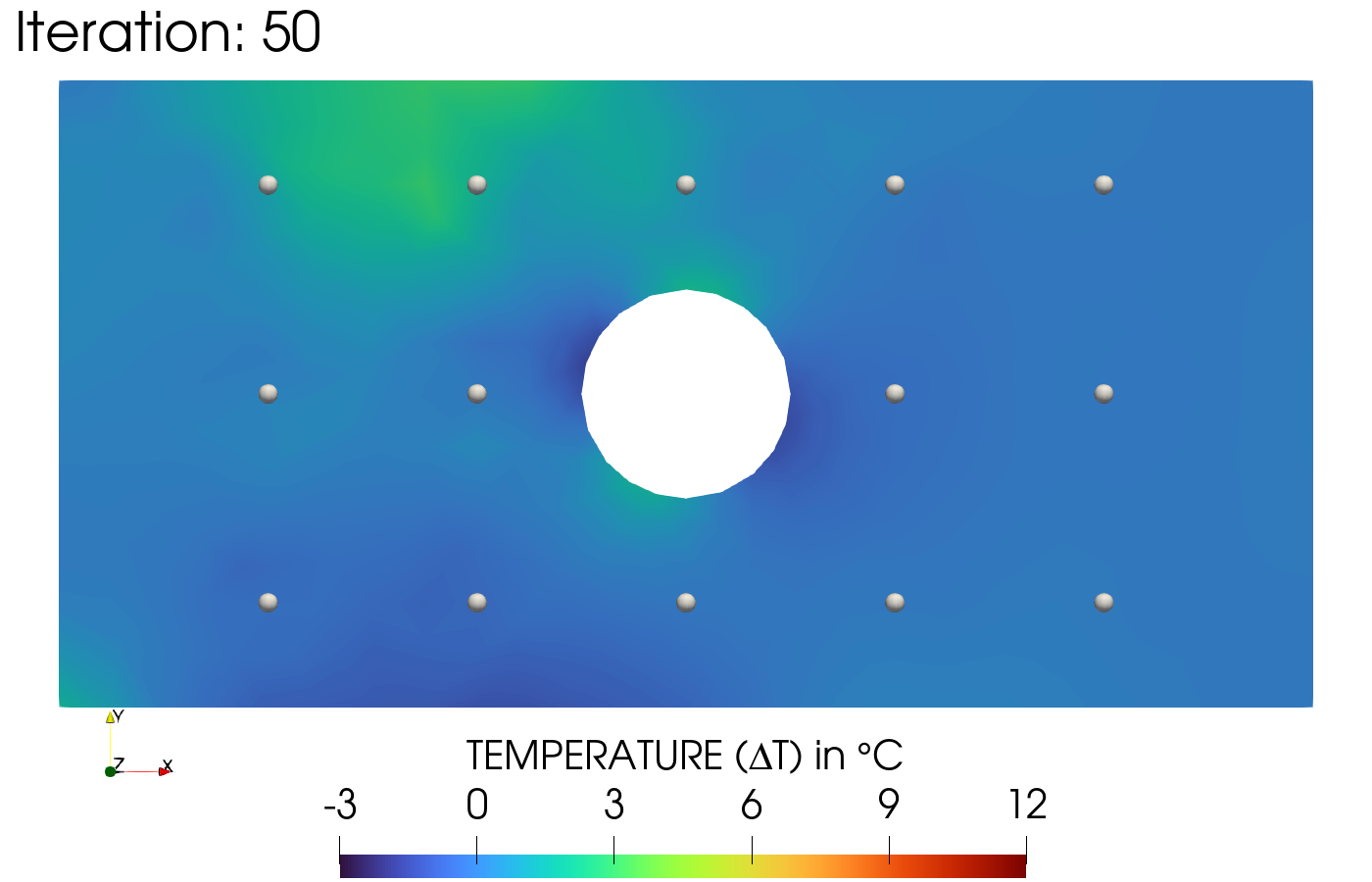}
  \subcaption{Without VM: Iteration 50}
   \label{fig:Figure_6a}
    \vspace{0.5em}
  \includegraphics[trim= 0 0 0 0, clip, width=0.3\paperwidth]{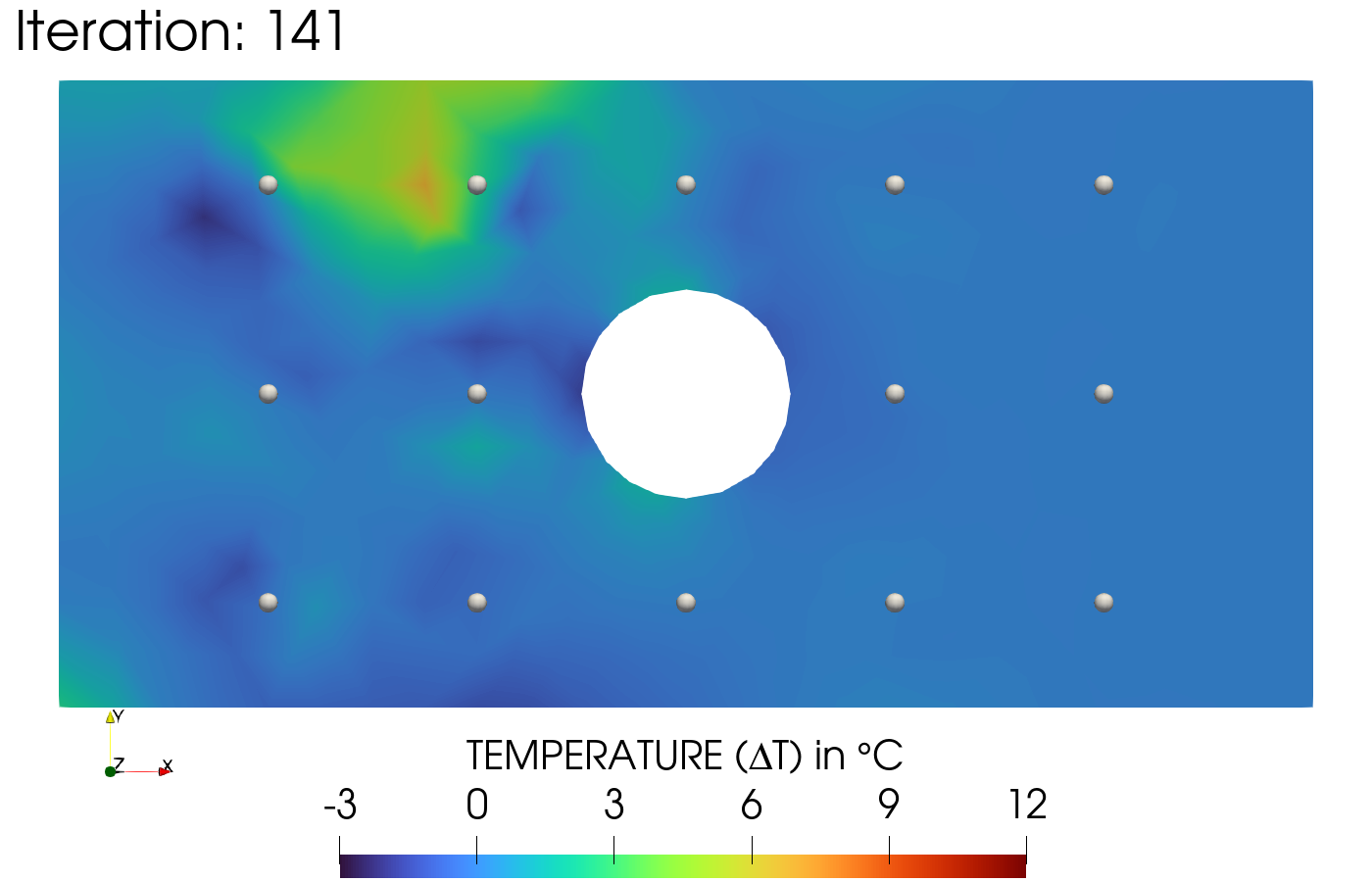}
  \subcaption{Without VM: Iteration 141}
   \label{fig:Figure_6b}
    \vspace{0.5em}
  \includegraphics[trim= 0 0 0 0, clip, width=0.3\paperwidth]{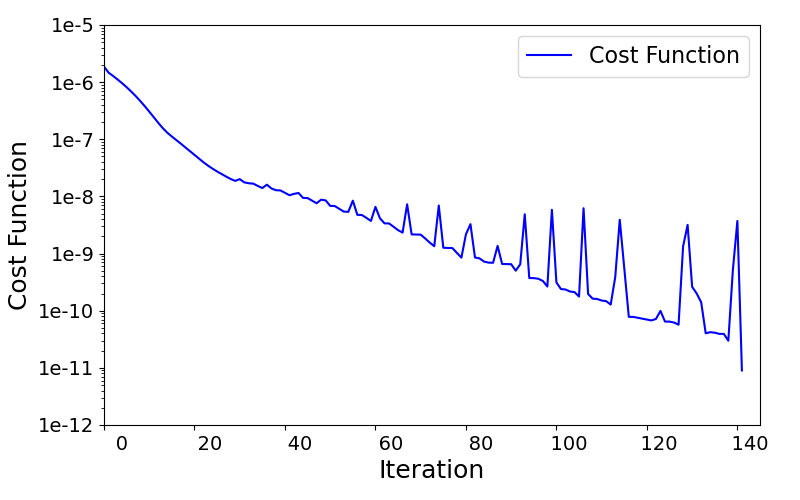}
  \subcaption{Without VM: Cost function convergence}
  \label{fig:Figure_6c}\par
\end{minipage}
\begin{minipage}[c][][t]{.5\textwidth}
  \vspace*{\fill}
  \centering
  \includegraphics[trim= 0 0 0 0, clip, width=0.3\paperwidth]{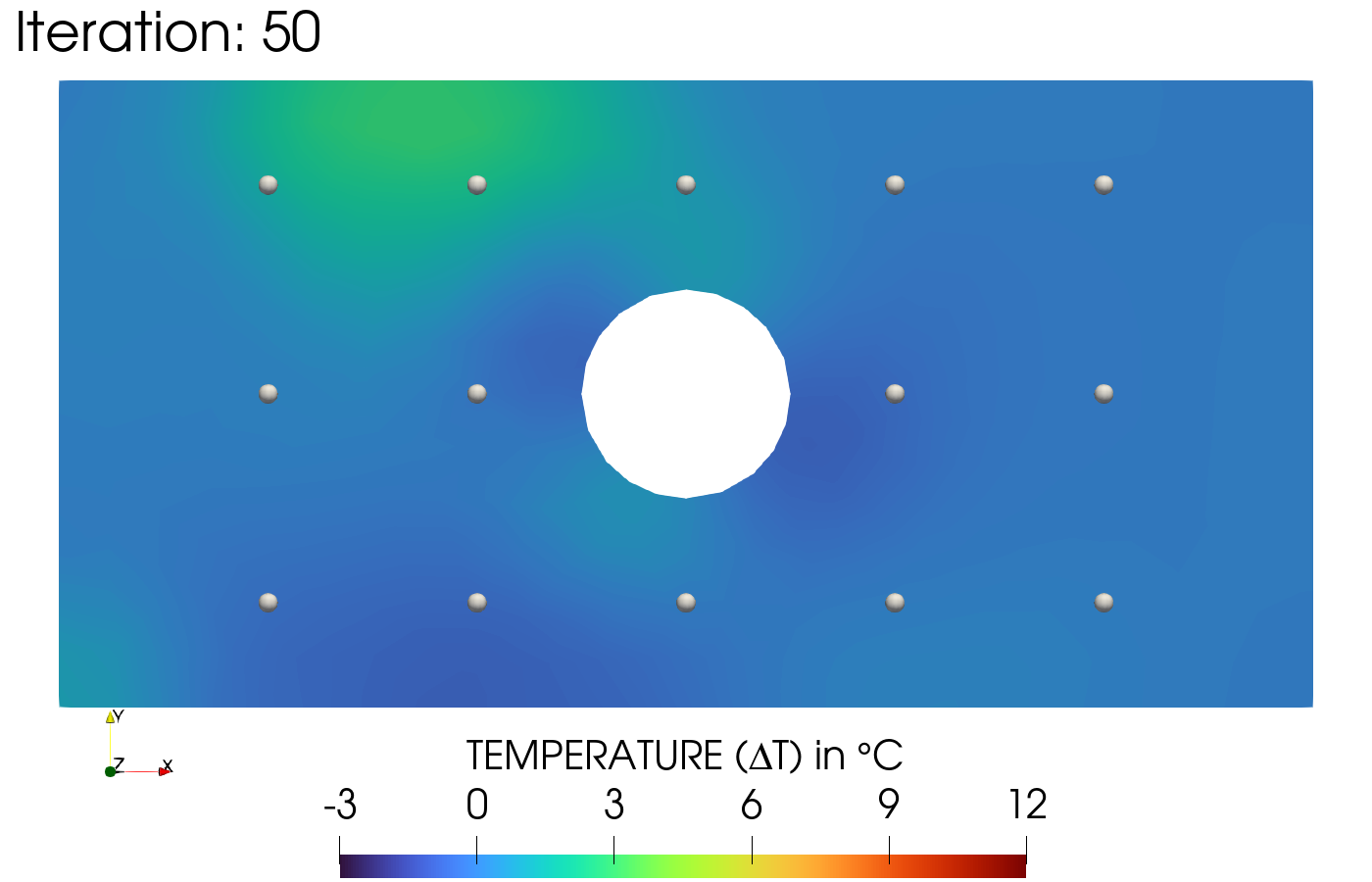}
  \subcaption{With VM ($r=5$): Iteration 50}
   \label{fig:Figure_6d}
    \vspace{0.5em}
  \includegraphics[trim= 0 0 0 0, clip, width=0.3\paperwidth]{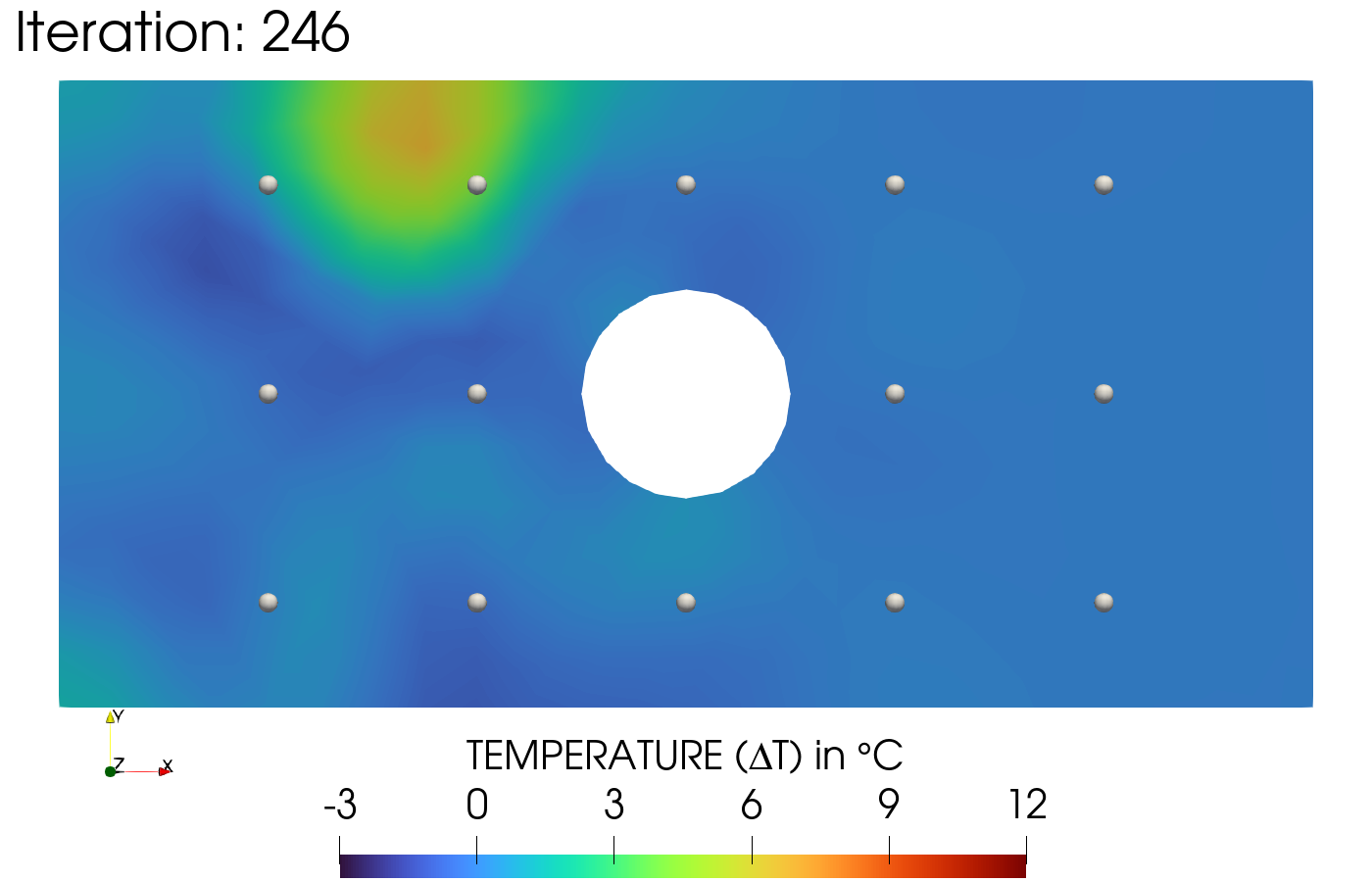}
  \subcaption{With VM ($r=5$): Iteration 246}
   \label{fig:Figure_6e}
    \vspace{0.5em}
  \includegraphics[trim= 0 0 0 0, clip, width=0.3\paperwidth]{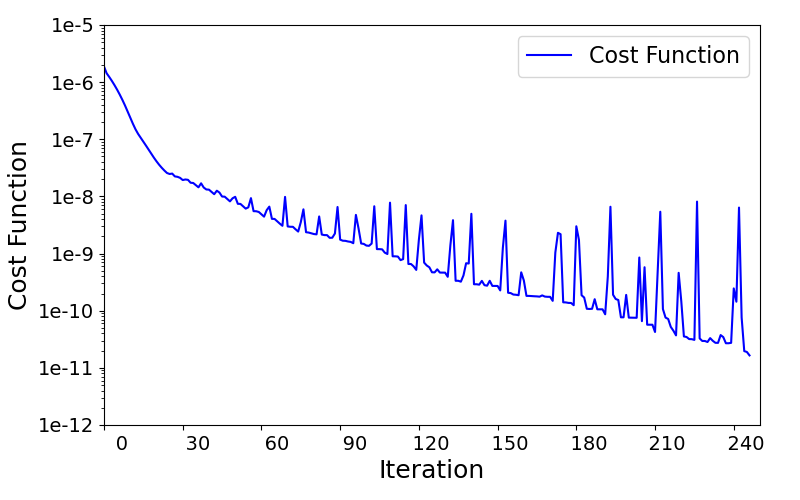}
  \subcaption{With VM ($r=5$): Cost function convergence}
  \label{fig:Figure_6f}\par
\end{minipage}
\caption{Plate With a Hole, 14 sensors configuration: Temperature distribution obtained without and with Vertex Morphing (VM) filtering at different iterations and the cost function convergence plots.}
\label{fig:Figure_6}
\end{figure}

\medskip
Figures \ref{fig:Figure_6a}-\ref{fig:Figure_6c} and \ref{fig:Figure_6d}-\ref{fig:Figure_6f} show the results obtained at iteration 50 and the final iterations along with the cost function convergence plots for the 14 sensors configuration, 'without' and 'with' Vertex Morphing filtering respectively. 
Similar to the 6 sensors configuration, it can be observed from Figures \ref{fig:Figure_6a} and \ref{fig:Figure_6d} that a hazy identification of the heated region has been made in about 50 iterations. Approximately 2 magnitudes drop in the cost function has been achieved until this iteration, with further cost function reduction contributing towards localization.

From Figures \ref{fig:Figure_6b} and \ref{fig:Figure_6e}, the conspicuous effect of Vertex Morphing can be observed. The sharp temperature gradients are abated and a smooth temperature distribution is obtained. Although not the exact prescribed temperature distribution is recovered, $\Delta T$ as high as $6.8$ $\degree C$ is seen in this case, which is a significant improvement over the peak identified in the 6 sensors configuration case. 
As observed earlier, the case with Vertex Morphing takes longer to converge compared to when no filtering is used.

Comparing the 6 and 14 sensors configurations, It is clearly evident that the greater the number of sensors, the better the identification and localization of the temperature field. Hence, sensor quantity, as well as sensor density, play a major role in determining the accuracy of the identification. In the 6 sensors configuration, the algorithm identifies the heated region to be somewhere between the left and the middle sensors of the top row (cf. Figure \ref{fig:Figure_5e}). The 14 sensors configuration further narrows down that region using information from the additional sensors (cf. Figure \ref{fig:Figure_6e}). 

The slightly higher cost function value at the start of the optimization in Figures \ref{fig:Figure_6c},\ref{fig:Figure_6f} when compared to Figures \ref{fig:Figure_5c},\ref{fig:Figure_5f}, is due to the higher number of sensors, thereby increasing the cost function. Therefore, for a fair comparison, the convergence criterion was chosen in terms of magnitude reduction rather than a target value due to the different initial cost function values.

Additionally, it can also be seen that due to the ill-posed nature of the optimization problem i.e., many more design variables than the responses, some negative $\Delta T$'s (up to $-2.2$ $\degree C$ in the 'without' Vertex Morphing cases) are also observed, for instance, in the vicinity of the sensors and the hole. Vertex Morphing helps in dampening out these artificial artifacts and also improves  the negative $\Delta T$'s by bringing them in the range of $-1$ $\degree C$. 

\medskip

\begin{figure}[!t]
\begin{minipage}[c][][t]{.5\textwidth}
  \vspace*{\fill}
  \centering
   \includegraphics[trim= 0 0 0 0, clip, width=0.3\paperwidth]{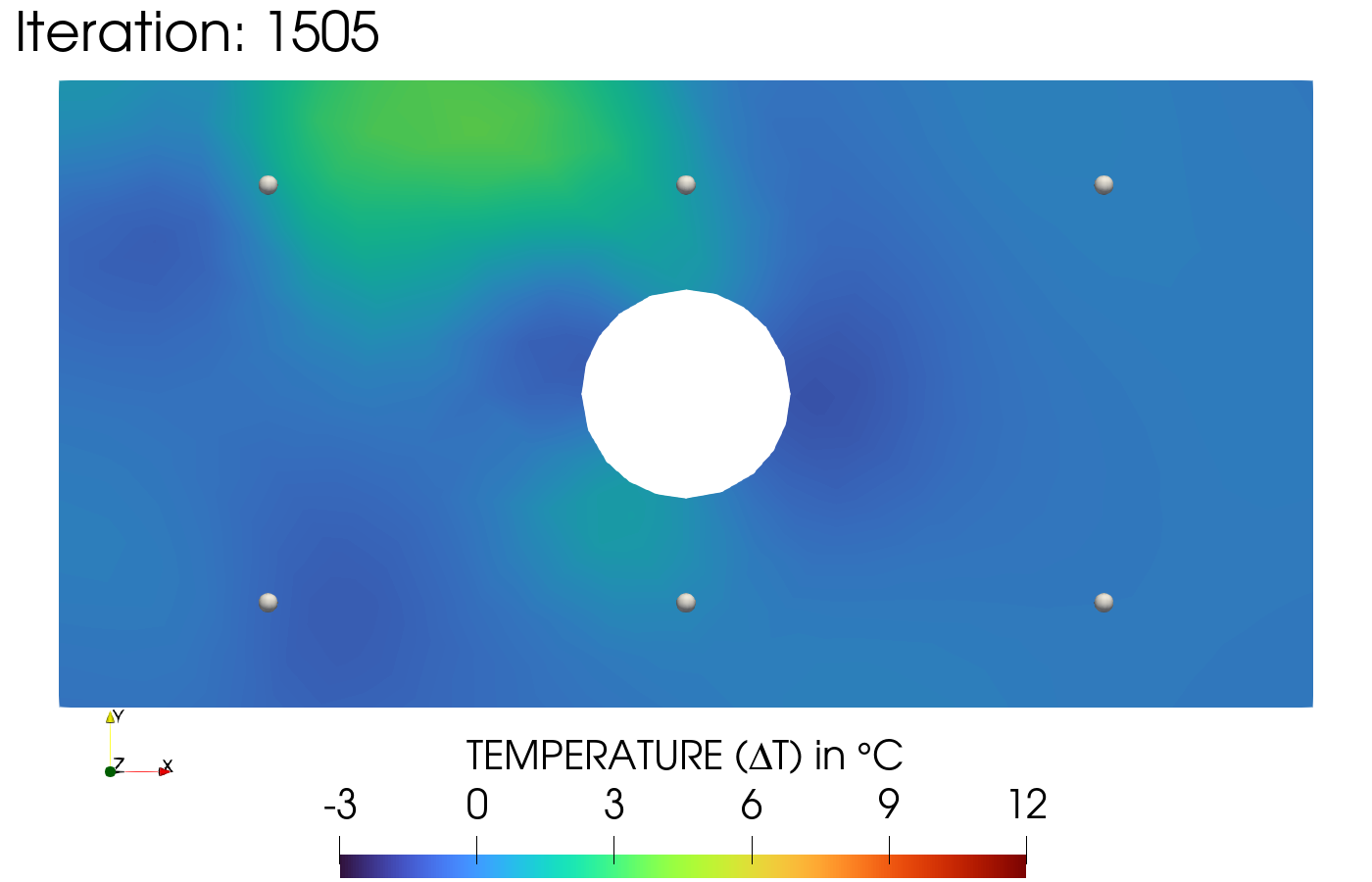}
  \subcaption{6 Sensors: Iteration 1505}
   \label{fig:Figure_7a}
    \vspace{0.5em}
  \includegraphics[trim= 0 0 0 0, clip, width=0.3\paperwidth]{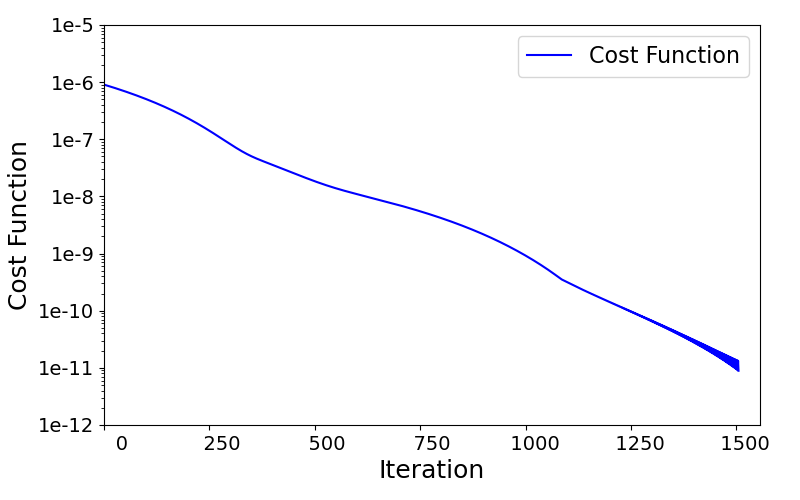}
  \subcaption{6 Sensors: Cost function convergence}
  \label{fig:Figure_7b}\par
\end{minipage}
\begin{minipage}[c][][t]{.5\textwidth}
  \vspace*{\fill}
  \centering
  \includegraphics[trim= 0 0 0 0, clip, width=0.3\paperwidth]{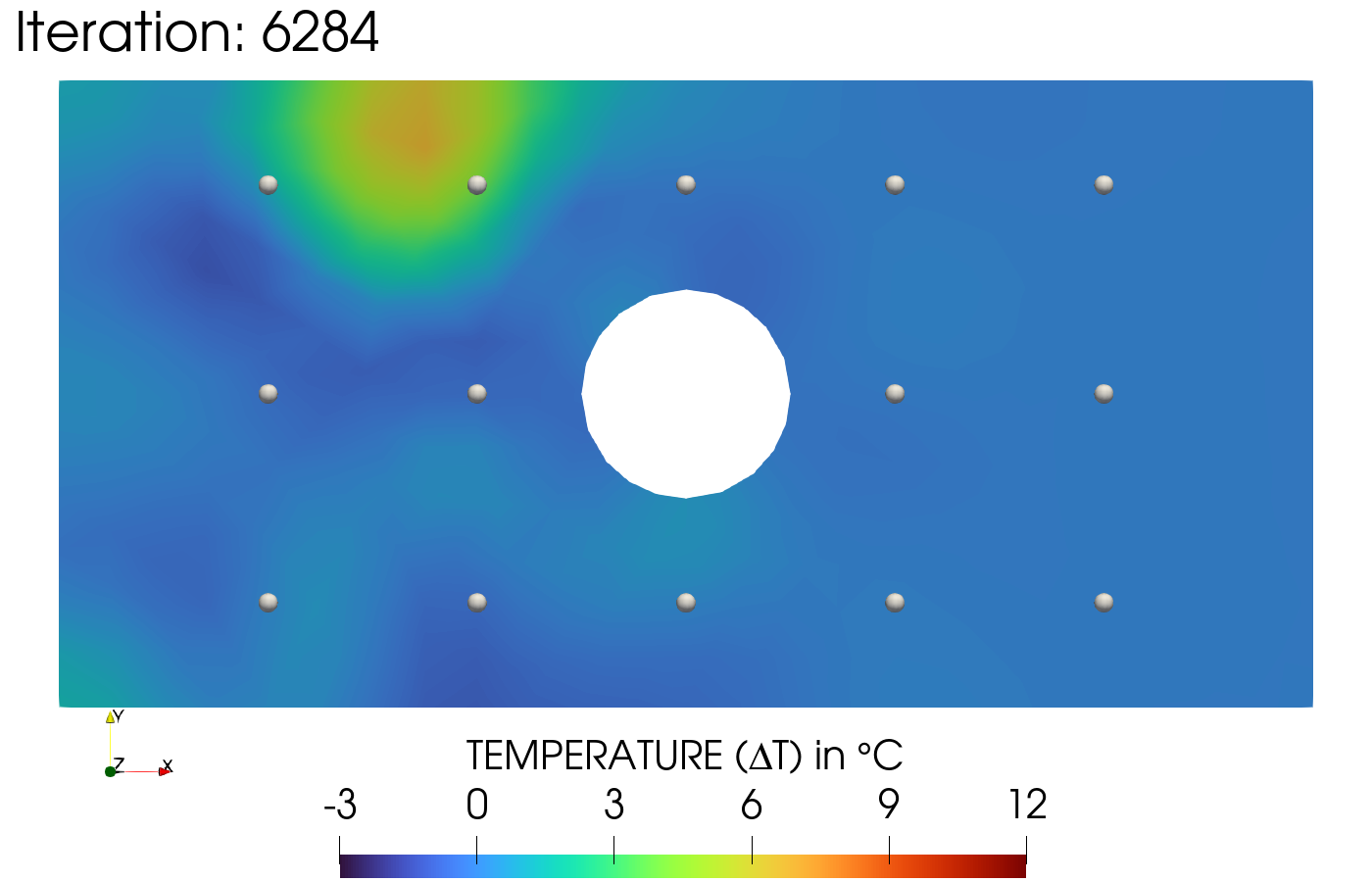}
  \subcaption{14 Sensors: Iteration 6284}
   \label{fig:Figure_7c}
    \vspace{0.5em}
  \includegraphics[trim= 0 0 0 0, clip, width=0.3\paperwidth]{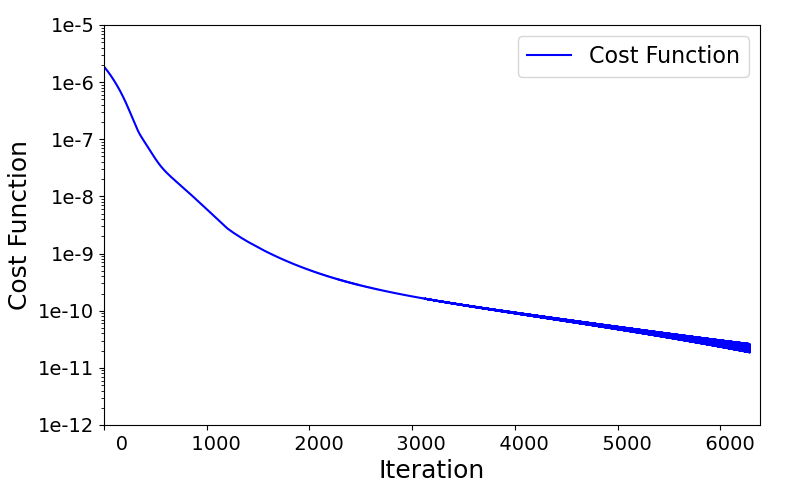}
  \subcaption{14 Sensors: Cost function convergence}
  \label{fig:Figure_7d}\par
\end{minipage}
\caption{Plate With a Hole, 6 and 14 sensors configurations: Temperature distribution obtained with Vertex Morphing ($radius=5$) filtering using steepest descent and constant step of $1\cdot 10^{-4}$ at final iterations and the cost function convergence plots.}
\label{fig:Figure_7}
\end{figure}

To further attest to the validity of the chosen \textit{steepest descent with \acrshort{BB} step} optimization algorithm, both the cases were also compared with the results obtained from the  \textit{steepest descent with small constant step} method. The small constant stepsize was set to $1\cdot 10^{-4}$ to obtain a smoother convergence. All other parameters, such as cost function, Vertex Morphing radius, convergence criteria, etc, were kept the same. The final iteration results and the cost function convergence plots using the constant stepsize method for the 6 and 14 sensors configurations with Vertex Morphing are shown in Figures \ref{fig:Figure_7a},\ref{fig:Figure_7b} and Figures \ref{fig:Figure_7c}, \ref{fig:Figure_7d} respectively. At a glance, it can be seen that the optimization takes significantly longer to converge for the constant stepsize case compared to the \acrshort{BB} step case. The cost function convergence plots look much smoother in the constant stepsize case. But more importantly, the results obtained from both the constant and \acrshort{BB} stepsize methods look almost exactly the same. Henceforth, the faster converging, steepest descent with the BB step method shall be used with confidence, considering the spiking phenomenon of the cost function convergence as an inherent behavior of the method.

In general, it can be observed from Figures \ref{fig:Figure_5}, \ref{fig:Figure_6}, and \ref{fig:Figure_7} that the exact prescribed temperature distribution is not obtained. The recovered temperatures are spread out over a slightly larger region, and the peak values are lower than the prescribed $\Delta \mathbf{T}$, but still, the distributions are reasonably close.

\FloatBarrier

\subsection{Bridge}
\label{sec:bridge}

\begin{figure}[!b]
    \centering
    \begin{subfigure}[t]{\textwidth}
        \centering
        \begin{minipage}[t]{0.62\textwidth}
            \centering
            \includegraphics[trim=0 77 0 85, clip, width=\textwidth]{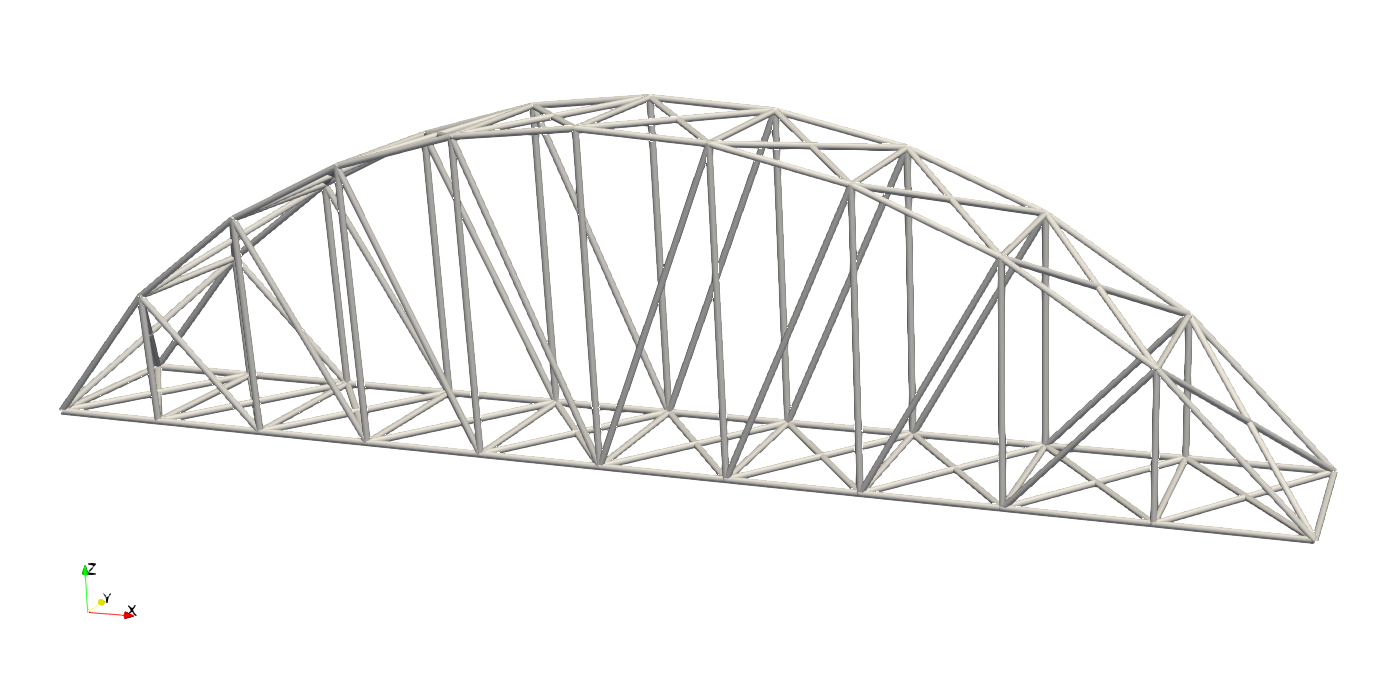}
        \end{minipage}
        \hfill
    \end{subfigure}
    \begin{subfigure}[t]{\textwidth}
        \centering
        \begin{minipage}[t]{0.62\textwidth}
            \centering
            \includegraphics[trim=0 0 0 85, clip, width=\textwidth]{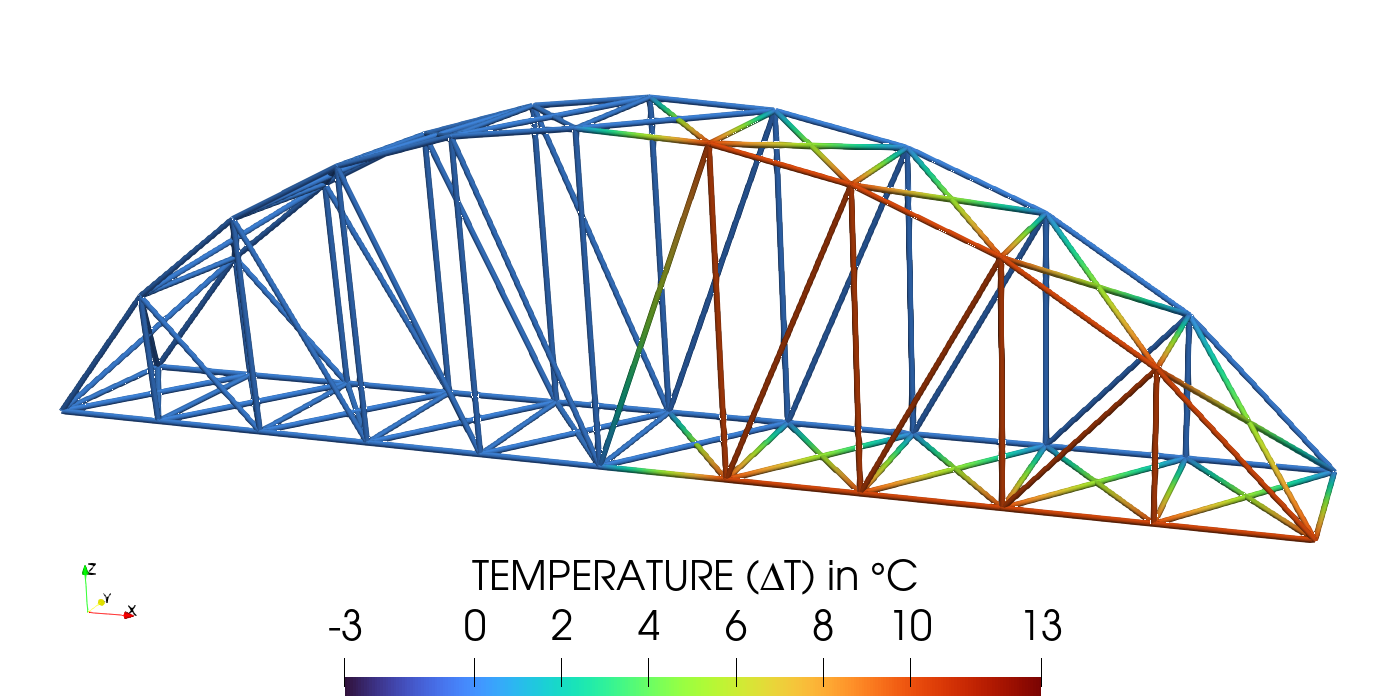}
        \end{minipage}
        \hfill
    \end{subfigure}
    \caption{Bridge example: Mesh used (top) and the target temperature distribution (bottom).}
    \label{fig:Figure_8}
\end{figure}

\begin{figure}[!t]
    \centering
    \begin{subfigure}[t]{\textwidth}
        \centering
        \begin{minipage}[t]{0.62\textwidth}
            \centering
            \includegraphics[trim=0 77 0 0, clip, width=\textwidth]{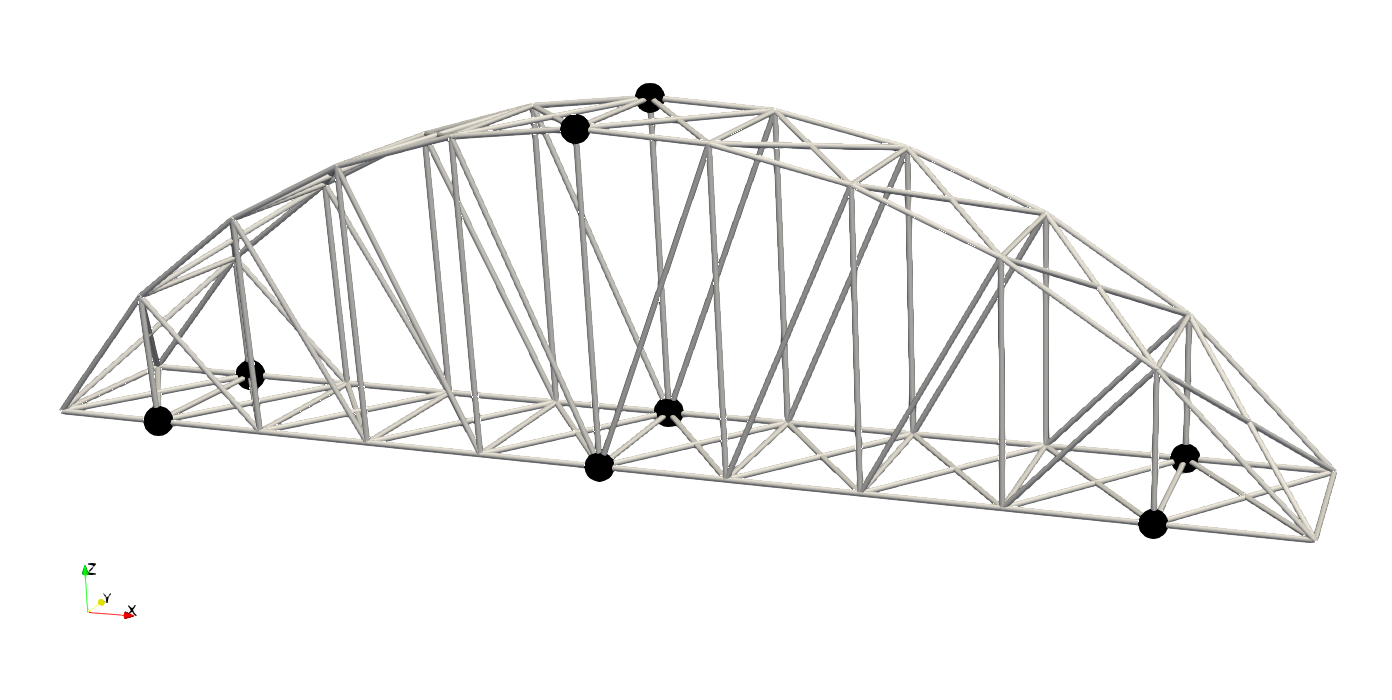}
        \end{minipage}
        \hfill
    \end{subfigure}
    \begin{subfigure}[t]{\textwidth}
        \centering
        \begin{minipage}[t]{0.62\textwidth}
            \centering
            \includegraphics[trim=0 77 0 0, clip, width=\textwidth]{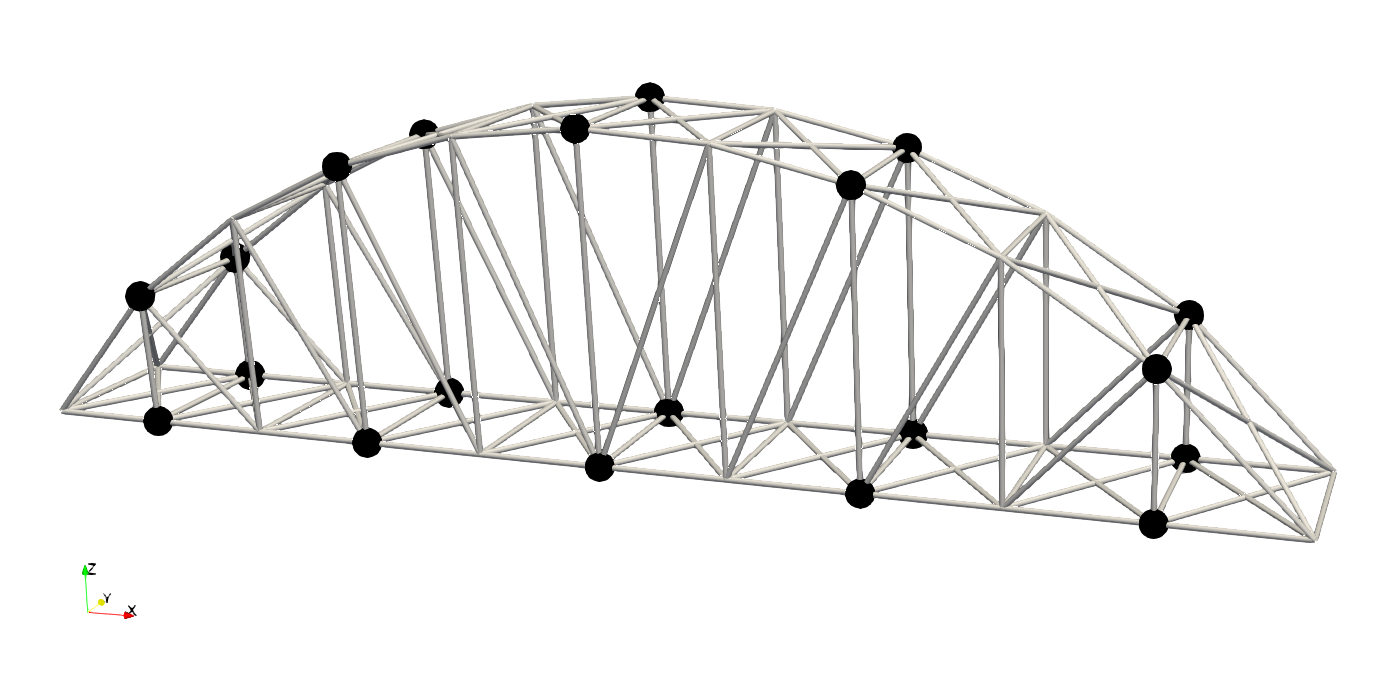}
        \end{minipage}
        \hfill
    \end{subfigure}
    \caption{Bridge example: Sensors distribution for the 8 (top) and 20 (bottom) sensors configurations.}
    \label{fig:Figure_9}
\end{figure}

\begin{figure}[!b]
    \centering
    \begin{subfigure}[t]{\textwidth}
        \centering
        \begin{minipage}[t]{0.62\textwidth}
            \centering
            \includegraphics[trim=0 0 0 0, clip, width=\textwidth]{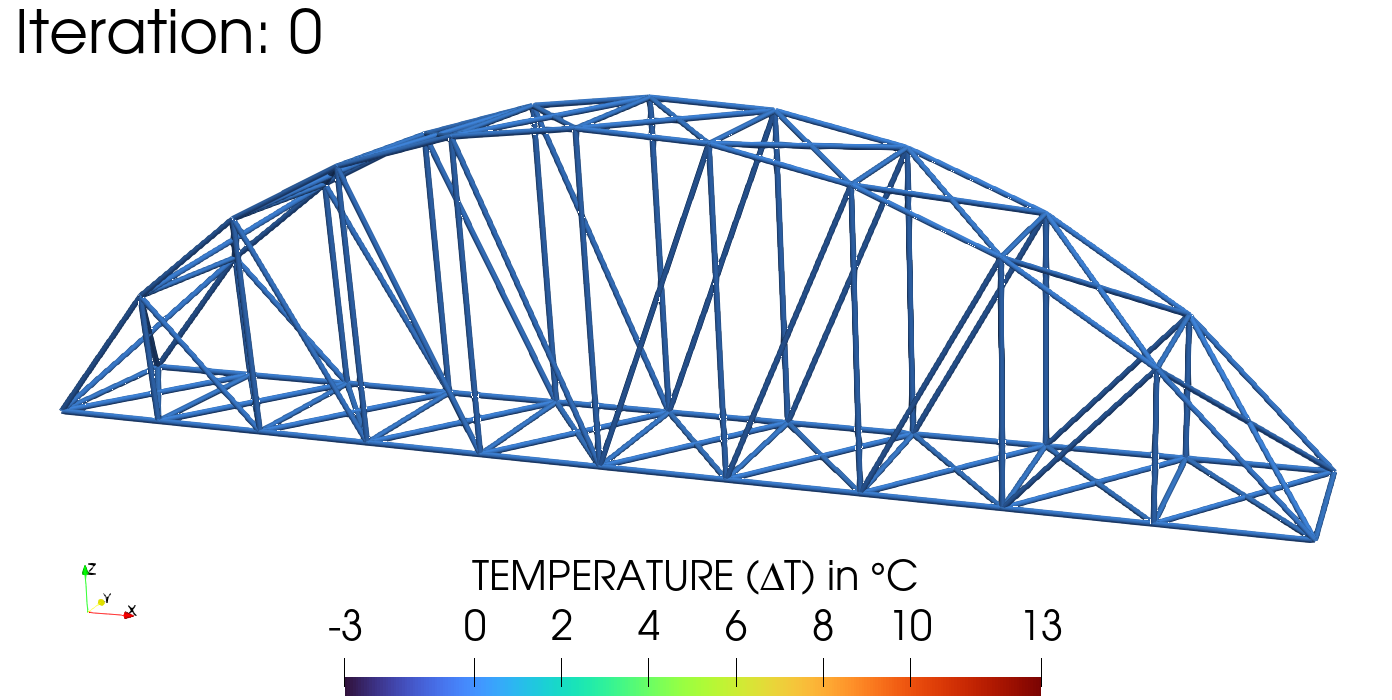}
        \end{minipage}
        \hfill
    \end{subfigure}
    \caption{Bridge example: Temperature Distribution at Optimization Start $\Delta \mathbf{T} = \mathbf{0}$.}
    \label{fig:Figure_10}
\end{figure}

\begin{figure}[!t]
    \centering
    \begin{subfigure}[t]{\textwidth}
        \centering
        \begin{minipage}[t]{0.62\textwidth}
            \centering
            \includegraphics[trim=0 0 0 0, clip, width=\textwidth]{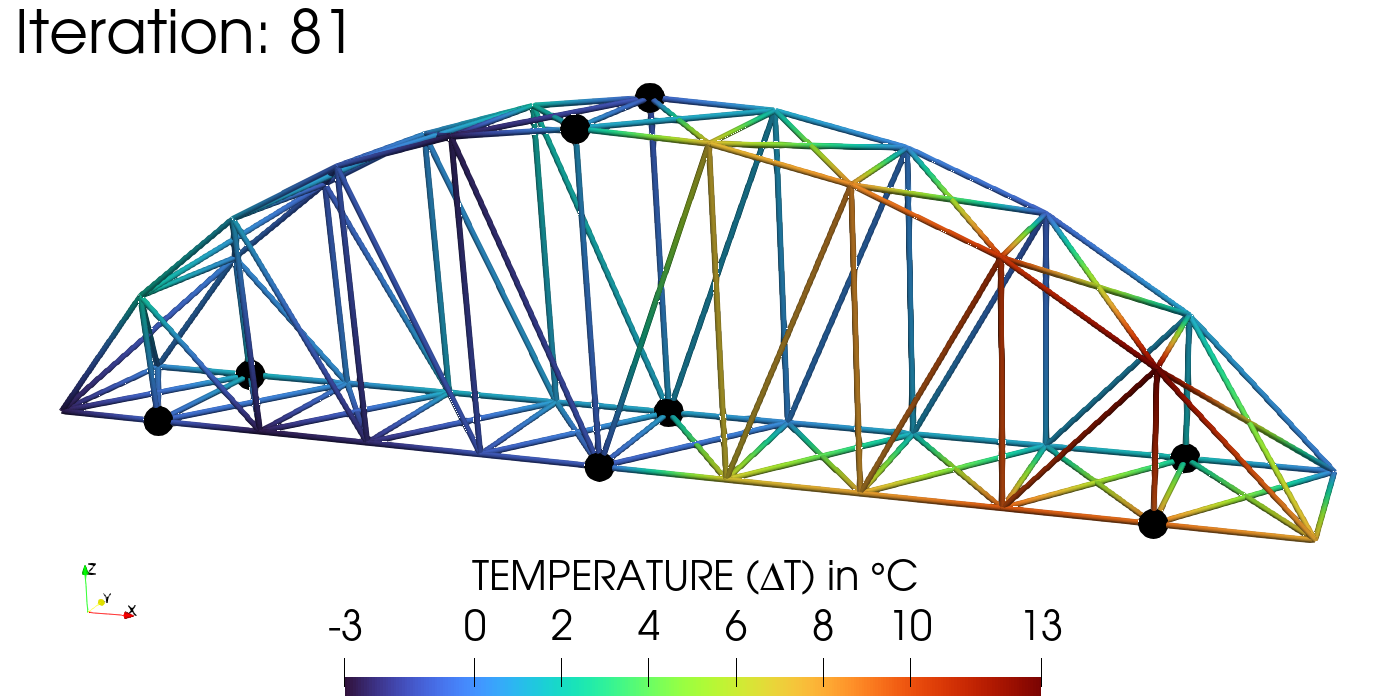}
        \end{minipage}
        \hfill
        \begin{minipage}[t]{0.37\textwidth}
            \centering
            \includegraphics[trim=0 0 0 0, clip, width=\textwidth]{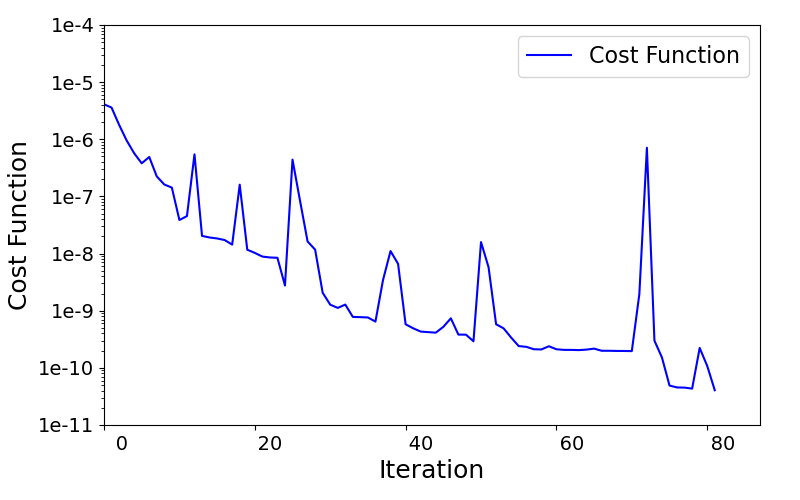}
        \end{minipage}
        \caption{Without Vertex Morphing filtering: Iteration 81}
        \label{fig:Figure_11a}
        \vspace{1em}
    \end{subfigure}
    \vspace{1em}
    \begin{subfigure}[t]{\textwidth}
        \centering
        \begin{minipage}[t]{0.62\textwidth}
            \centering
            \includegraphics[trim=0 0 0 -20, clip, width=\textwidth]{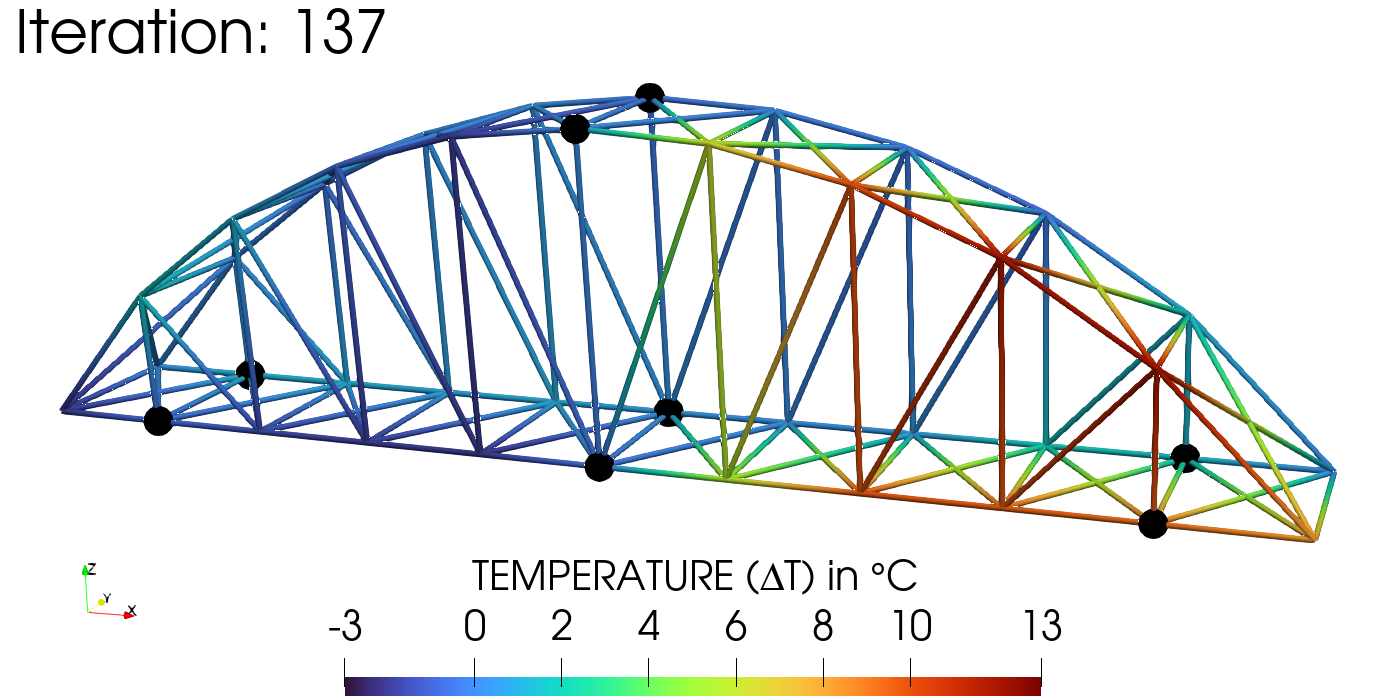}
        \end{minipage}
        \hfill
        \begin{minipage}[t]{0.37\textwidth}
            \centering
            \includegraphics[trim=0 0 0 0, clip, width=\textwidth]{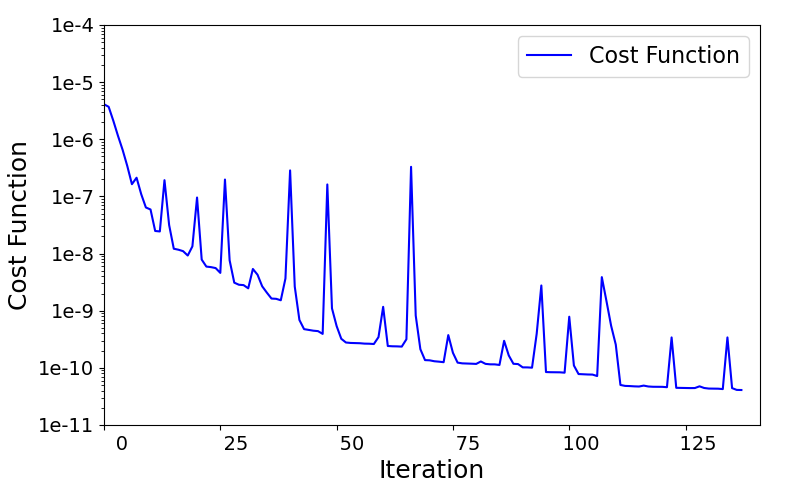}
        \end{minipage}
        \caption{With Vertex Morphing filtering (radius = $6$): Iteration 137}
        \label{fig:Figure_11b}
    \end{subfigure}
    \caption{Bridge example, 8 sensors configuration: Temperature distributions (left) and the cost function convergence plots (right) obtained without and with Vertex Morphing filtering.}
    \label{fig:Figure_11}
\end{figure}

\begin{figure}[!t]
    \centering
    \begin{subfigure}[t]{\textwidth}
        \centering
        \begin{minipage}[t]{0.62\textwidth}
            \centering
            \includegraphics[trim=0 0 0 0, clip, width=\textwidth]{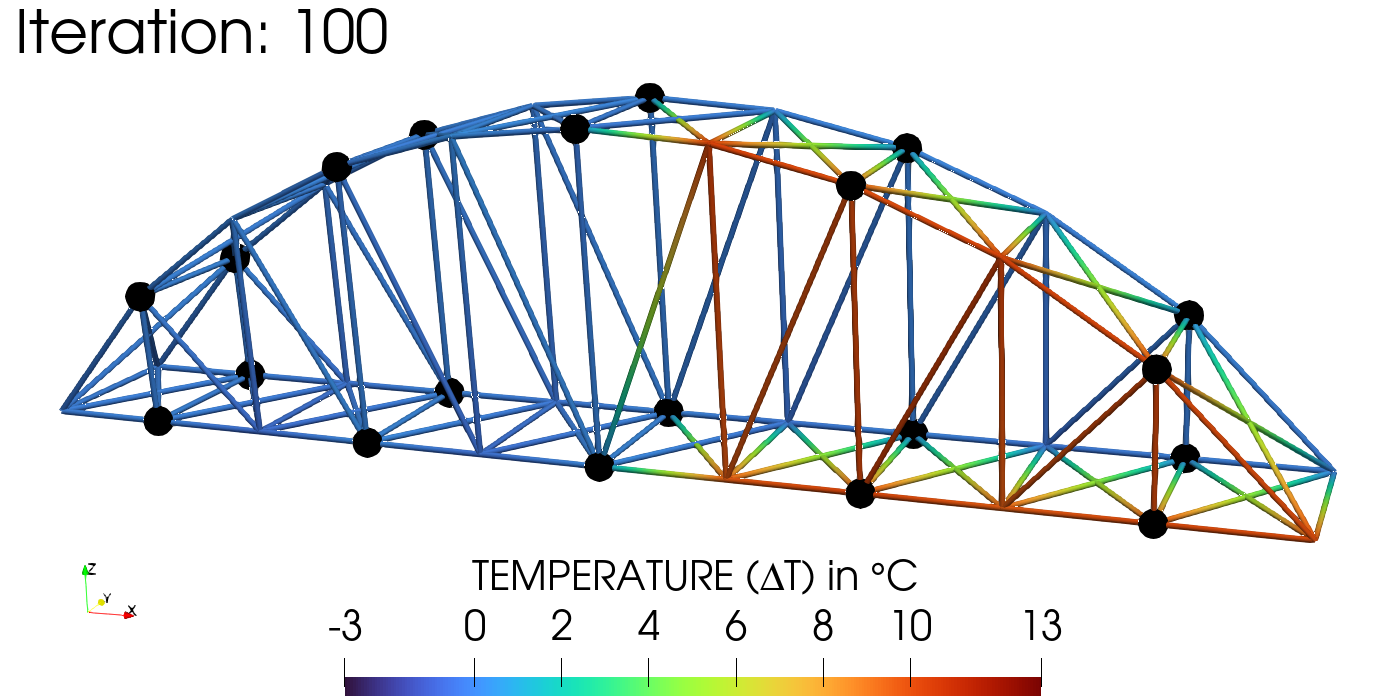}
        \end{minipage}
        \hfill
        \begin{minipage}[t]{0.37\textwidth}
            \centering
            \includegraphics[trim=0 0 0 0, clip, width=\textwidth]{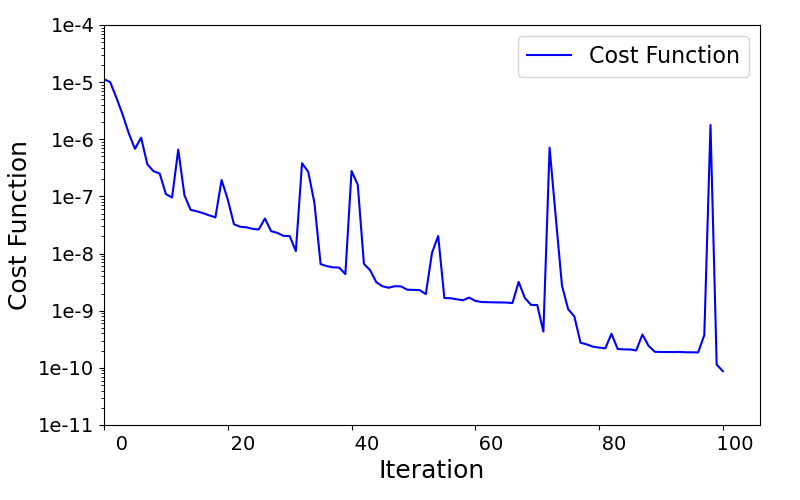}
        \end{minipage}
        \caption{Without Vertex Morphing filtering: Iteration 100}
        \label{fig:Figure_12a}
        \vspace{1em}
    \end{subfigure}
    \vspace{1em}
    \begin{subfigure}[t]{\textwidth}
        \centering
        \begin{minipage}[t]{0.62\textwidth}
            \centering
            \includegraphics[trim=0 0 0 -20, clip, width=\textwidth]{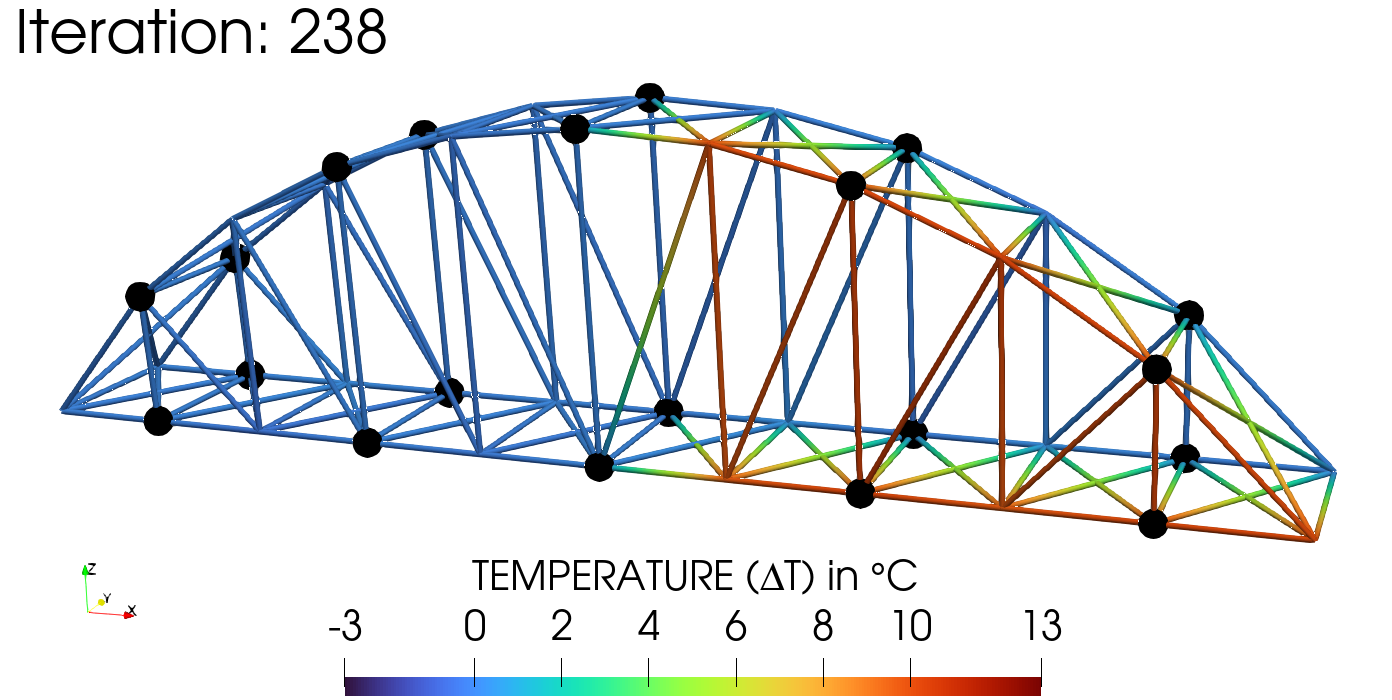}
        \end{minipage}
        \hfill
        \begin{minipage}[t]{0.37\textwidth}
            \centering
            \includegraphics[trim=0 0 0 0, clip, width=\textwidth]{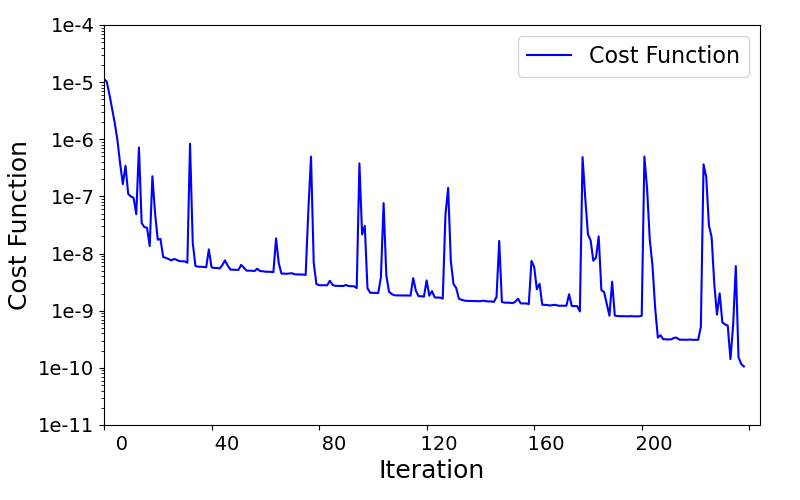}
        \end{minipage}
        \caption{With Vertex Morphing filtering (radius = $6$): Iteration 238}
        \label{fig:Figure_12b}
    \end{subfigure}
    \caption{Bridge example, 20 sensors configuration: Temperature distributions (left) and the cost function convergence plots (right) obtained without and with Vertex Morphing filtering.}
    \label{fig:Figure_12}
\end{figure}

The case is shown in Figures \ref{fig:Figure_8}-\ref{fig:Figure_12}. The bridge has a span of 40 $meters$ (in $x$-direction with range $-20 \leq x \leq 20$), has a maximum height of 10 $meters$ (in $z$-direction with range $0 \leq z \leq 10$), and a width of 5 $meters$ (in $y$-direction with range $0 \leq y \leq 5$). The trusses have areas ranging from 1 $cm^2$ to 100 $cm^2$. The bridge is symmetric length-wise about the $x=0$ plane. The front and rear length-wise sections of the bridge are located at $y=0$ and $y=5$ planes respectively (refer Figure \ref{fig:Figure_8}). 134 linear truss elements were used for the finite element discretization, which is shown in Figure \ref{fig:Figure_8} (top). The left and right edges were assumed clamped ($\mathbf{u}= \mathbf{0}$) and load due to self-weight and gravity of $9.81$ $m/s^2$ was considered acting in the negative $z$-direction.
Density, Young's modulus, Poisson ratio, and thermal expansion coefficient were set to $\rho = 7800$, $E = 2\cdot 10^{11}$, $\nu = 0.3$, $\alpha = 1.0 \cdot 10^{-5} /\degree K$ respectively in SI units. A temperature difference of $\Delta T = 10$ $\degree C$ was prescribed on the front-right section of the bridge and is shown in Figure \ref{fig:Figure_8} (bottom). Everywhere else, the temperature difference was prescribed to be $\Delta T = 0$ $\degree C$.

Two sensor configurations consisting of 8 and 20 displacement sensors distributed over the bridge were analyzed. Figure \ref{fig:Figure_9} shows the sensor locations for both sensor configurations. For both the sensor configurations, the steepest descent optimization algorithm with the Barzilai-Borwein stepsize method (Eqn.\eqref{eq:bb}) with a maximum stepsize of $2.5 \cdot 10^{-2}$ was used. The optimization starts with the presumption that the temperature $\Delta T = 0$ $\degree C$ everywhere on the bridge. The temperature distribution at the start of the optimization is shown in Figure \ref{fig:Figure_10}.
The same stopping criteria as for the Plate With a Hole example was set: 5 magnitudes reduction of the initial cost function i.e., the algorithm stops when the $\text{current cost function} \leq (1\cdot 10^{-5} * \text{initial cost function})$. 
Vertex Morphing with radius $r = 6$ was used to smoothen the gradients and subsequently get smoother results. Considering the width of the bridge and the length of the trusses, a radius of $r=6$ encompassed approximately $3-4$ neighboring nodes, allowing for a good regularization.
For better visualization of the results, the color bar is shown in the range $[-3,13]$ $\degree C$ as the temperature ($\Delta T$) is above $-3$ $\degree C$ in all but one case, and in all cases it is below $13$ $\degree C$. 

\medskip

The temperature distributions obtained at the final iterations and the cost function convergence plots for the 8 sensors configuration 'without' and 'with' Vertex Morphing filtering are shown in Figures \ref{fig:Figure_11a} and \ref{fig:Figure_11b} respectively. 
In general, it can be observed that a fairly good reconstruction of the thermal field is obtained in both, 'with' and 'without' Vertex Morphing filtering cases. The heated and non-heated regions are localized quite well.

For the case without Vertex Morphing (Figure \ref{fig:Figure_11a}), some nodes are observed to have larger positive/negative deviations from the prescribed temperatures. The case with Vertex Morphing has slightly better temperature distribution due to the smoothing of gradients and abating sharp peak temperatures.
The temperature distribution in the case without Vertex Morphing is in the range $[-2.65, 12.89]$ $\degree C$, whereas it is in the range $[-2.05, 11.9]$ $\degree C$ for the case with Vertex Morphing, indicating the effect of regularization. 
From the cost function convergence plots, it can be seen that the case with Vertex Morphing requires more iterations to converge compared to the case with no filtering.

\medskip

The temperature distributions obtained at the final iterations and the cost function convergence plots for the 20 sensors configuration 'without' and 'with' Vertex Morphing filtering are shown in Figures \ref{fig:Figure_12a} and \ref{fig:Figure_12b} respectively. Similar to the 8 sensors configuration case, the temperature distribution is successfully reconstructed with good localization of the heated region for both 'with' and 'without' filtering cases.

The initial cost function value seen in the convergence plots is higher for the 20 sensors configuration (Figure \ref{fig:Figure_12}) relative to the 8 sensors configuration (Figure \ref{fig:Figure_11}) because of the aggregation of sensor errors in the objective function (Eqn.\eqref{eq:cost_function}).

As expected, the addition of 12 sensors to the 8 sensors configuration substantially improves the thermal field obtained, especially in reducing the sharp peak temperatures. 
The temperature distribution in the case without Vertex Morphing is in the range $[-0.54,10.24]$ $\degree C$, whereas it is in the range $[-0.35,10.16]$ $\degree C$ for the case with Vertex Morphing indicating a marginal improvement due to filtering. 
From the cost function convergence plots, it can be seen that the case with Vertex Morphing requires a much higher number of optimization iterations to converge than the case without Vertex Morphing.

\medskip

Compared to the Plate With a Hole example (consisting of 646 elements and 369 nodes), the Bridge example (consisting of 134 truss elements and 40 nodes) is relatively smaller and simpler, as the number of design parameters (i.e., nodal temperatures), are many-folds fewer. Due to this, the ill-conditioning effect, as seen in the Plate With a Hole example, is not very prominent, and sufficiently good thermal fields are obtained even without the use of regularization or filtering. However, in all cases, Vertex Morphing improved the overall temperature distribution, also evident by the root-mean-squared errors discussed later in Section \ref{sec:comparison_error}.

Comparing the temperature distribution ranges for the 8 and 20 sensors configuration cases, it can be observed that the 20 sensors configuration temperature range is much closer (narrower) to the actual prescribed temperatures of $\Delta T = 0$ and $10$ $\degree C$. This is in line with the observations for the Plate With a Hole example, where the 14 sensors configuration also gave a significantly better thermal field than the 6 sensors case.

It is noted here as well that the temperature distribution is not recovered exactly due to smoothing and filtering effects, but a very close approximation of the prescribed temperature field is obtained.

\medskip

The heated region is identified and localized fairly quickly within 20\% of the total iterations to converge. This can be observed from the cost function convergence plots in Figures \ref{fig:Figure_11} and \ref{fig:Figure_12}, where approximately $3-4$ magnitudes decrease in the cost function (compared to the initial value), is realized within the initial 20\% of the optimization progress. Within this 20\% of the optimization progress, the thermal field throughout the structure is roughly identified signified by a steep decrease in the cost function. This is followed by a gradual reduction in the slope of the cost function as the optimizer refines the temperature distribution further.

\begin{figure}[!t]
    \centering
    \begin{subfigure}[t]{\textwidth}
        \centering
        \begin{minipage}[t]{0.62\textwidth}
            \centering
            \includegraphics[trim=0 0 0 0, clip, width=\textwidth]{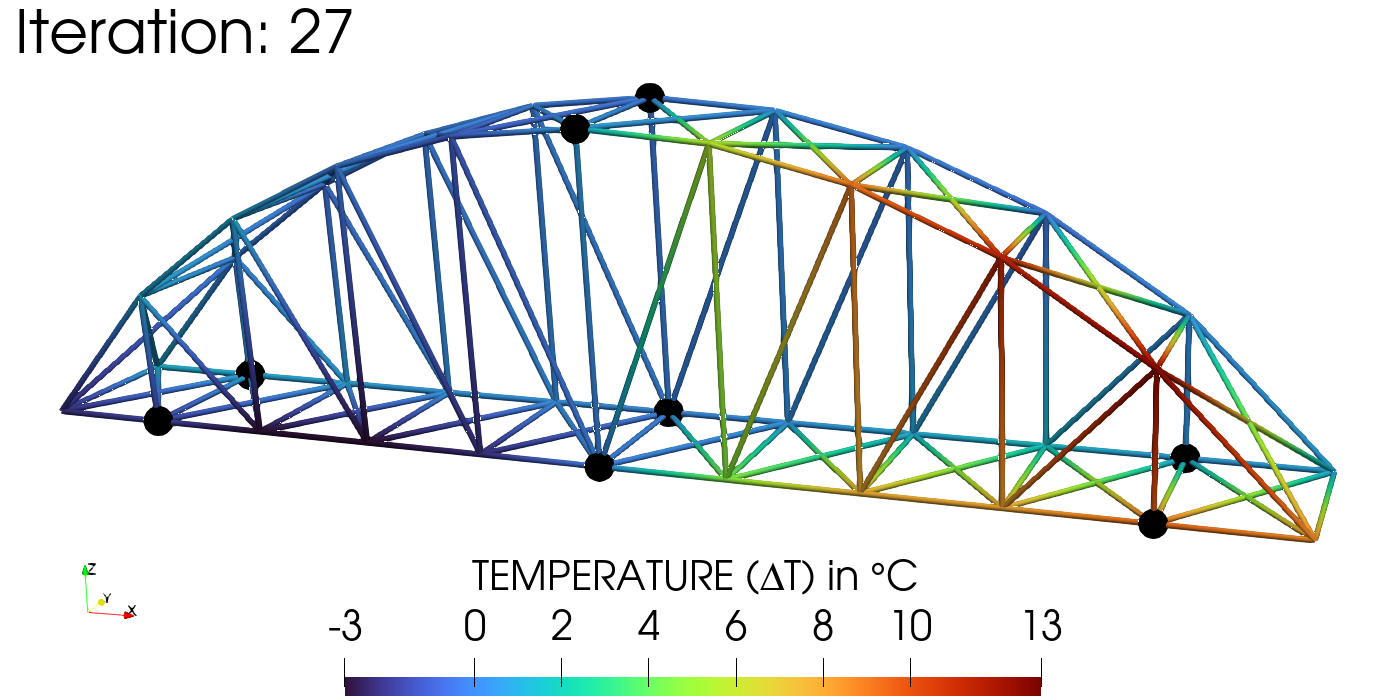}
        \end{minipage}
        \hfill
        \caption{8 sensors: Iteration 27}
        \label{fig:Figure_13a}
        \vspace{1em}
    \end{subfigure}
    \vspace{1em}
    \begin{subfigure}[t]{\textwidth}
        \centering
        \begin{minipage}[t]{0.62\textwidth}
            \centering
            \includegraphics[trim=0 0 0 -20, clip, width=\textwidth]{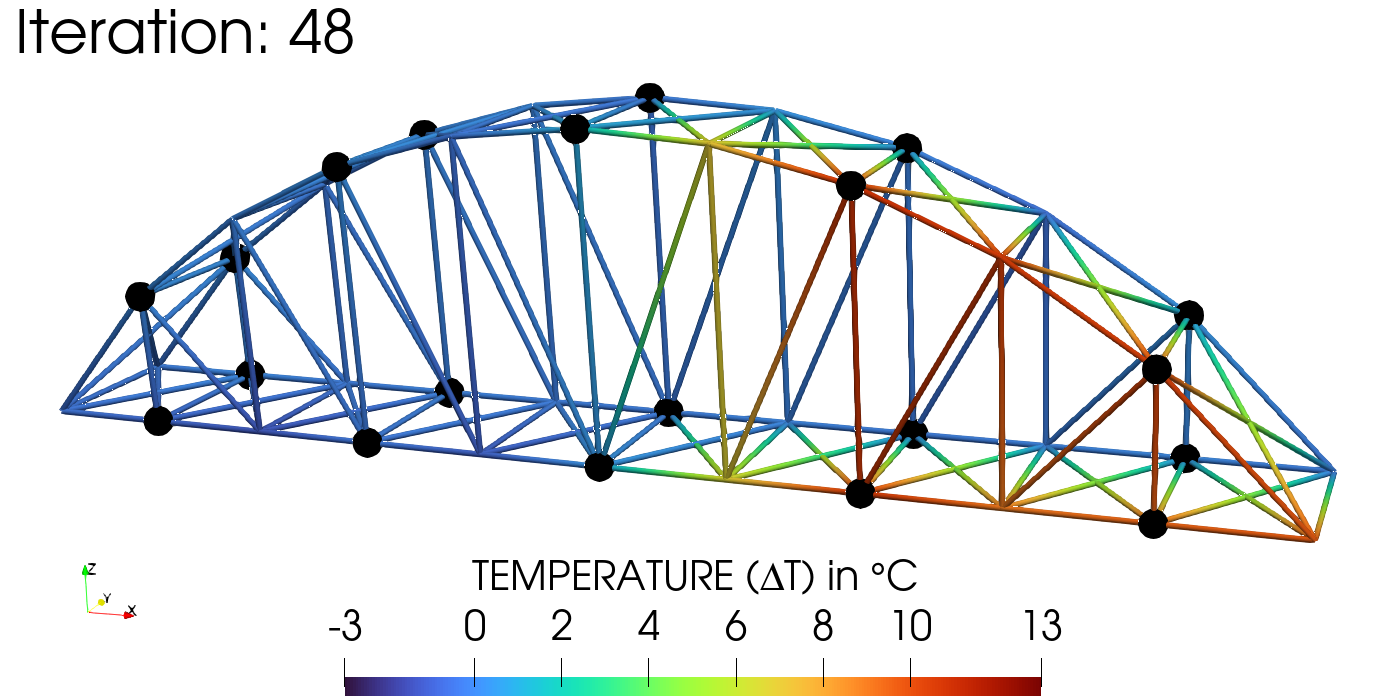}
        \end{minipage}
        \hfill
        \caption{20 sensors: Iteration 48}
        \label{fig:Figure_13b}
    \end{subfigure}
    \caption{Bridge example, 8 and 20 sensors configurations: Intermediate temperature distributions obtained at (approx.) 20\% of the optimization progress i.e., at Iteration = 20\% * (Iterations to converge) with Vertex Morphing ($radius = 6$) filtering.}
    \label{fig:Figure_13}
\end{figure}

To examine this behavior visually, the intermediate temperature distributions obtained at around 20\% of the optimization progress for the 8 and 20 sensors configurations using Vertex Morphing filtering are shown in Figures \ref{fig:Figure_13a} and \ref{fig:Figure_13b} respectively.
From the figures, it is evident that at 20\% progress, most of the nodes have roughly been correctly identified.
However, from Figure \ref{fig:Figure_13a}, it can be seen that the temperature at some nodes deviates largely from the prescribed temperature distribution. It reaches as low as $-3.03$ $\degree C$, while the upper bound exceeds $12$ $\degree C$. 
While the thermal field for the 20 sensors configuration at 20\% progress, shown in Figure \ref{fig:Figure_13b}, appears to be much better, it is also only roughly identified with temperatures in the range of $[-1.2,10.82]$ $\degree C$. As mentioned above, the remaining optimization iterations (approximately 80\% of the total iterations) are used to fine-tune the thermal field and achieve the last 1-2 magnitudes drop in the cost function.

\FloatBarrier
\subsection{Simplified Hoover Dam}
\label{sec:Hoover dam}

\begin{figure}[!t]
\begin{minipage}[c][][t]{.32\textwidth}
  \vspace*{\fill}
  \centering
   \includegraphics[trim= 0 0 0 0, clip, width=\textwidth]{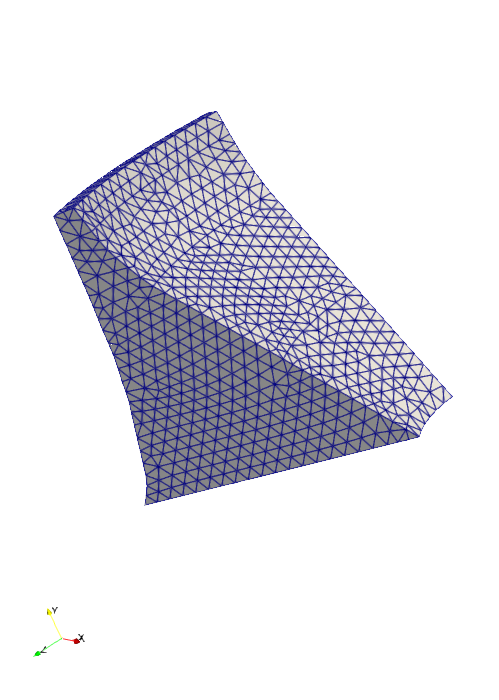}
\end{minipage}
\begin{minipage}[c][][t]{.67\textwidth}
  \vspace*{\fill}
  \centering
  \includegraphics[trim= 0 0 510 0, clip,width=\textwidth]{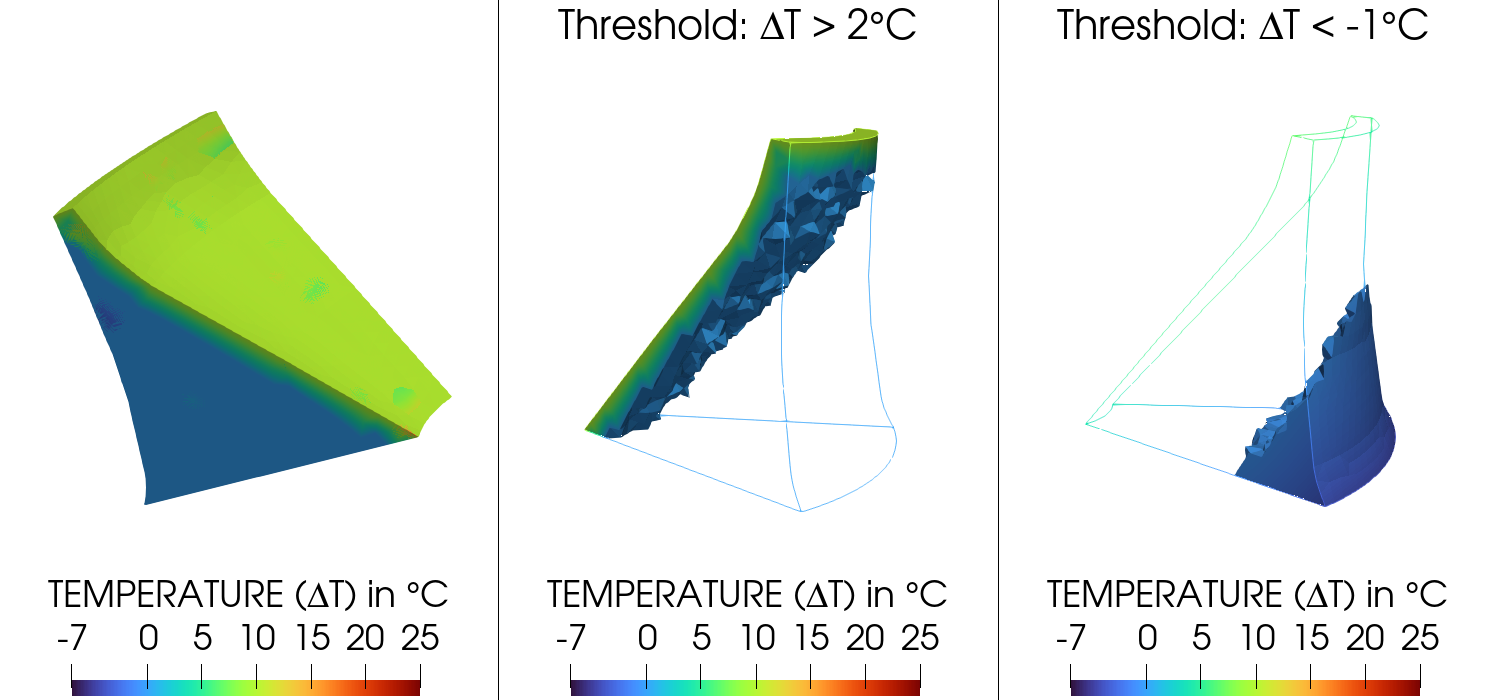} 
\end{minipage}
  \caption{Dam example: Mesh used (left), the target temperature distribution (middle), and a threshold distribution of $\Delta T > 2$ $\degree C$ (right).}
  \label{fig:Figure_14}
\end{figure}

A simplified Finite-Element model of the Hoover dam is shown in Figure \ref{fig:Figure_14}. 
Symmetry is assumed along the middle section, and thus, only half of the dam is simulated with appropriate boundary conditions to reflect symmetry. 
The model is aligned such that the dam spans from the right surface (viewed from downstream), coinciding with the $xy$ plane at $z=0$, to the left surface (viewed from downstream), coinciding with the vertical plane at $z=-x$ line.

The downstream bottom edge is at a radius of 65 $meters$ from the origin and at an elevation of 0 $meters$, the upstream bottom edge is at a radius of 266.863 $meters$ from the origin and at an elevation of 6.6612 $meters$. The elevation at the highest section of the dam is 219.24 $meters$. 

Density, Young's modulus, Poisson ratio, and thermal expansion coefficient were set to $\rho = 2400$, $E = 3\cdot 10^{10}$, $\nu = 0.15$, $\alpha = 1.0 \cdot 10^{-5} /\degree K$ respectively in SI units. 
The bottom surface and left surface of the dam were assumed to be fixed ($\mathbf{u} = \mathbf{0}$), while the right surface restricted movements only in $z$-direction due to the symmetry ($u_z = 0$). A hydrostatic pressure of 45,000 psf (Pounds per Square Foot) was applied at the upstream bottom edge, which linearly reduced to zero at the water level considered to be 8.2296 $meters$ (27 feet) below the top edge. Load due to self-weight and gravity of $9.81$ $m/s^2$ was also considered. 

\begin{figure}[!b]
  \centering
        \includegraphics[trim= 0 250 0 225, clip, width=0.2\paperwidth]{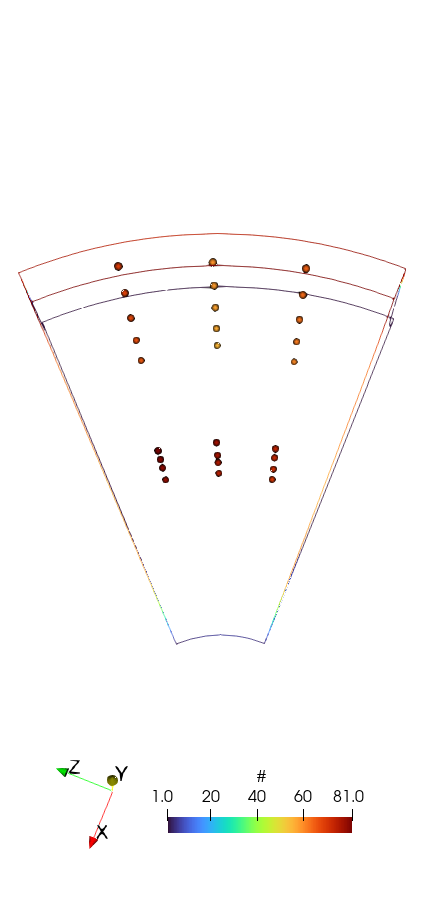} 
        \includegraphics[trim= 0 250 0 170, clip,width=0.2\paperwidth]{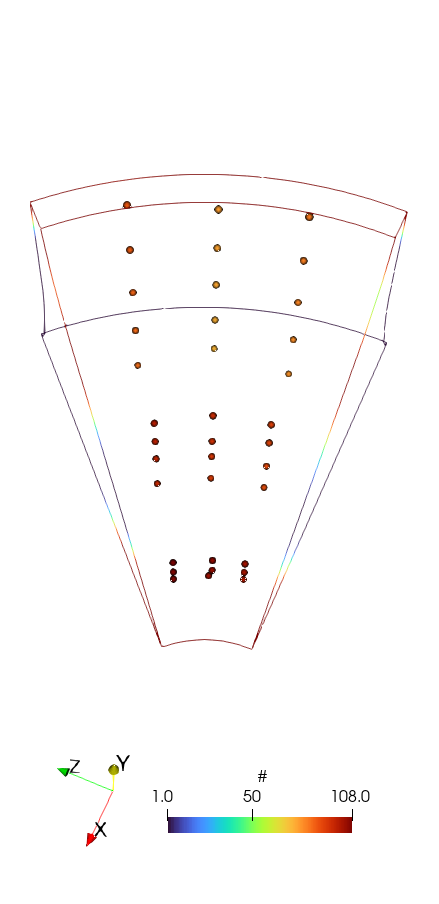} 
        \includegraphics[trim= 0 250 0 225, clip,width=0.2\paperwidth]{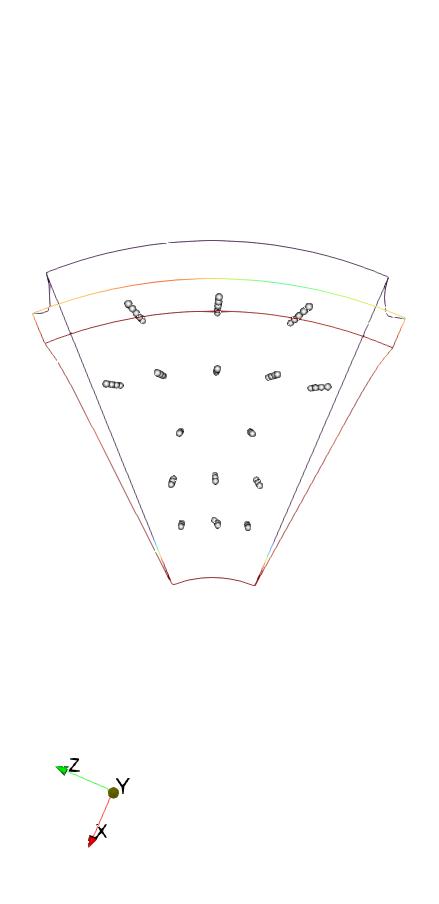} \hfill
  \caption{Dam example: Sensors distribution for the 27, 36, and 59 sensors configuration (left to right).}
  \label{fig:Figure_15}
\end{figure}

24,205 linear tetrahedral small displacement elements were used in the finite-element mesh. The temperature difference was prescribed to be $\Delta T  = 10$ $\degree C$ on the downstream surface and the top surface. Due to the slow heat transfer in concrete, it was assumed that the heat only penetrated a depth of 14 $meters$ from the irradiated surface and thus the temperature in this sub-surface zone was set to $\Delta T  = 5$ $\degree C$; everywhere else the temperature difference of $\Delta T = 0$ $\degree C$ was set. This distribution is shown in Figure \ref{fig:Figure_14} (middle) along with a threshold distribution of $\Delta T > 2$ $\degree C$ to better visualize the depth of the heated region. 
It should be noted that the threshold distributions shown in this paper display an element if at least one of its nodes satisfies the threshold condition. Due to this, the presented threshold distributions are not an exact representation of the threshold conditions but rather a close approximation as some nodes may be outside the threshold but still visible due to its connectivity. A true nodal threshold distribution is not shown as the visualization is rather difficult.

Three different sensor configurations with 27, 36, and 59 displacement sensors were tested. For the 27 sensors configuration, 12 sensors were placed at a radius of 155 $meters$ from origin on three vertical planes at $11.25\degree, 22.5\degree, $ and $33.75\degree$ counterclockwise (\acrshort{CCW}) to the right surface (viewed from downstream), with each plane consisting of 4 sensors at different elevations; the remaining 15 sensors were placed at a radius of 222 $meters$ from origin on the same three vertical planes, with each plane consisting of 5 sensors at different elevations.

For the 36 sensors configuration, the 27 sensors configuration was modified to add 9 additional sensors on the same vertical planes at a radius of 100 $meters$ from the origin, with each plane containing 3 sensors at different elevations.

The 59 sensors configuration comprised of sensors distributed along five different radii or rows. 9 sensors each were placed at radii of 100 $meters$ and 127 $meters$, and 15 sensors were placed at a radius of 222 $meters$, with each row equally divided into three vertical planes at $11.25\degree, 22.5\degree, $ and $33.75\degree$ CCW to the right surface. Another 20 sensors were placed at a radius of 188.5 $meters$ from origin on five vertical planes at $5\degree, 13\degree, 22.5\degree,32\degree,$ and $40\degree$ CCW to the right surface (viewed from downstream), with each plane consisting of 4 sensors at different elevations; the remaining 6 sensors were placed at a radius of 155 $meters$ from origin on the two vertical planes at $15\degree,$ and $30\degree$ CCW to the right surface, with each plane consisting of 3 sensors at different elevations.
The three sensor configurations viewed from above are shown in Figure \ref{fig:Figure_15}.

\begin{figure}[!t]
  \centering
        \includegraphics[trim= 0 0 0 0, clip, width=0.25\paperwidth]{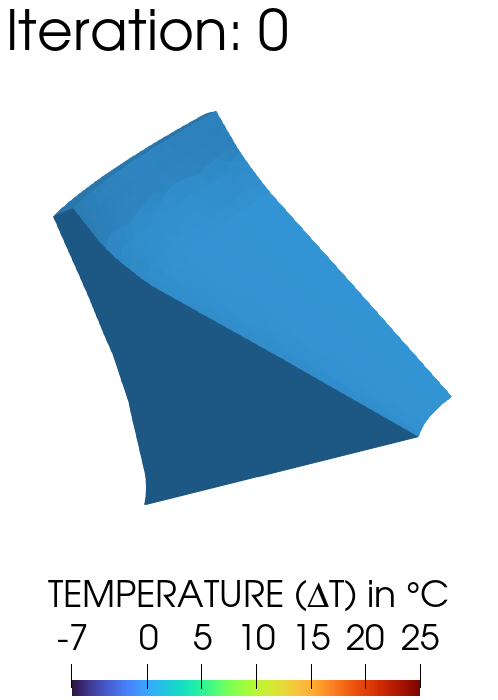} \hfill
  \caption{Dam example: Temperature distribution at optimization start $\Delta \mathbf{T} = \mathbf{0}$.}
  \label{fig:Figure_16}
\end{figure}

\begin{figure}[!b]
    \centering
    \begin{subfigure}[t]{\textwidth}
        \centering
        \begin{minipage}[t]{0.62\textwidth}
            \centering
            \includegraphics[trim=0 0 0 0, clip, width=\textwidth]{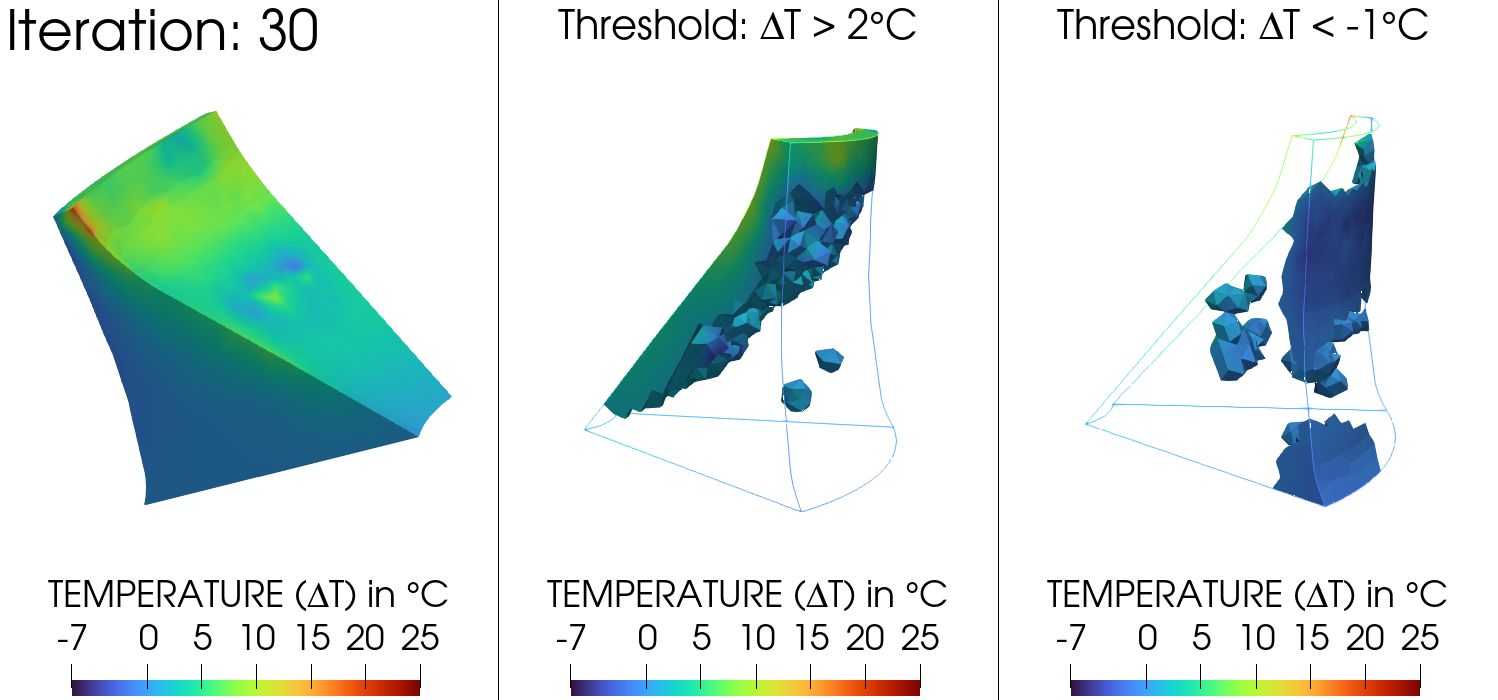}
        \end{minipage}
        \hfill
        \begin{minipage}[t]{0.37\textwidth}
            \centering
            \includegraphics[trim=0 0 0 0, clip, width=\textwidth]{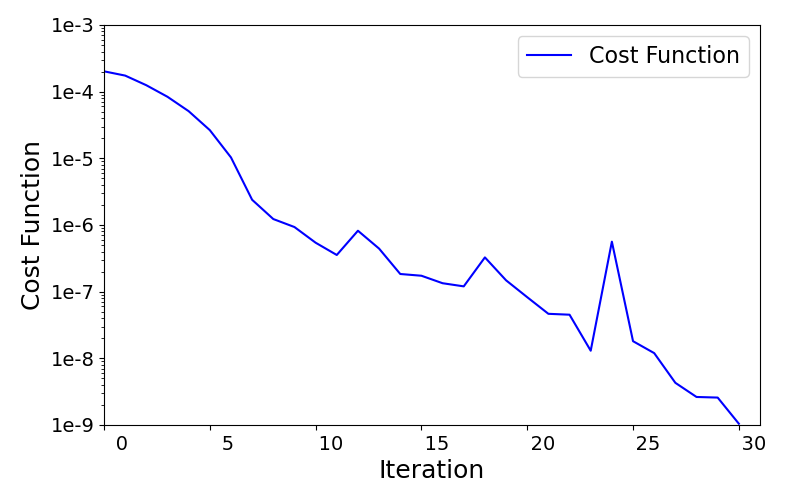}
        \end{minipage}
        \caption{Without Vertex Morphing filtering: Iteration 30}
        \label{fig:Figure_17a}
        \vspace{1em}
    \end{subfigure}
    \vspace{1em}
    \begin{subfigure}[t]{\textwidth}
        \centering
        \begin{minipage}[t]{0.62\textwidth}
            \centering
            \includegraphics[trim=0 0 0 0, clip, width=\textwidth]{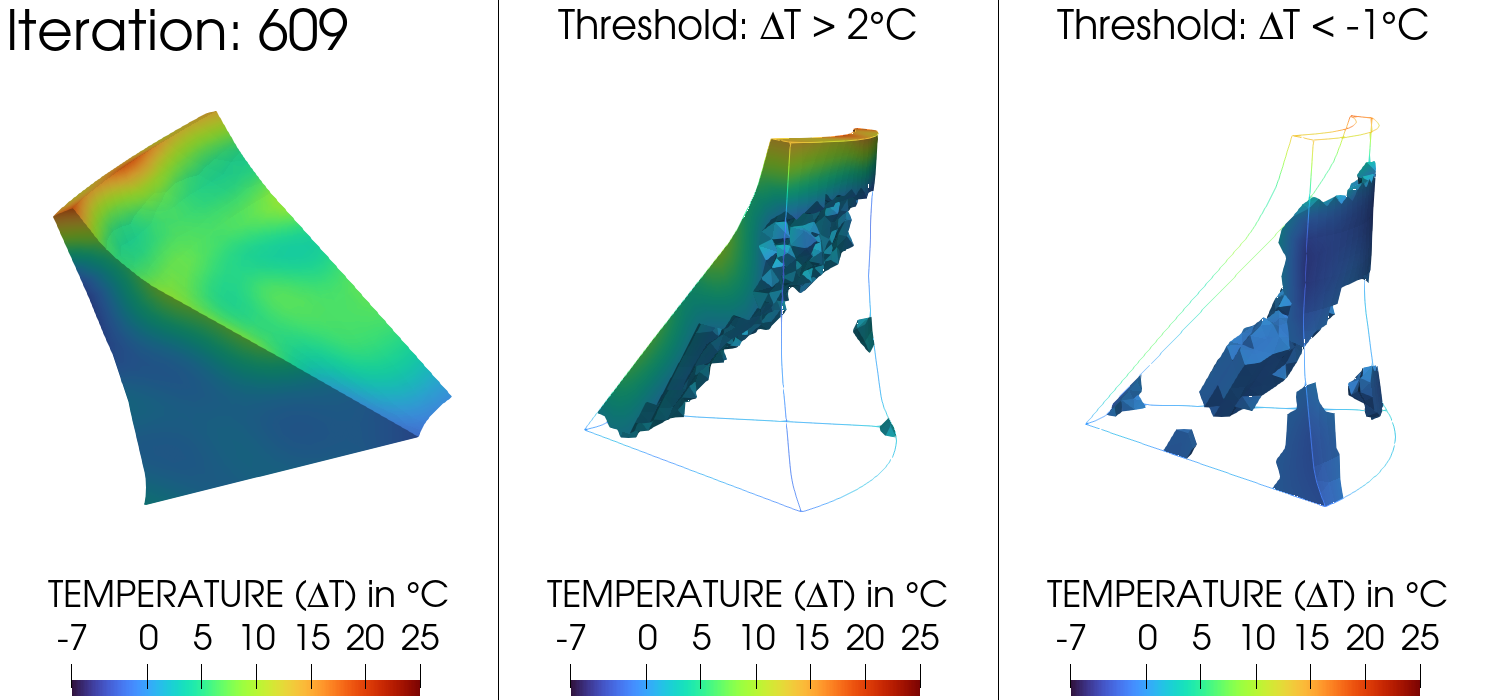}
        \end{minipage}
        \hfill
        \begin{minipage}[t]{0.37\textwidth}
            \centering
            \includegraphics[trim=0 0 0 0, clip, width=\textwidth]{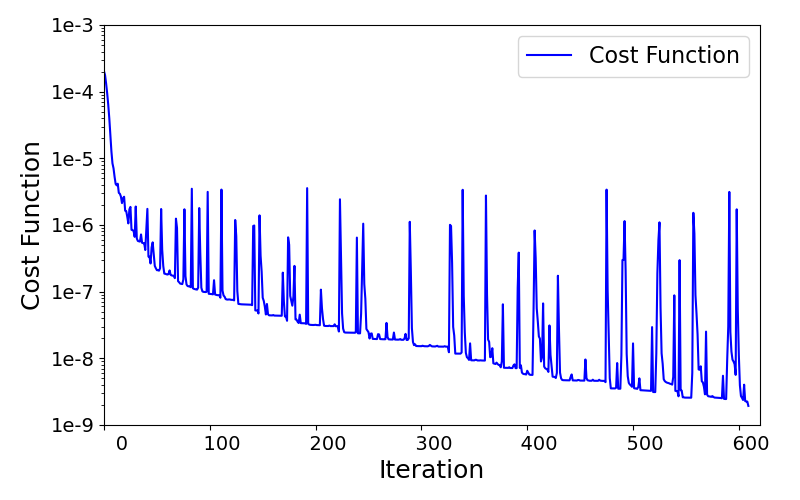}
        \end{minipage}
        \caption{With Vertex Morphing filtering (radius = $50$): Iteration 609}
        \label{fig:Figure_17b}
    \end{subfigure}
    \caption{Dam example, 27 sensors configuration: Temperature distributions obtained (leftmost) along with the $\Delta T > 2$ $\degree C$ thresholds, the $\Delta T < -1$ $\degree C$ thresholds, and the cost function convergence plots (right) without and with Vertex Morphing filtering.}
    \label{fig:Figure_17}
\end{figure}

For all three sensor configurations, the steepest descent optimization algorithm with the Barzilai-Borwein stepsize method (Eqn.\eqref{eq:bb}) was used. A maximum stepsize of $1\cdot 10^{-1}$ was set. The optimization starts with the presumption that the $\Delta T = 0$ $\degree C$ everywhere on the dam. The temperature distribution at the start of the optimization is shown in Figure \ref{fig:Figure_16}.

The same stopping criteria as for the Plate With a Hole and Bridge examples was set: 5 magnitudes reduction of the initial cost function i.e., the algorithm stops when the $\text{current cost function} \leq (1\cdot 10^{-5} * \text{initial cost function})$. 
Vertex morphing with radius $r = 50$ was used to smoothen the gradients and subsequently get smoother results.
For better visualization of the results, the color bar is shown in the range $[-7,25]$ $\degree C$ as the temperature ($\Delta T$) in no case falls below $-7$ $\degree C$, neither does it exceed $25$ $\degree C$. 

The temperature distributions obtained at the final iterations, threshold distributions of $\Delta T > 2$ $\degree C$ and $\Delta T < -1$ $\degree C$, and the cost function convergence plots for the 27 sensors configuration, 'without' and 'with' Vertex Morphing filtering are shown in Figures \ref{fig:Figure_17a} and \ref{fig:Figure_17b} respectively. 
In general, it can observed that due to a lack of sensors in the lower downstream section of the dam, the temperature identification is weak in that region. Although Vertex Morphing produces a smoother filtered field, the lowest downstream section of the dam is still incorrectly identified. This is no surprise because that region is devoid of sensors, hence the sensitivity in that region is very low. 

The $\Delta T > 2$ $\degree C$ threshold distributions indicate that a good localization in terms of depth of the heat penetration tangential to the downstream surface has been achieved. On the other hand, the $\Delta T < -1$ $\degree C$ threshold distributions show some negative $\Delta T$'s, particularly in regions of low sensitivities usually at corners and edges. These are regions where the change in temperature has a negligible effect on the responses of the sensors. Hence, these regions act freely in the sense that their value will not affect the cost function. This is a strong argument for using a larger number of sensors placed at optimal locations and responses recorded at different loading conditions such that such 'blind spots' of negligible sensitivities are reduced.

From the cost function convergence plots, it can be seen that the case with Vertex Morphing requires a lot more iterations to converge than the one using no filtering, but visually, it is evident that the temperature profile obtained using Vertex Morphing is closer to the prescribed temperature field than the no filtering case.

\medskip

\begin{figure}[!b]
    \centering
    \begin{subfigure}[t]{\textwidth}
        \centering
        \begin{minipage}[t]{0.62\textwidth}
            \centering
            \includegraphics[trim=0 0 0 0, clip, width=\textwidth]{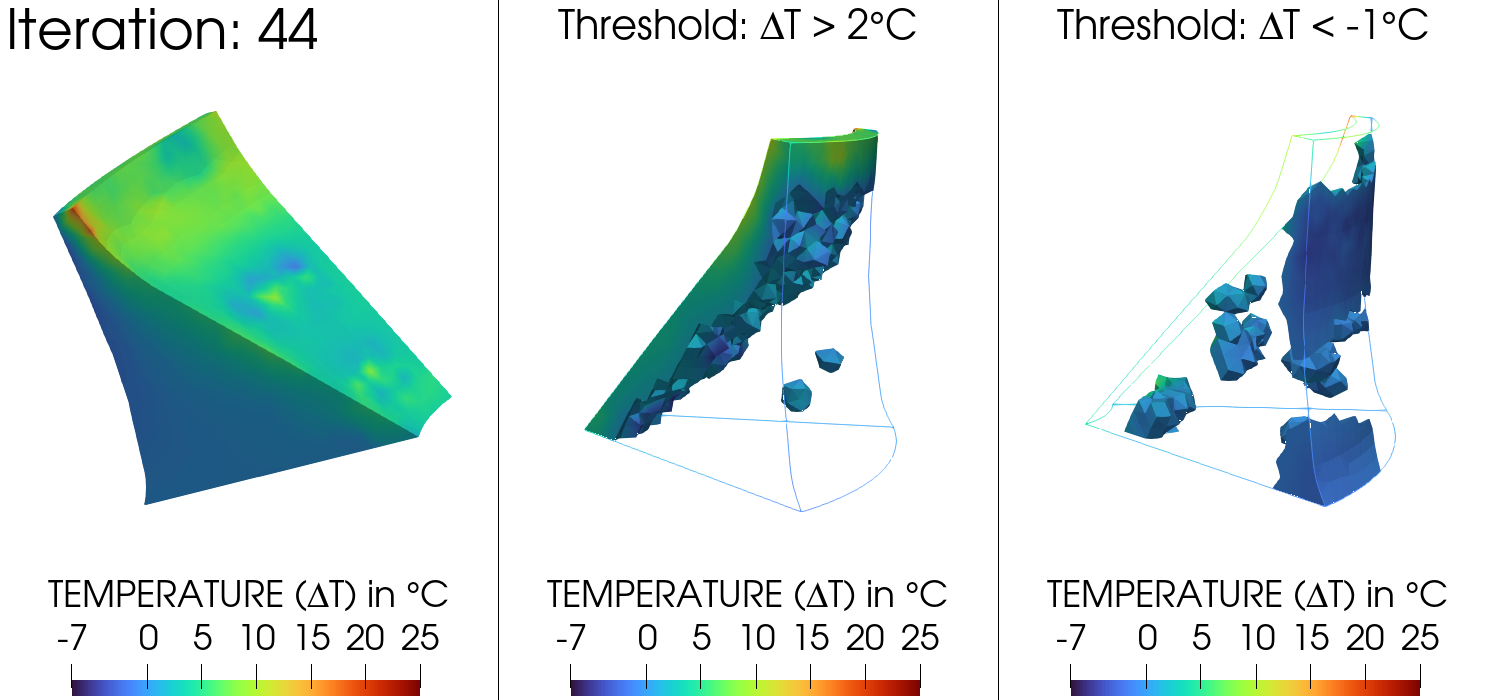}
        \end{minipage}
        \hfill
        \begin{minipage}[t]{0.37\textwidth}
            \centering
            \includegraphics[trim=0 0 0 0, clip, width=\textwidth]{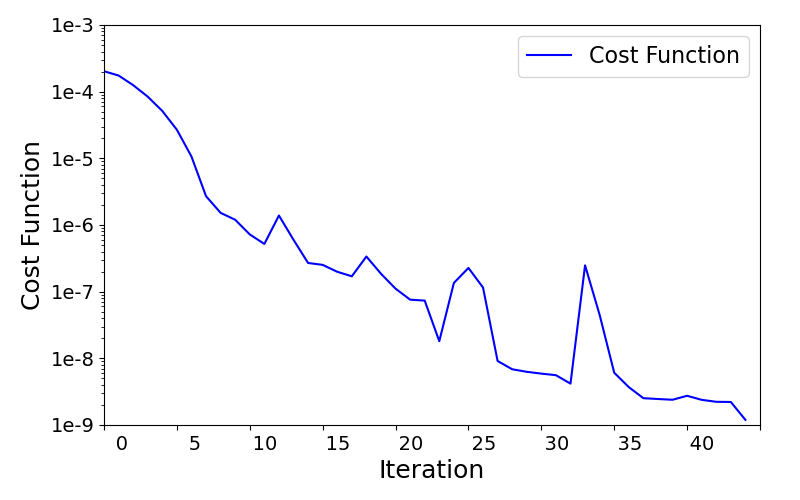}
        \end{minipage}
        \caption{Without Vertex Morphing filtering: Iteration 44}
        \label{fig:Figure_18a}
        \vspace{1em}
    \end{subfigure}
    \vspace{1em}
    \begin{subfigure}[t]{\textwidth}
        \centering
        \begin{minipage}[t]{0.62\textwidth}
            \centering
            \includegraphics[trim=0 0 0 0, clip, width=\textwidth]{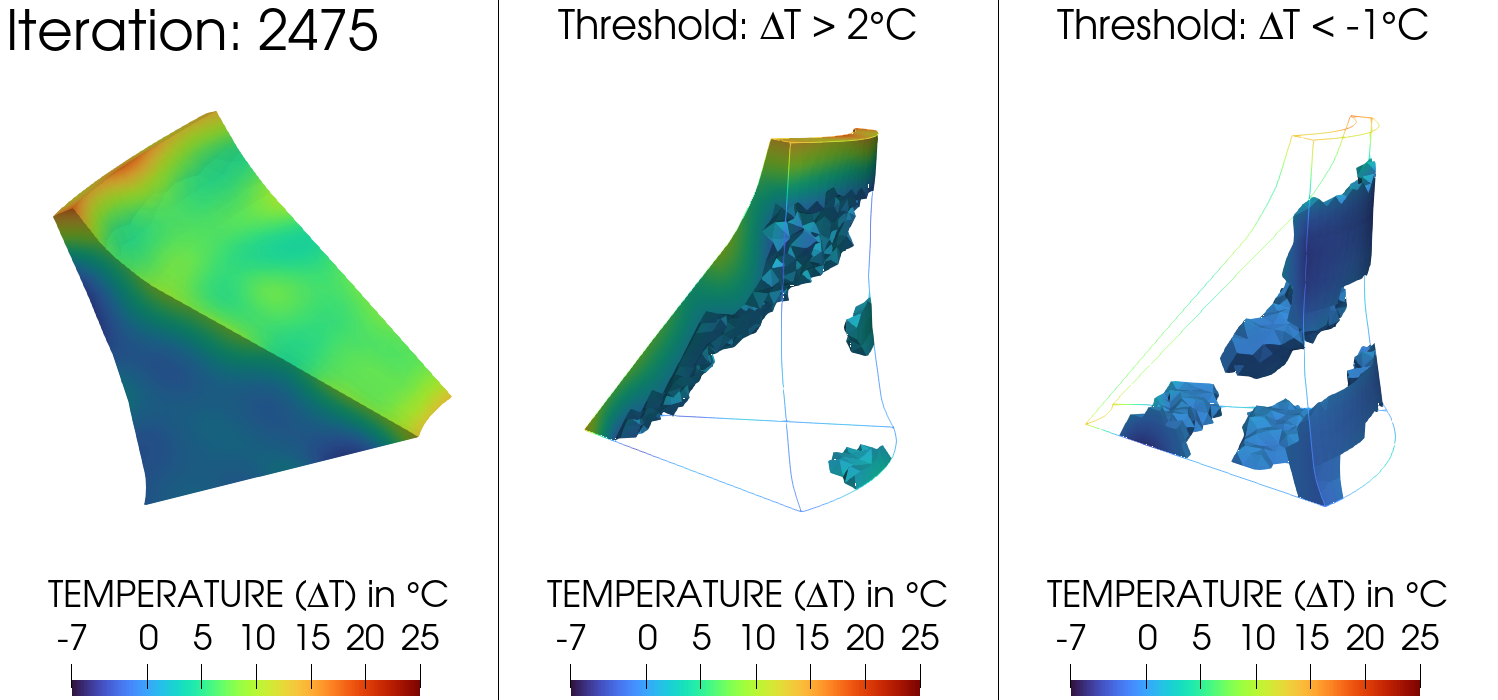}
        \end{minipage}
        \hfill
        \begin{minipage}[t]{0.37\textwidth}
            \centering
            \includegraphics[trim=0 0 0 0, clip, width=\textwidth]{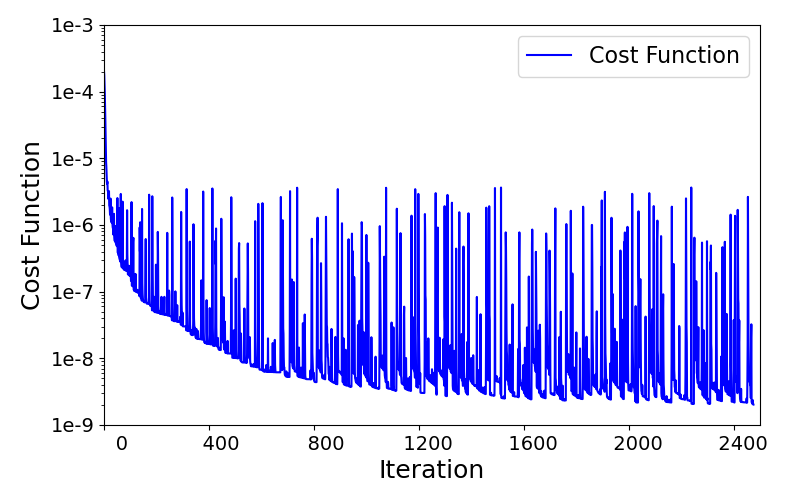}
        \end{minipage}
        \caption{With Vertex Morphing filtering (radius = $50$): Iteration 2475}
        \label{fig:Figure_18b}
    \end{subfigure}
    \caption{Dam example, 36 sensors configuration: Temperature distributions obtained (leftmost) along with the $\Delta T > 2$ $\degree C$ thresholds, the $\Delta T < -1$ $\degree C$ thresholds, and the cost function convergence plots (right) without and with Vertex Morphing filtering.}
    \label{fig:Figure_18}
\end{figure}

The temperature distributions obtained at the final iterations, threshold distributions of $\Delta T > 2$ $\degree C$ and $\Delta T < -1$ $\degree C$, and the cost function convergence plots for the 36 sensors configuration, 'without' and 'with' Vertex Morphing filtering are shown in Figures \ref{fig:Figure_18a} and \ref{fig:Figure_18b} respectively. 
The addition of 9 sensors, i.e., 27 to 36 sensors, does not cause a significant increase in the number of iterations for convergence in the case 'without' Vertex Morphing as seen in Figure \ref{fig:Figure_18a}, but the number of iterations required to converge with Vertex Morphing has almost quadrupled. 
Nevertheless, it can immediately be seen that the 9 additional sensors located in the lower downstream region of the dam in the 36 sensors configuration have a significant improvement in the quality of the thermal field identification. 

Without the use of Vertex Morphing, areas close to sensors experience quick and sharp temperature ($\Delta T$) rise, thereby satisfying the convergence criteria swiftly in 44 iterations. In contrast, Vertex Morphing produces smoother design ($\Delta T$) updates resulting in a smoother thermal gradient in the structure. Since the sharp gradients get smoothed out over a region, the optimization consumes more iterations to achieve convergence.   

Upon close inspection of Figure \ref{fig:Figure_18b}, a 'hotspot' pattern related to the sensor locations can be observed in the temperature distribution on the downstream surface. This is an indication towards the region of observation of a sensor which should be studied for optimal sensor placement.
The identification of the depth of heat penetration and large positive and artificial artifacts causing negative $\Delta T$'s are similar to the 27 sensors configuration.

\begin{figure}[!b]
    \centering
    \begin{subfigure}[t]{\textwidth}
        \centering
        \begin{minipage}[t]{0.62\textwidth}
            \centering
            \includegraphics[trim=0 0 0 0, clip, width=\textwidth]{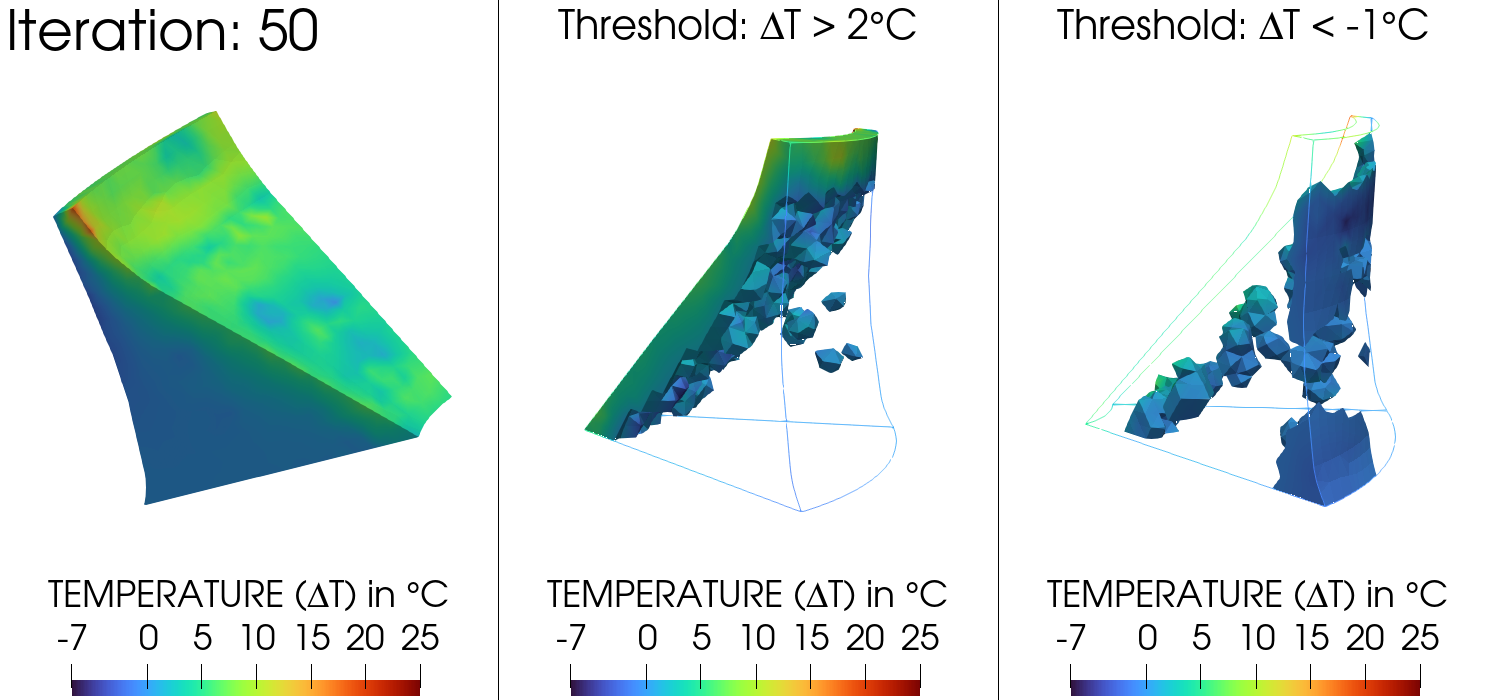}
        \end{minipage}
        \hfill
        \begin{minipage}[t]{0.37\textwidth}
            \centering
            \includegraphics[trim=0 0 0 0, clip, width=\textwidth]{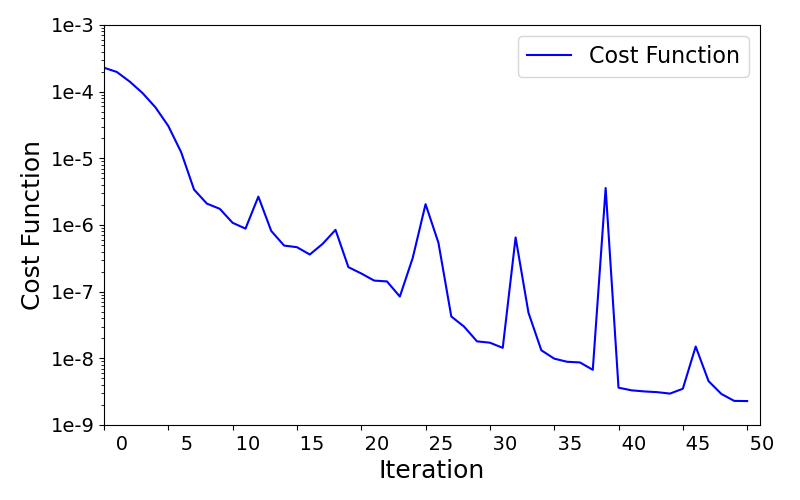}
        \end{minipage}
        \caption{Without Vertex Morphing filtering: Iteration 50}
        \label{fig:Figure_19a}
        \vspace{1em}
    \end{subfigure}
    \vspace{1em}
    \begin{subfigure}[t]{\textwidth}
        \centering
        \begin{minipage}[t]{0.62\textwidth}
            \centering
            \includegraphics[trim=0 0 0 0, clip, width=\textwidth]{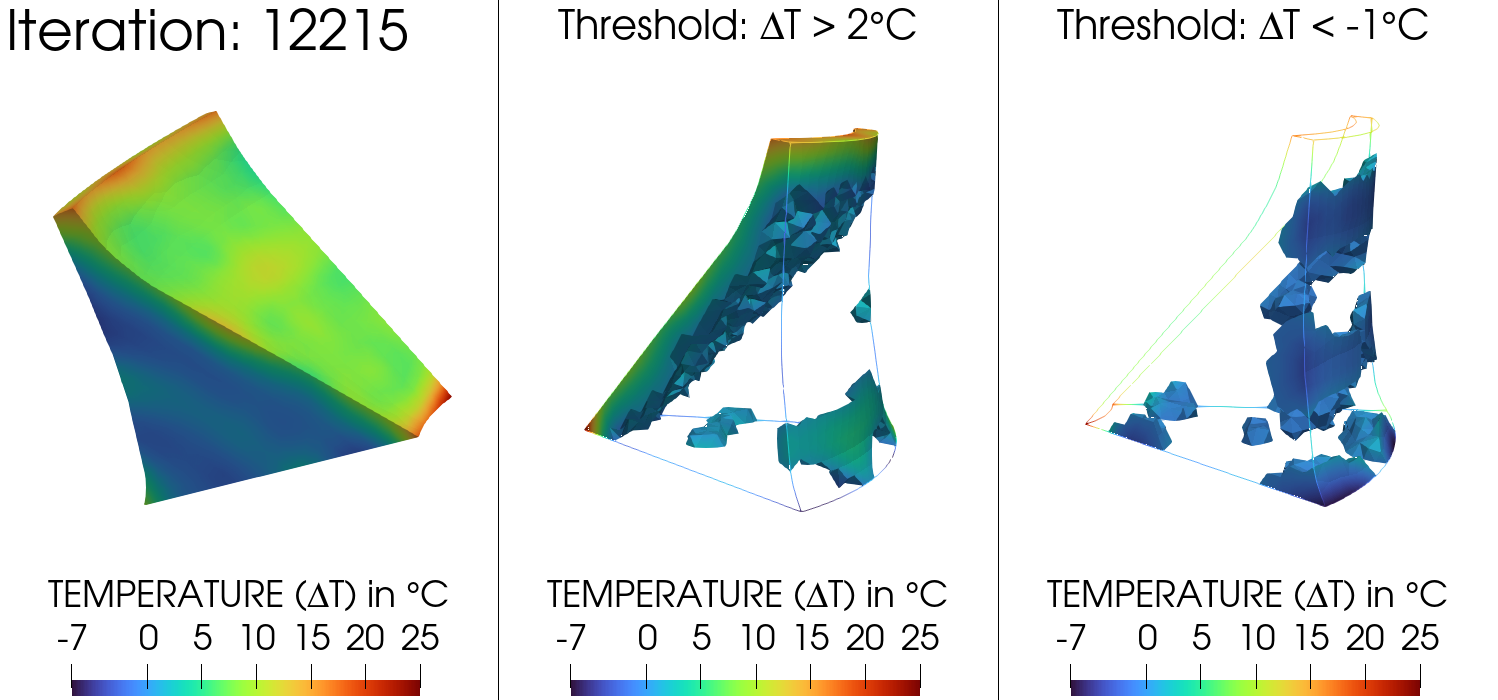}
        \end{minipage}
        \hfill
        \begin{minipage}[t]{0.37\textwidth}
            \centering
            \includegraphics[trim=0 0 0 0, clip, width=\textwidth]{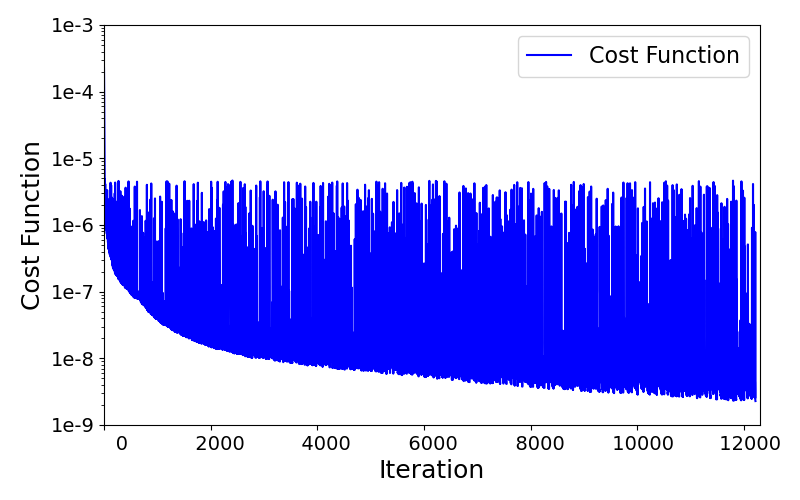}
        \end{minipage}
        \caption{With Vertex Morphing filtering (radius = $50$): Iteration 12215}
        \label{fig:Figure_19b}
    \end{subfigure}
    \caption{Dam example, 59 sensors configuration: Temperature distributions obtained (leftmost) along with the $\Delta T > 2$ $\degree C$ thresholds, the $\Delta T < -1$ $\degree C$ thresholds, and the cost function convergence plots (right) without and with Vertex Morphing filtering.}
    \label{fig:Figure_19}
\end{figure}

\medskip

The temperature distributions obtained at the final iterations, threshold distributions of $\Delta T > 2$ $\degree C$ and $\Delta T < -1$ $\degree C$, and the cost function convergence plots for the 59 sensors configuration, 'without' and 'with' Vertex Morphing filtering are shown in Figures \ref{fig:Figure_19a} and \ref{fig:Figure_19b} respectively. 
Due to the distribution of 59 sensors in five different radii, the temperature field identification is even more accurate and very close to the target distribution. Nevertheless, there are deficiencies. 
The case without Vertex Morphing shown in Figure \ref{fig:Figure_19a} converges quickly but has a few small regions where the temperature identified is grossly wrong. These are the regions that do not fall in the area of observation of the sensors, even with the presence of 59 sensors. The case with Vertex Morphing smoothens these regions based on the temperature identification of the neighboring area. 

The 59 sensors configuration with Vertex Morphing (Figure \ref{fig:Figure_19b}) is the only case where the temperature reaches around $24$ $\degree C$ compared to around $18$ $\degree C$ in all other cases. But as it can be seen from the figure, this extreme temperature is an anomaly that occurs at the downstream bottom corner of the dam, which is a low-sensitivity region, hence allowing the temperature to go high without much effect observed by the sensors. 

The localization quality using the 59 sensors configuration is also superior to that of the 27 and 36 sensors configurations. This is seen by the depth of the $\Delta T > 2$ $\degree C$ threshold distribution on the downstream surface where few elements (and nodes) appear to satisfy the threshold as expected from the prescribed temperature distribution in Figure \ref{fig:Figure_14}. 

The number of iterations required to converge increases many folds when using Vertex Morphing but it can be reduced by using advanced optimization algorithms with better stepsize calculation. 

Due to the 3-D temperature distribution, the temperature profile has different scales in different directions (cf. Figure \ref{fig:Figure_14}). The inclined vertical distance from the downstream surface bottom to the top is approximately 265 $meters$, the horizontal length of the downstream surface varies from around 50 $meters$ at the bottom to around 165 $meters$ at the top edge, and the depth or tangential direction to the downstream surface is up to 170 $meters$ with the heat penetrating only the initial 14 $meters$. Due to these different scales in the thermal field to be identified, a large Vertex Morphing filter radius (radius = 50) was used. However, it is noteworthy that using such a large filtering radius did not obscure the local feature identification. This is seen in the tangential direction localization of the heat penetration, which has a smaller scale (14 $meters$ i.e., approximately 28\% of the filter radius) (cf. the $\Delta T>2$ $\degree C$ threshold distributions in Figures \ref{fig:Figure_17}, \ref{fig:Figure_18}, and \ref{fig:Figure_19}).

Ignoring any anomalous outliers causing very localized temperature peaks, in general, the temperature identified on the downstream surface is between $4-10$ $\degree C$ for the 27 and 36 sensors configurations and between $6-11$ $\degree C$ for the 59 sensors configuration. The 27 sensors case also has low temperatures between $0$ $\degree C$ and $-1$ $\degree C$ due to a sensor deficit in that area.

Similar to the Bridge example presented above, in the dam case also, the neighborhood of the heated region is identified rather quickly within 20\% of the total iterations to converge. This can be observed from the cost function convergence plots in Figures \ref{fig:Figure_17b}, \ref{fig:Figure_18b}, and \ref{fig:Figure_19b}, where approximately within the initial 20\% of the optimization progress the cost function has already dropped $3-4$ magnitudes compared to the initial value. In the initial 20\% of the optimization progress, the cost function has a steep decrease, after which the cost function convergence plot starts to flatten, and the improvement in the solution reduces. This slow improvement is seen in the last 1-2 magnitudes drop in the cost function, which consumes approximately 80\% of the total iterations required to converge. 

\begin{figure}[!b]
\begin{minipage}[c][][t]{.327\textwidth}
  \vspace*{\fill}
  \centering
  \includegraphics[trim= 0 0 998 0, clip, width=0.8\textwidth]{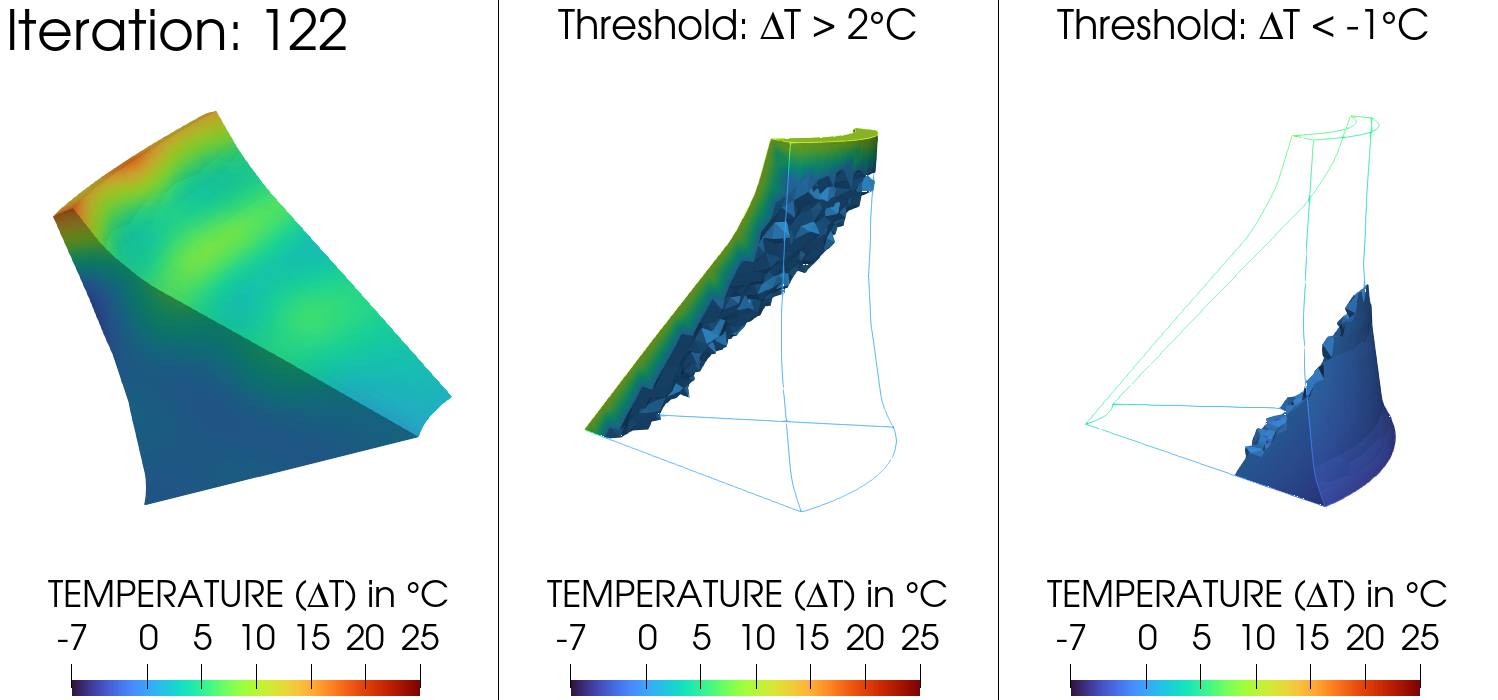}
  \subcaption{27 Sensors: Iteration 122}
   \label{fig:Figure_20a}
\end{minipage}
\begin{minipage}[c][][t]{.327\textwidth}
  \vspace*{\fill}
  \centering
  \includegraphics[trim= 0 0 998 0, clip, width=0.8\textwidth]{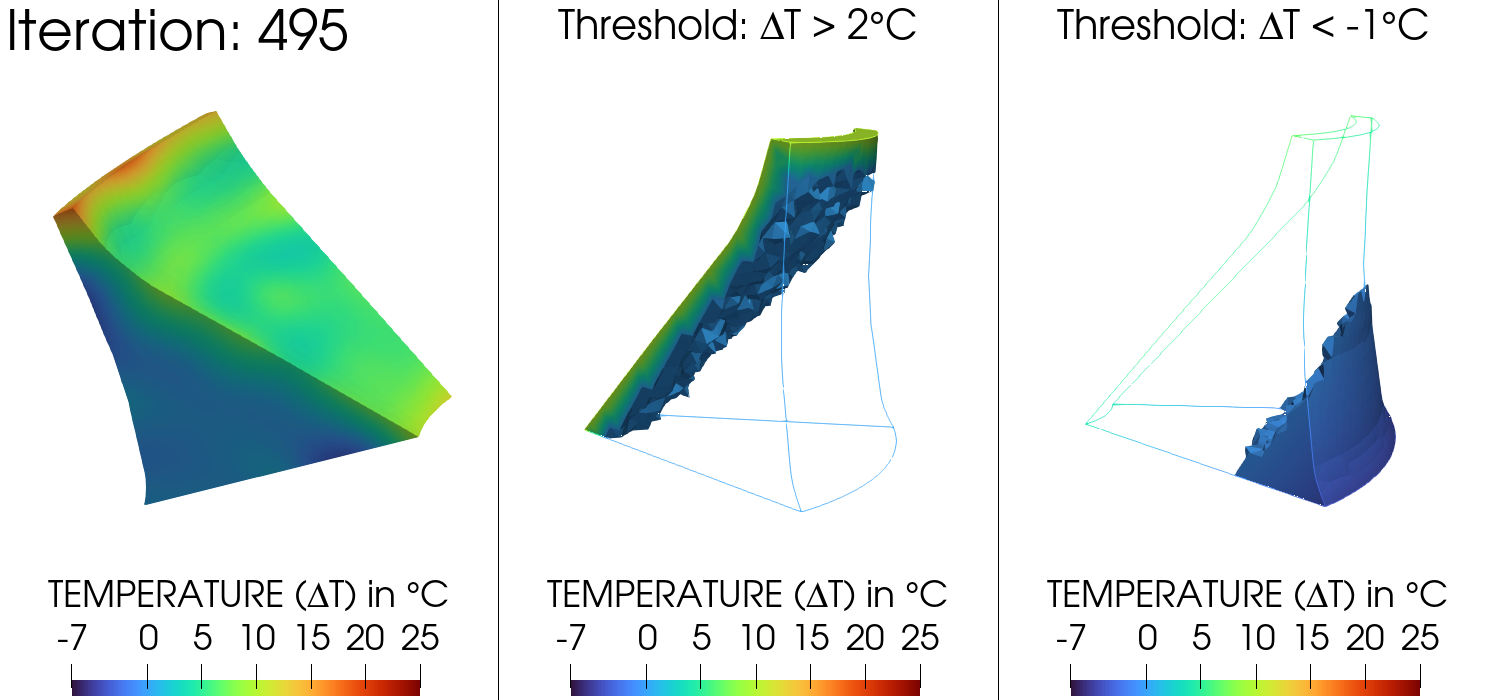}
  \subcaption{36 Sensors: Iteration 495}
  \label{fig:Figure_20b}\par
\end{minipage}
\begin{minipage}[c][][t]{.327\textwidth}
  \vspace*{\fill}
  \centering
  \includegraphics[trim= 0 0 1002 0, clip, width=0.8\textwidth]{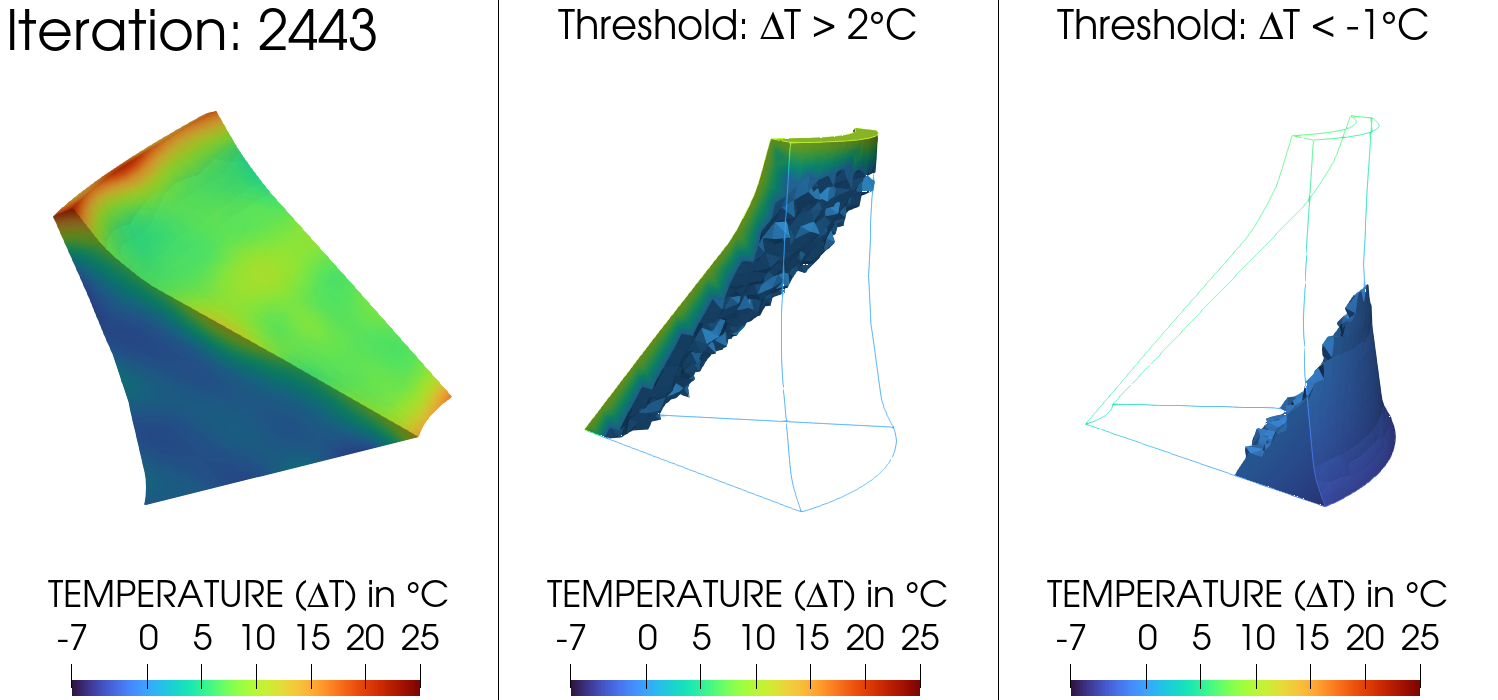}
  \subcaption{59 Sensors: Iteration 2443}
   \label{fig:Figure_20c}
\end{minipage}
\caption{Dam example, 27, 36, and 59 sensors configurations: Intermediate temperature distributions obtained at (approx.) 20\% of the optimization progress i.e., at Iteration = 20\% * (Iterations to converge) with Vertex Morphing filtering.}
\label{fig:Figure_20}
\end{figure}

To visually inspect this behavior, the intermediate temperature distributions obtained at around 20\% of the optimization progress for the 27, 36, and 59 sensors configurations using Vertex Morphing filtering are shown in Figures \ref{fig:Figure_20a}, \ref{fig:Figure_20b}, and \ref{fig:Figure_20c} respectively.
Comparing the thermal fields at 20\% optimization progress with the converged thermal fields for each of the three sensor configurations, the above observation supports that most of the progress is achieved within the initial 20\% of the optimization. This behavior will obviously vary for different cases and optimization algorithms, but it attests to the general optimization converge behavior, and the user can define the balance between accuracy and computational expense.

It is clear that the higher the number of sensors, the better the reconstruction of the thermal field. However, since sensors cannot be placed everywhere in the structure, there is always going to be some level of error that is present when using a limited number of sensors. Nevertheless, the results obtained are a close approximation of the true temperature distribution.

\FloatBarrier
\subsection{Comparison with Spatial Interpolation}
\label{sec:spatial interpolation}

To observe the benefit of the proposed approach, the temperature distribution was obtained using spatial interpolation techniques with the same sensor configurations (number and location) but using temperature sensors instead of displacement sensors. k-Nearest Neighbor and Kriging (Ordinary and Universal) interpolations were used for comparison here.

\subsubsection{k-Nearest Neighbors (\acrshort{kNN}) interpolation}

k-Nearest Neighbor uses \textit{k} nearest neighbors from the training dataset to interpolate the quantity at the point of interest. The \textit{k} neighbors can have equal contributions to the point of interest, i.e., uniform weighting, or the contribution of each neighbor can be proportional to the inverse of the distance between the neighbor and the point of interest, i.e., inverse distance weighting; or the weighting method can be user-defined.
A number of neighbors \textit{k} and the weighting method are hyperparameters that can be tuned based on the case.

\begin{figure}[!h]
\begin{minipage}[c][][t]{.5\textwidth}
  \vspace*{\fill}
  \centering
  \includegraphics[trim= 0 0 0 0, clip, width=0.3\paperwidth]{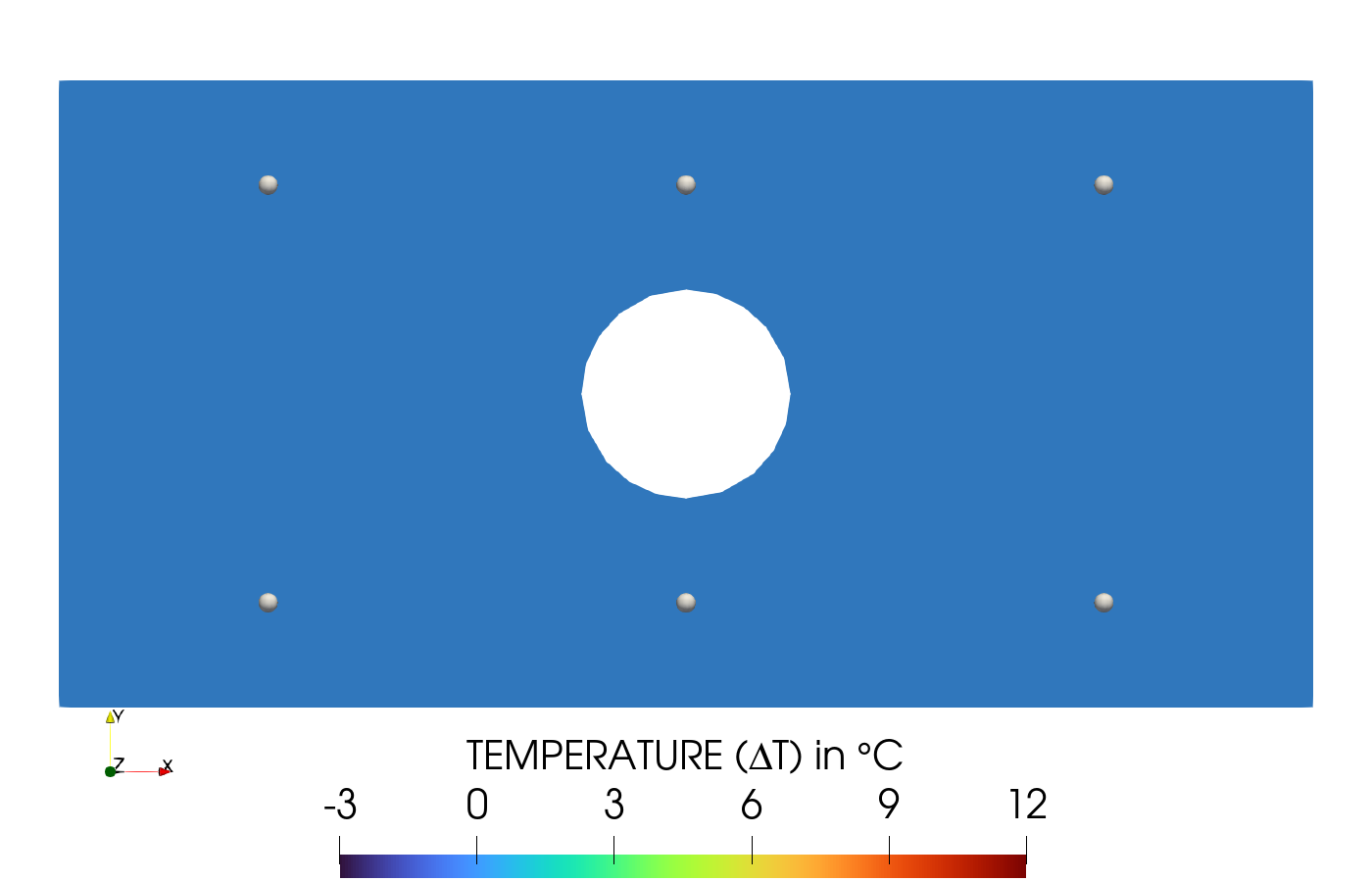}
  \subcaption{kNN: 6 temperature sensors configuration}
   \label{fig:Figure_21a}
\end{minipage}
\begin{minipage}[c][][t]{.5\textwidth}
  \vspace*{\fill}
  \centering
  \includegraphics[trim= 0 0 0 0, clip, width=0.3\paperwidth]{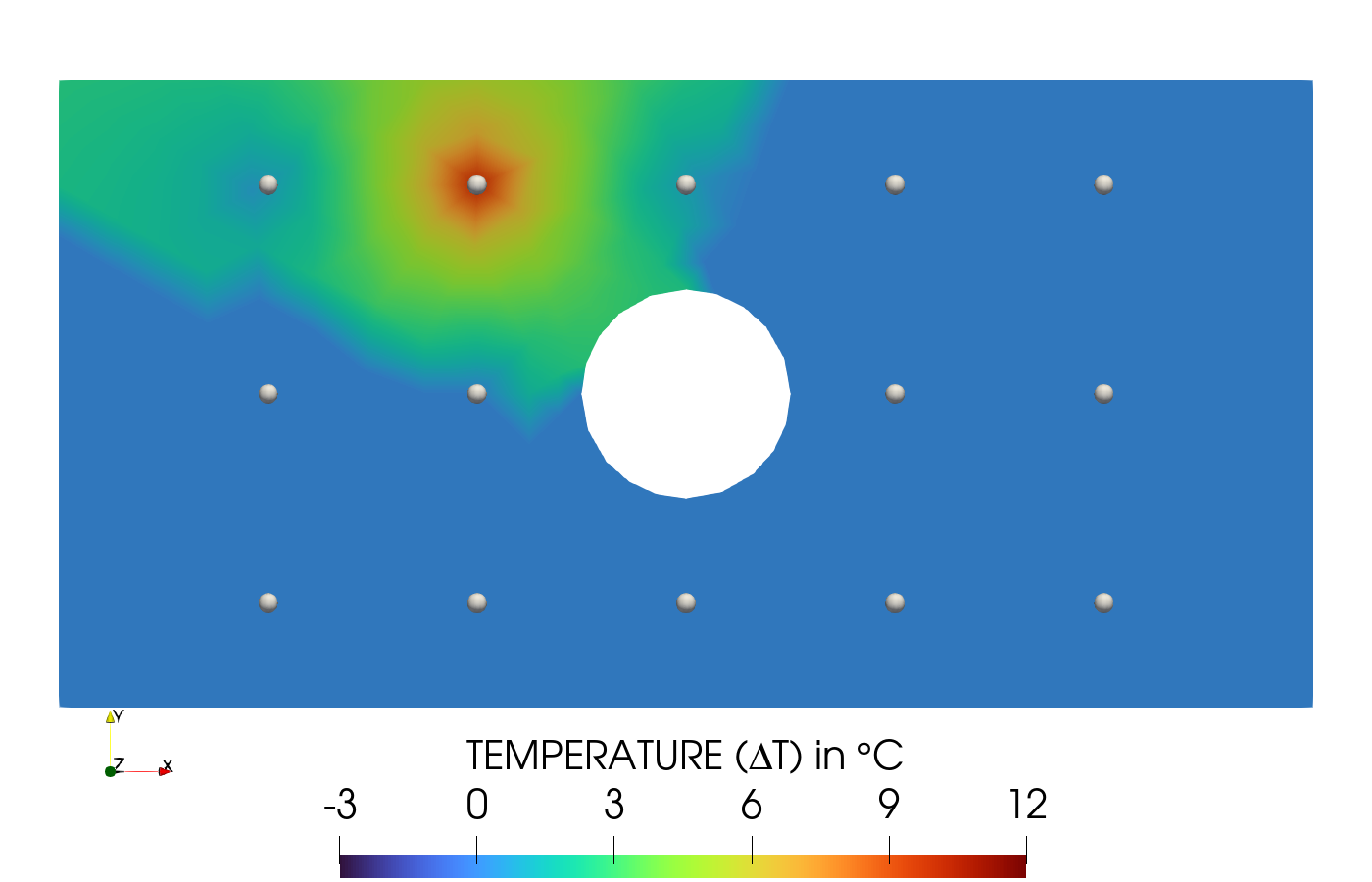}
  \subcaption{kNN: 14 temperature sensors configuration}
  \label{fig:Figure_21b}\par
\end{minipage}
\caption{Plate With a Hole, 6 and 14 temperature sensors configurations: k-Nearest Neighbor interpolation (kNN) with $\textit{k}=3$ and inverse distance weighting.}
\label{fig:Figure_21}
\end{figure}

\begin{figure}[!h]
    \centering
    \begin{subfigure}[t]{\textwidth}
        \centering
        \begin{minipage}[t]{0.62\textwidth}
            \centering
            \includegraphics[trim=0 0 0 0, clip, width=\textwidth]{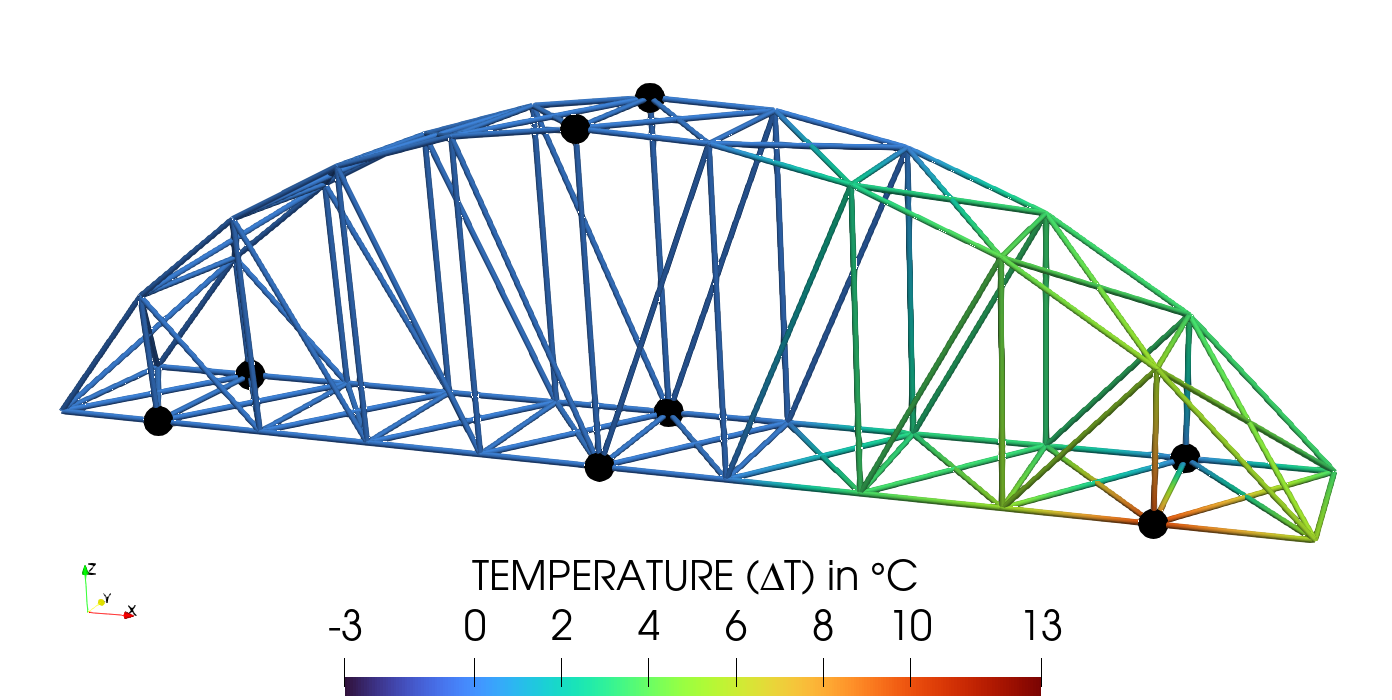}
        \end{minipage}
        \hfill
        \caption{kNN: 8 temperature sensors configuration}
        \label{fig:Figure_22a}
    \end{subfigure}
    \begin{subfigure}[t]{\textwidth}
        \centering
        \begin{minipage}[t]{0.62\textwidth}
            \centering
            \includegraphics[trim=0 0 0 0, clip, width=\textwidth]{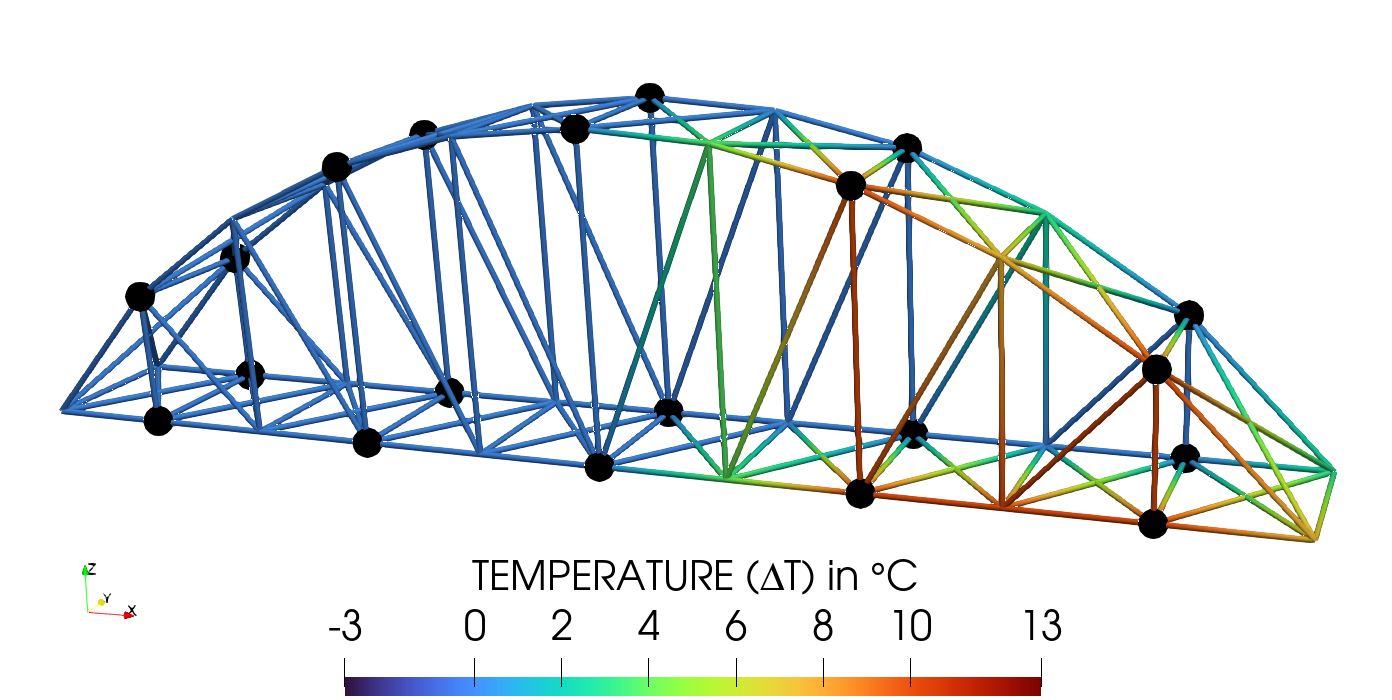}
        \end{minipage}
        \hfill
        \caption{kNN: 20 temperature sensors configuration}
        \label{fig:Figure_22b}
    \end{subfigure}
    \caption{Bridge example, 8 and 20 temperature sensors configurations: k-Nearest Neighbor interpolation (kNN) with $\textit{k}=3$ and inverse distance weighting.}
    \label{fig:Figure_22}
\end{figure}

Since the interpolation just uses sensor temperature data as the training dataset, it is paramount that the temperature sensors detect the temperature change. If the training dataset is poor, i.e., it cannot fully capture all the temperature variations, then the interpolated values will also be deficient.

Figure \ref{fig:Figure_21} shows the k-nearest neighbor interpolated temperature distributions for the 6 and 14 temperature sensors configurations for the Plate With a Hole example. The temperature sensor data was extracted at the measurement locations from the prescribed temperature distribution specified in Section \ref{sec:Plate With a Hole} (shown in Figure \ref{fig:Figure_3} (right)). For both sensor configurations, 3 nearest neighbors with inverse distance weighting were specified.
It can be observed from Figure \ref{fig:Figure_21a} that due to low number of sensors, none of the temperature sensors are able to detect the heated region, hence the interpolated temperature field is $\Delta T = 0$ everywhere.

\begin{figure}[!b]
    \centering
    \begin{subfigure}[t]{\textwidth}
        \centering
        \includegraphics[trim=0 0 0 0, clip, width=0.62\textwidth]{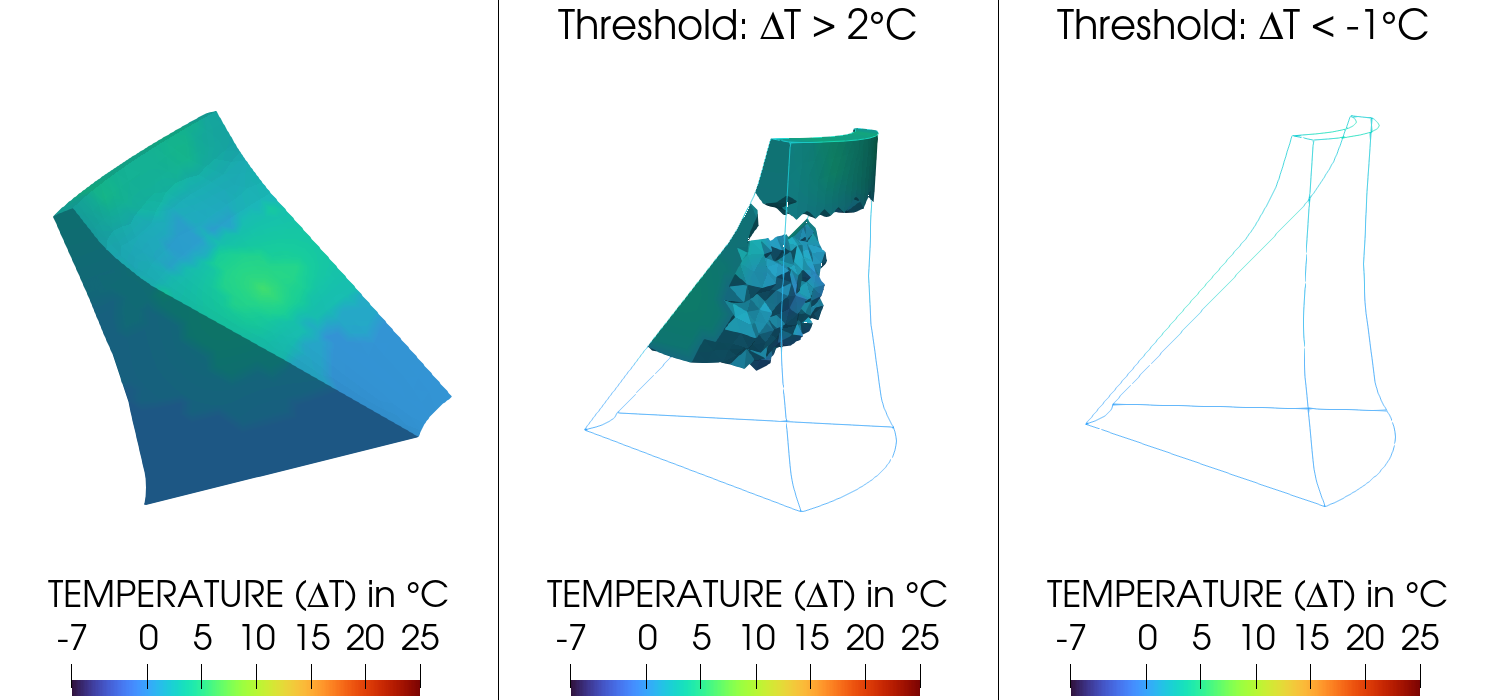}
        \caption{kNN: 27 temperature sensors configuration}
        \label{fig:Figure_23a}
        \vspace{1em}
    \end{subfigure}
    \begin{subfigure}[t]{\textwidth}
        \centering
        \includegraphics[trim=0 0 0 0, clip, width=0.62\textwidth]{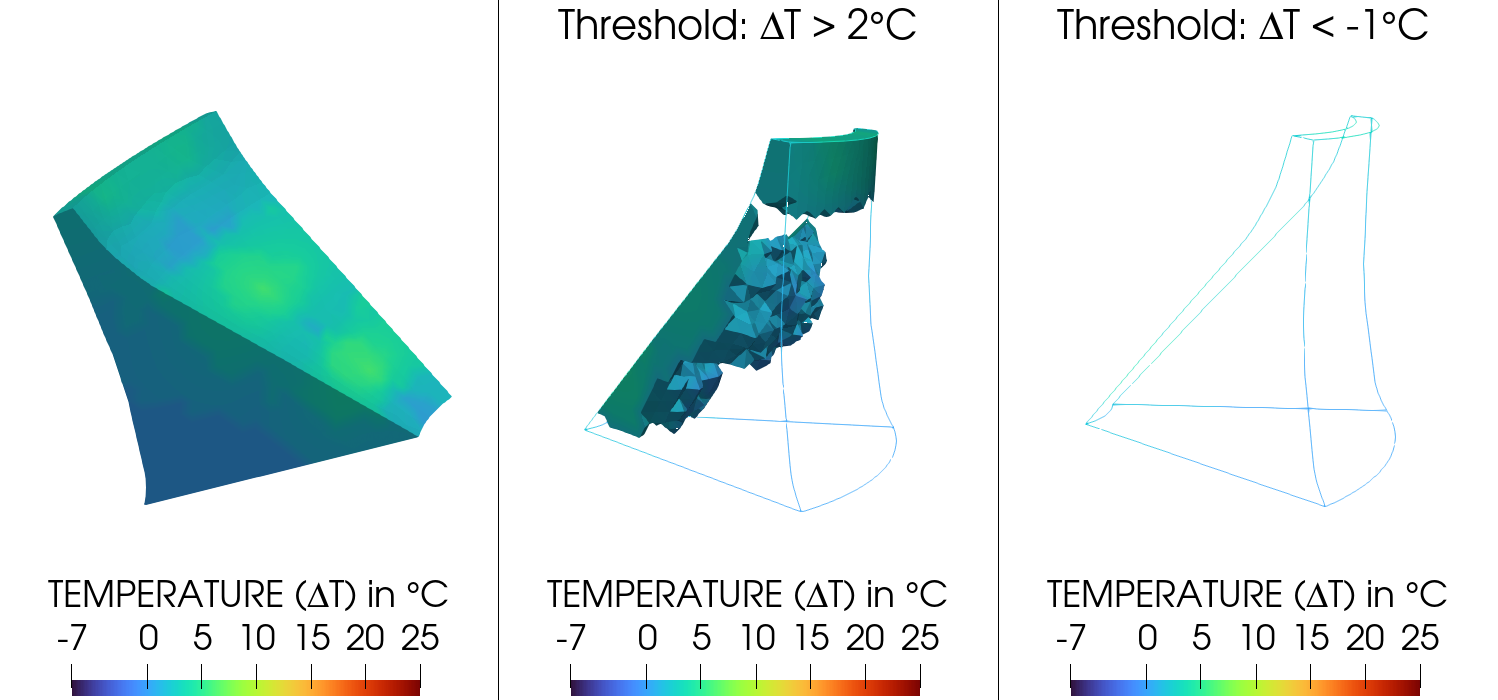}
        \caption{kNN: 36 temperature sensors configuration}
        \label{fig:Figure_23b}
        \vspace{1em}
    \end{subfigure}
    \begin{subfigure}[t]{\textwidth}
        \centering
        \includegraphics[trim=0 0 0 0, clip, width=0.62\textwidth]{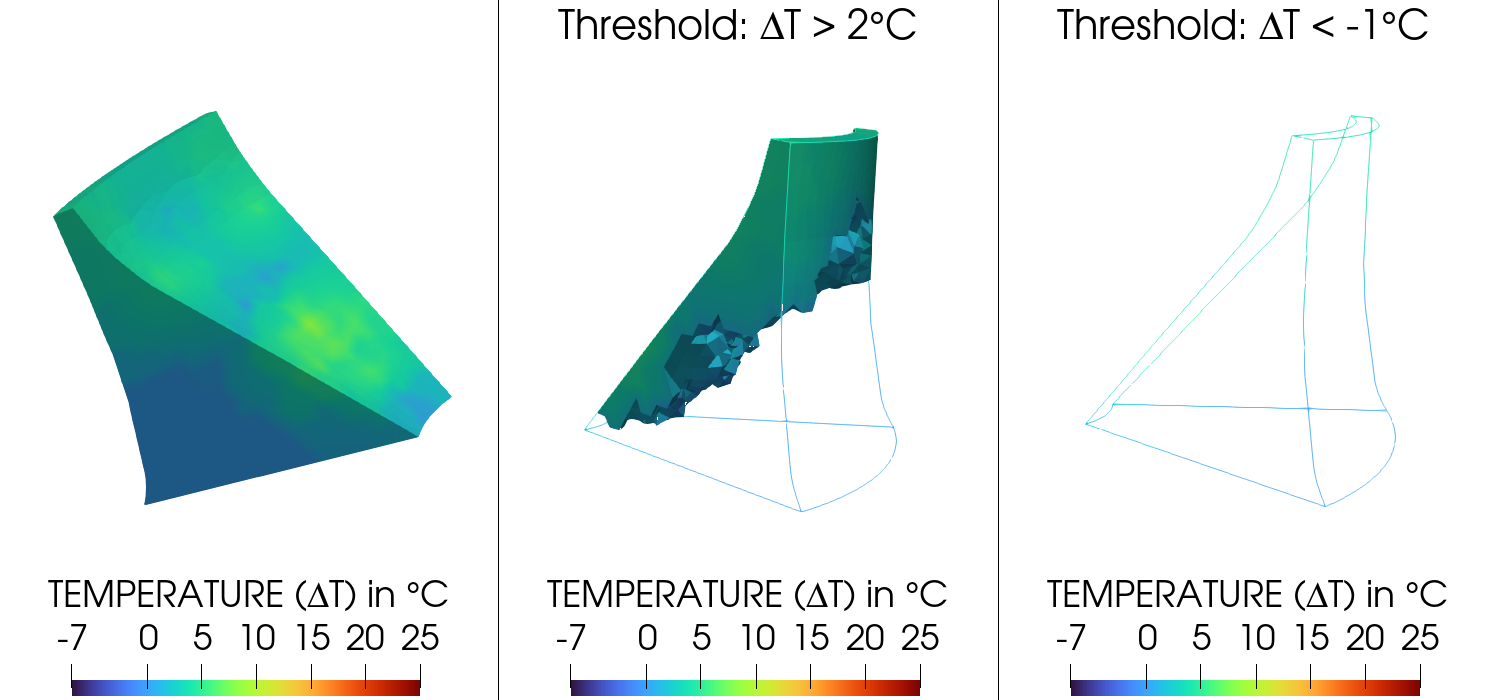}
        \caption{kNN: 59 temperature sensors configuration}
        \label{fig:Figure_23c}
    \end{subfigure}
    \caption{Dam example with 27, 36, and 59 temperature sensors configurations: k-Nearest Neighbor interpolation (kNN) with $\textit{k}=5$ and inverse distance weighting: Interpolated temperature distributions (left) along with the $\Delta T > 2$ $\degree C$ thresholds (middle), and the $\Delta T < -1$ $\degree C$ thresholds (right).}
    \label{fig:Figure_23}
\end{figure}

For the 14 temperature sensors case as shown in Figure \ref{fig:Figure_21b}, only one of the sensors (top row second sensor) captures the $\Delta T = 10$ $\degree C$ temperature rise. Due to this, a broad, poorly interpolated temperature distribution is obtained, which is mainly centered around the sensor location rather than the actual heated region.  

\medskip

Figure \ref{fig:Figure_22} shows the k-nearest neighbor interpolated temperature distributions for the 8 and 20 temperature sensors configurations for the Bridge example. The temperature sensor data was extracted at the measurement locations from the prescribed temperature distribution specified in Section \ref{sec:bridge} (shown in Figure \ref{fig:Figure_8} (top)). For both the sensor configurations, 3 nearest neighbors with inverse distance weighting were specified.
It can be observed from Figures \ref{fig:Figure_22a} and \ref{fig:Figure_22b} that the temperature in the left section of the bridge is accurate as the sensors correctly detect no change in temperature in this section. Whilst, for the right-section of the bridge, the temperature interpolation in this region is quite deficient. Many nodes on the rear-right section of the bridge, which are prescribed to be unheated ($\Delta T = 0$ $\degree C$), have interpolated temperatures up to $3.43$ $\degree C$. Nevertheless, as expected, due to the higher number of sensors, the temperature interpolation for the 20 temperature sensors configuration (Figure \ref{fig:Figure_22b}) is markedly better than that obtained for the 8 temperature sensors configuration (Figure \ref{fig:Figure_22a}). 

\medskip

Figure \ref{fig:Figure_23} shows the k-nearest neighbor interpolated temperature distributions and threshold distributions of $\Delta T > 2$ $\degree C$ and $\Delta T < -1$ $\degree C$ for the 27, 36, and 59 temperature sensors configurations for (simplified) Hoover dam example. The temperature sensor data was extracted at the measurement locations from the prescribed temperature distribution specified in Section \ref{sec:Hoover dam} (shown in Figure \ref{fig:Figure_14} (middle)). For all three sensor configurations, 5 nearest neighbors with inverse distance weighting were specified.
It can be observed from Figure \ref{fig:Figure_23a} that due to a lack of temperature sensors in the lower downstream portion of the dam, no temperature rise is detected after interpolation. 
Due to the addition of 9 temperature sensors in the lower downstream section of the dam, the interpolated temperature distribution obtained from the 36 sensors configuration as shown in Figure \ref{fig:Figure_23b} is better than the 27 sensors configuration.

The 59 sensors configuration is able to get an even better temperature interpolation than the 27 and 36 sensors configurations, but still, the temperature distribution is quite different from prescribed $\Delta \mathbf{T}$ as shown in Figure \ref{fig:Figure_14}. Temperature hotspots (with maximum $\Delta T$ of 8.2\degree C) are identified on the downstream surface close to the temperature sensors but quickly vanish away from the sensors producing highly varying temperatures on the surface. Since there are no negative $\Delta T$'s in the prescribed temperature distribution, the \acrshort{kNN} interpolation also does not produce any negative $\Delta T$. 

\FloatBarrier
\subsubsection{Kriging interpolation}

Kriging interpolation works on the basis of spatial correlation between points to interpolate a quantity spatially based on limited sample / measurement data.
The spatial correlation, i.e., how the similarity decreases with distance, is described using models called \textit{variograms}. Common variograms include Gaussian, linear, power, exponential, etc models.
The kriging weights are obtained by fitting the sampled data and the variogram. The kriging weights are then used to interpolate the quantity at unsampled locations.
When the underlying spatial process is considered to be stationary, i.e., the mean is constant, then it is called \textit{ordinary kriging}, whereas when the mean variation contains a trend, it is called \textit{universal kriging}. The deterministic trend can be 'regional linear', meaning the trend is a linear function of the spatial coordinates.
The choice of the type of kriging depends on the underlying spatial process being interpolated. An explanation of the 3D temperature field reconstruction using kriging interpolation can be found in \cite{lin20213d}.
In this work, the standard kriging implementation from the Python PyKrige library (\cite{pykrige}) was used.

\begin{figure}[!b]
\begin{minipage}[c][][t]{.5\textwidth}
  \vspace*{\fill}
  \centering
  \includegraphics[trim= 0 0 0 0, clip, width=0.3\paperwidth]{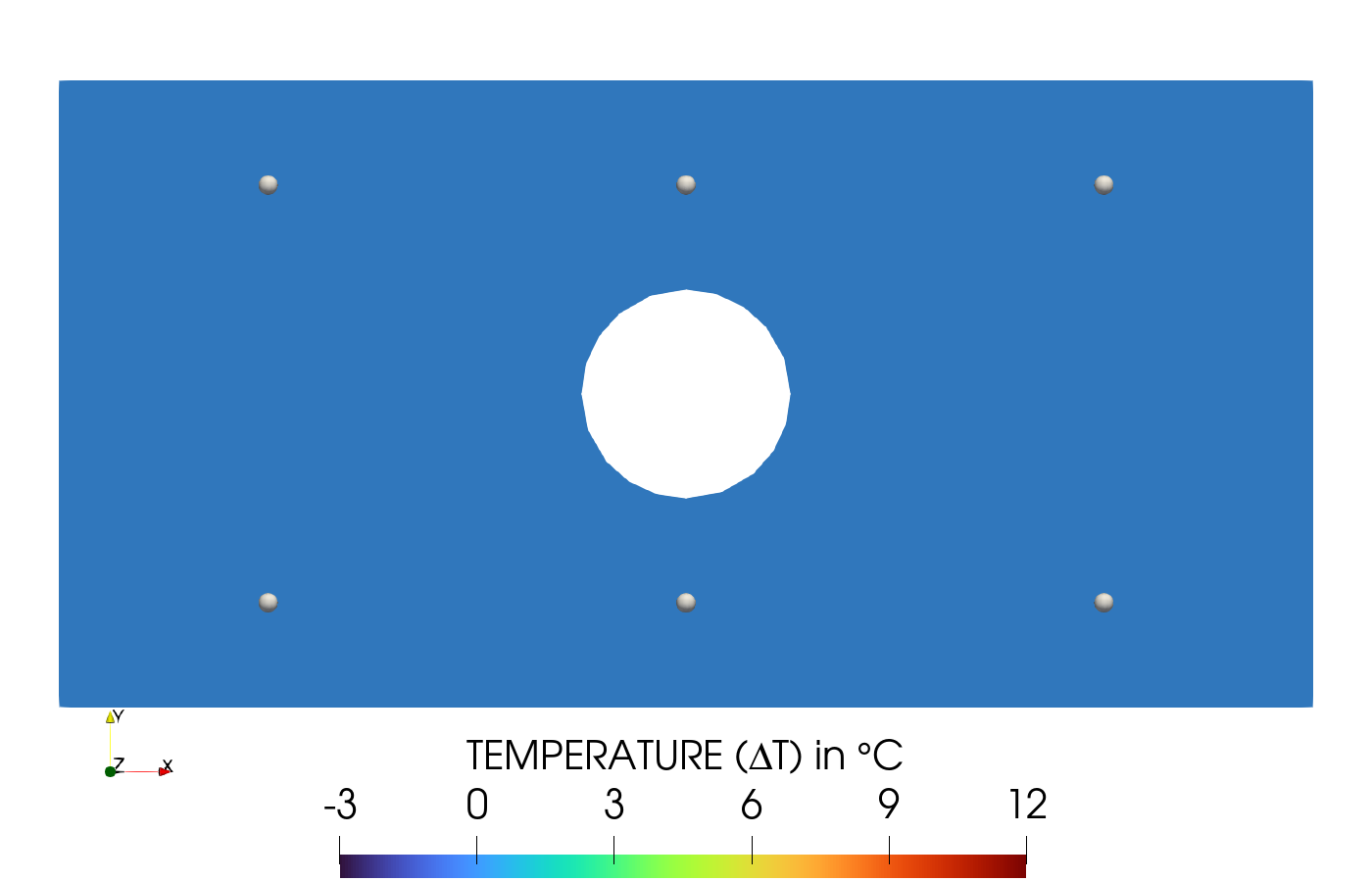}
  \subcaption{OK: 6 temperature sensors configuration}
   \label{fig:Figure_24a}
\end{minipage}
\begin{minipage}[c][][t]{.5\textwidth}
  \vspace*{\fill}
  \centering
  \includegraphics[trim= 0 0 0 0, clip, width=0.3\paperwidth]{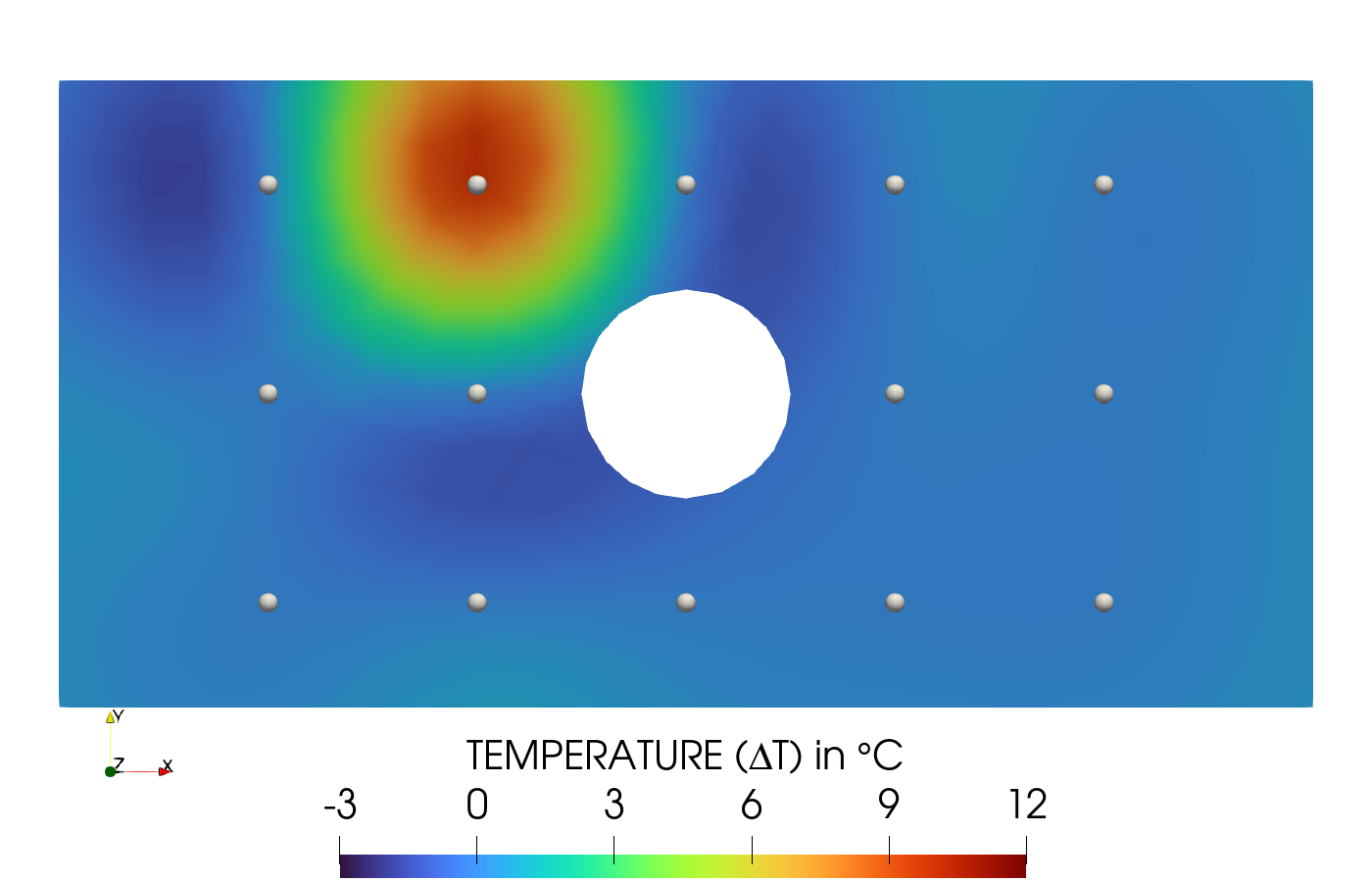}
  \subcaption{OK: 14 temperature sensors configuration}
  \label{fig:Figure_24b}\par
\end{minipage}
\caption{Plate With a Hole, 6 and 14 temperature sensors configurations: Ordinary kriging (OK) interpolation with Gaussian variogram.}
\label{fig:Figure_24}
\end{figure}

\medskip

\begin{figure}[!t]
\begin{minipage}[c][][t]{.5\textwidth}
  \vspace*{\fill}
  \centering
  \includegraphics[trim= 0 0 0 0, clip, width=0.3\paperwidth]{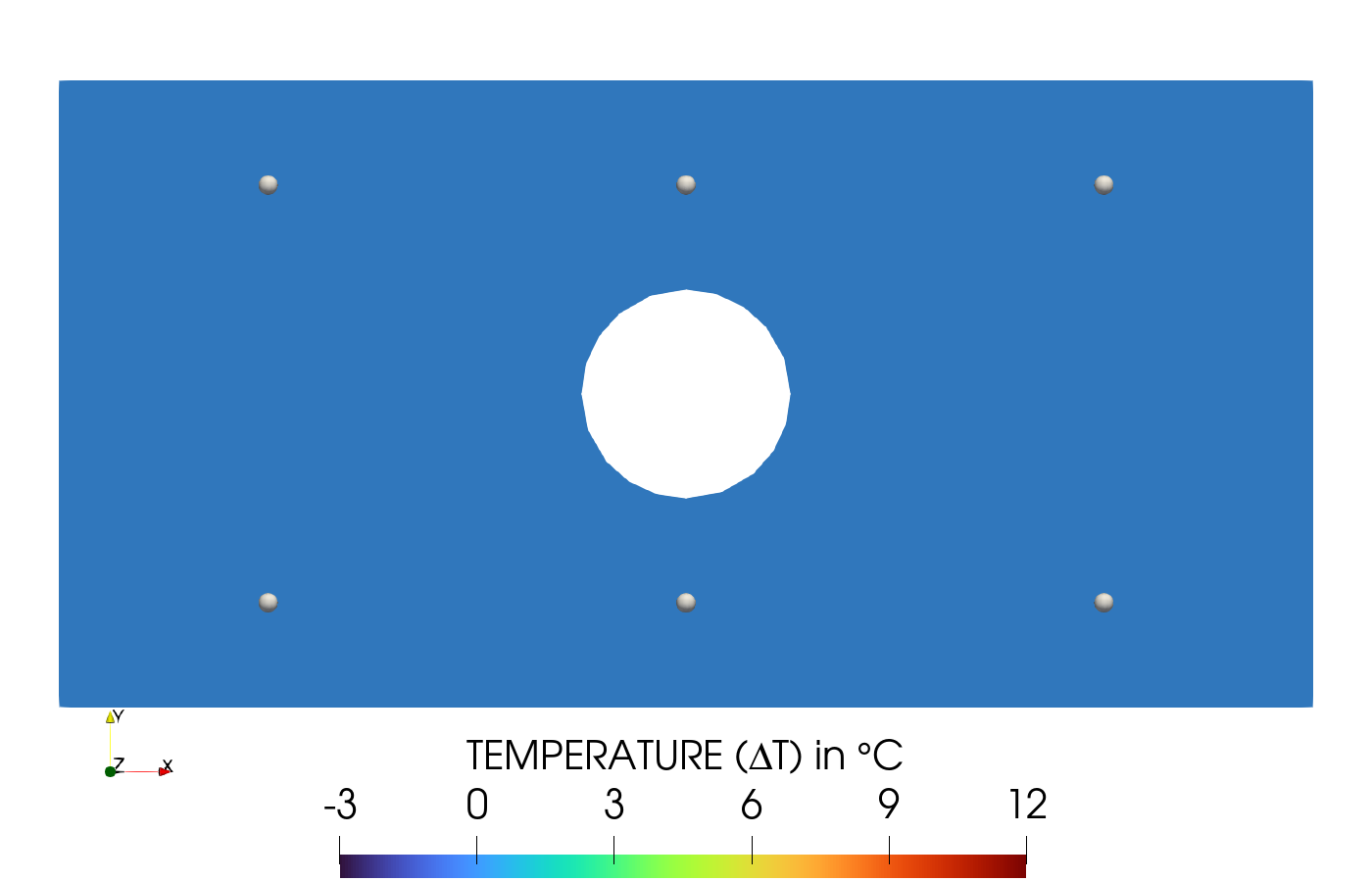}
  \subcaption{UK: 6 temperature sensors configuration}
   \label{fig:Figure_25a}
\end{minipage}
\begin{minipage}[c][][t]{.5\textwidth}
  \vspace*{\fill}
  \centering
  \includegraphics[trim= 0 0 0 0, clip, width=0.3\paperwidth]{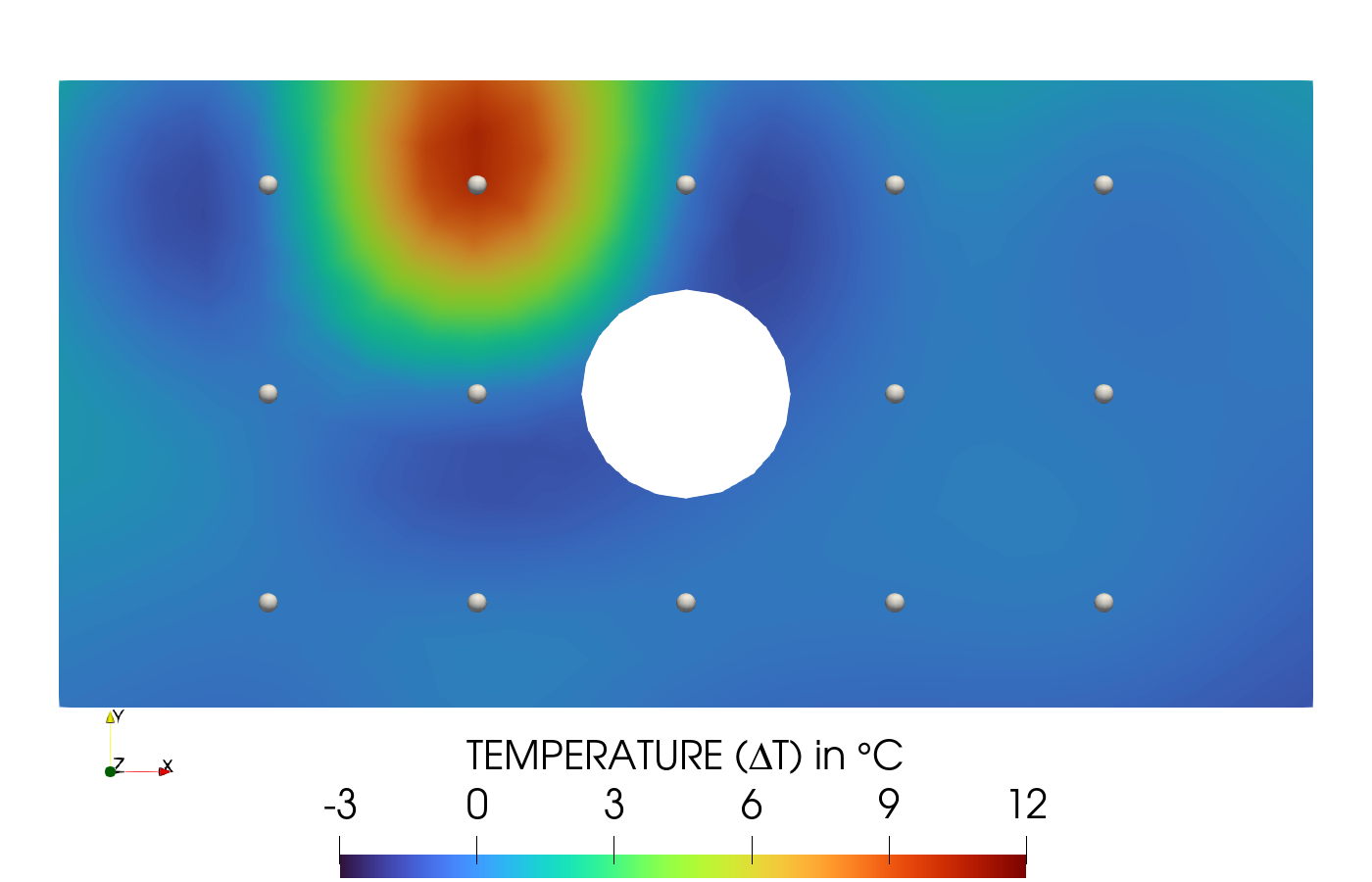}
  \subcaption{UK: 14 temperature sensors configuration}
  \label{fig:Figure_25b}\par
\end{minipage}
\caption{Plate With a Hole, 6 and 14 temperature sensors configurations: Universal kriging (UK) interpolation with Gaussian variogram and 'regional linear' trend.}
\label{fig:Figure_25}
\end{figure}

\begin{figure}[!b]
    \centering
    \begin{subfigure}[t]{\textwidth}
        \centering
        \begin{minipage}[t]{0.62\textwidth}
            \centering
            \includegraphics[trim=0 0 0 80, clip, width=\textwidth]{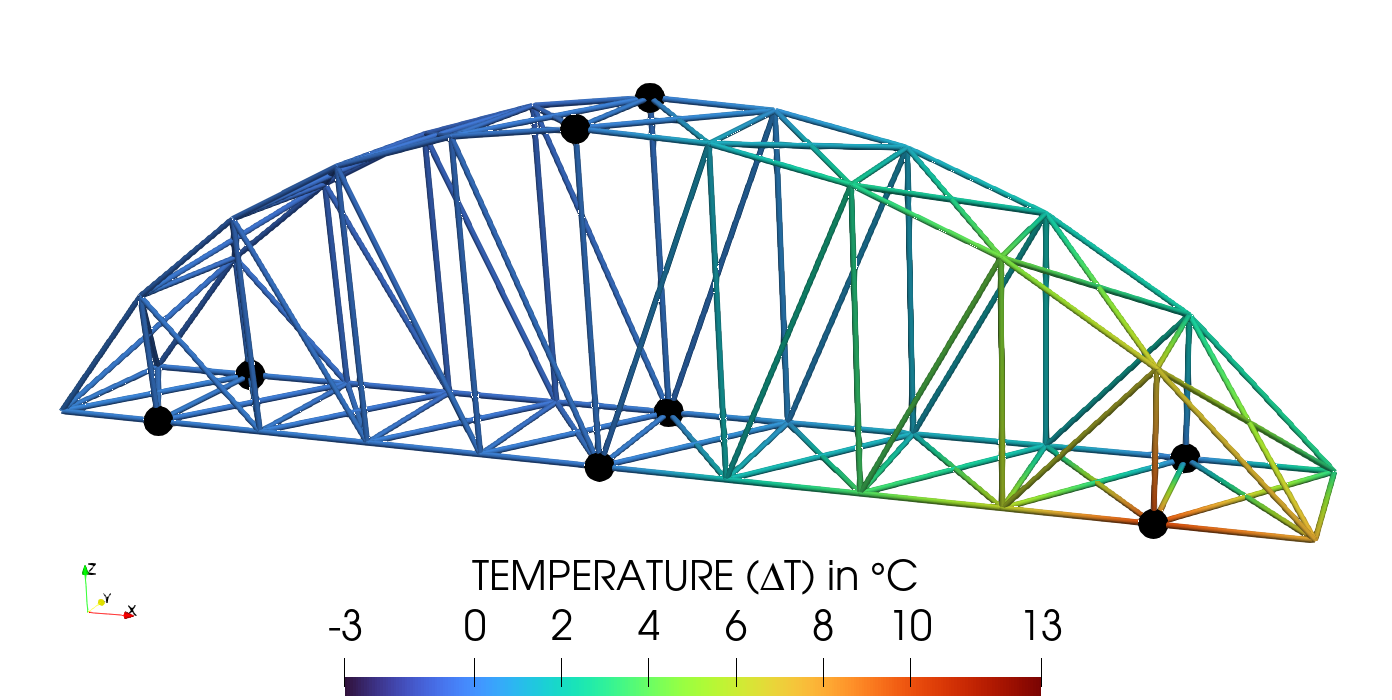}
        \end{minipage}
        \hfill
        \caption{OK: 8 temperature sensors configuration}
        \label{fig:Figure_26a}
    \end{subfigure}
    \begin{subfigure}[t]{\textwidth}
        \centering
        \begin{minipage}[t]{0.62\textwidth}
            \centering
            \includegraphics[trim=0 0 0 0, clip, width=\textwidth]{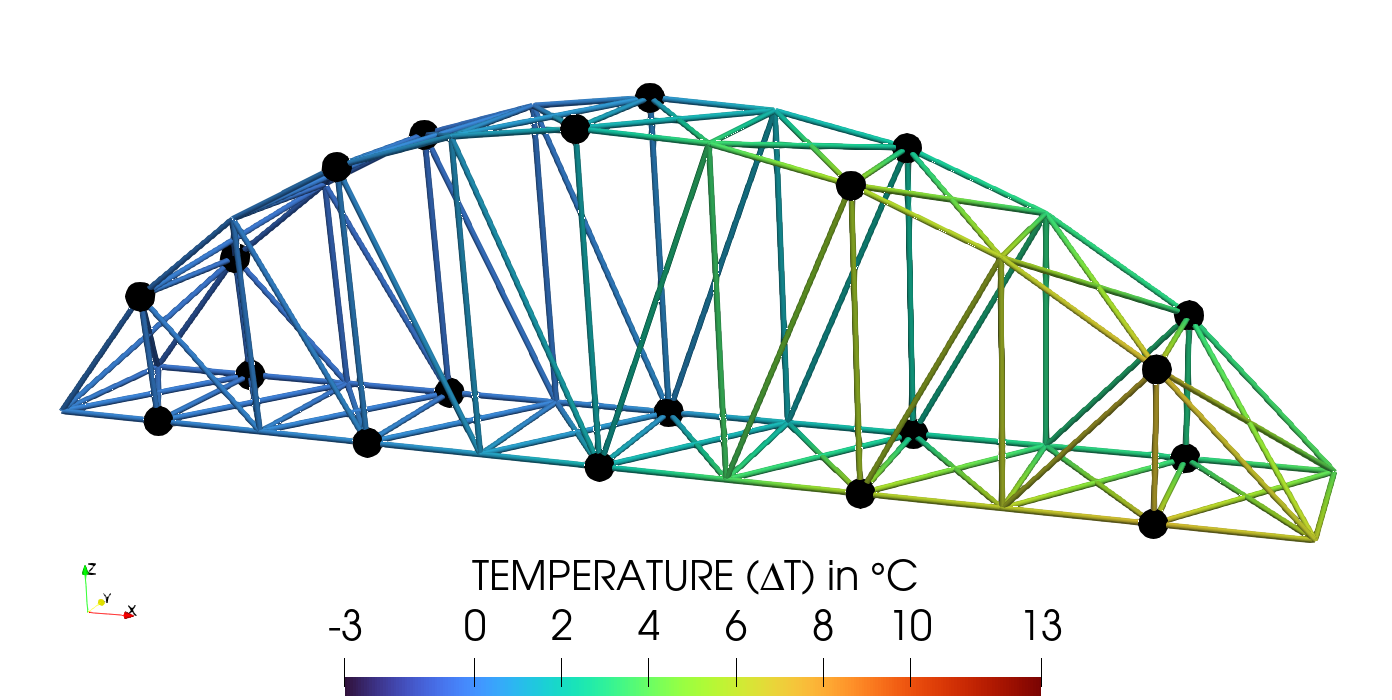}
        \end{minipage}
        \hfill
        \caption{OK: 20 temperature sensors configuration}
        \label{fig:Figure_26b}
    \end{subfigure}
    \caption{Bridge example, 8 and 20 temperature sensors configurations: Ordinary kriging (OK) interpolation with Gaussian variogram.}
    \label{fig:Figure_26}
\end{figure}

For the Plate With a Hole example with 6 and 14 temperature sensors configurations, Figures \ref{fig:Figure_24} and \ref{fig:Figure_25} show the ordinary kriging and universal kriging interpolated temperature distributions, respectively. For both, ordinary and universal kriging, the Gaussian variogram was specified to quantify the spatial correlation; and a regional linear trend was considered for universal kriging.
The temperature sensor data was extracted from the same prescribed temperature distribution specified in Section \ref{sec:Plate With a Hole} (shown in Figure \ref{fig:Figure_3} (right)). 

For the 6 temperature sensors configurations, a very small noise (in the order of $10^{-7}$) was added to two sensors (top row; left and middle sensor) measurement data since all 6 temperature sensors recorded $\Delta T = 0$ and the PyKrige implementation required at least very small variation in the measurement dataset. This did not change the outcome as from Figures \ref{fig:Figure_24a} and \ref{fig:Figure_25a} it can be seen that the interpolated temperature distribution all over the plate is $\Delta T \xrightarrow{} 0$, consistent with the \acrshort{kNN} results. 

The 14 temperature sensors configurations gave better and smoother temperature distributions based on the information of one sensor detecting the high temperature. Similar to the \acrshort{kNN} case, the high-temperature zone is centered around this one sensor rather than the actual heating area. No information is available regarding whether the heated zone is to the left, right, below, or above the sensor, hence the circular pattern.
In this case, no significant difference is observed between the ordinary and universal kriging interpolations.

\begin{figure}[!t]
    \centering
    \begin{subfigure}[t]{\textwidth}
        \centering
        \begin{minipage}[t]{0.62\textwidth}
            \centering
            \includegraphics[trim=0 0 0 80, clip, width=\textwidth]{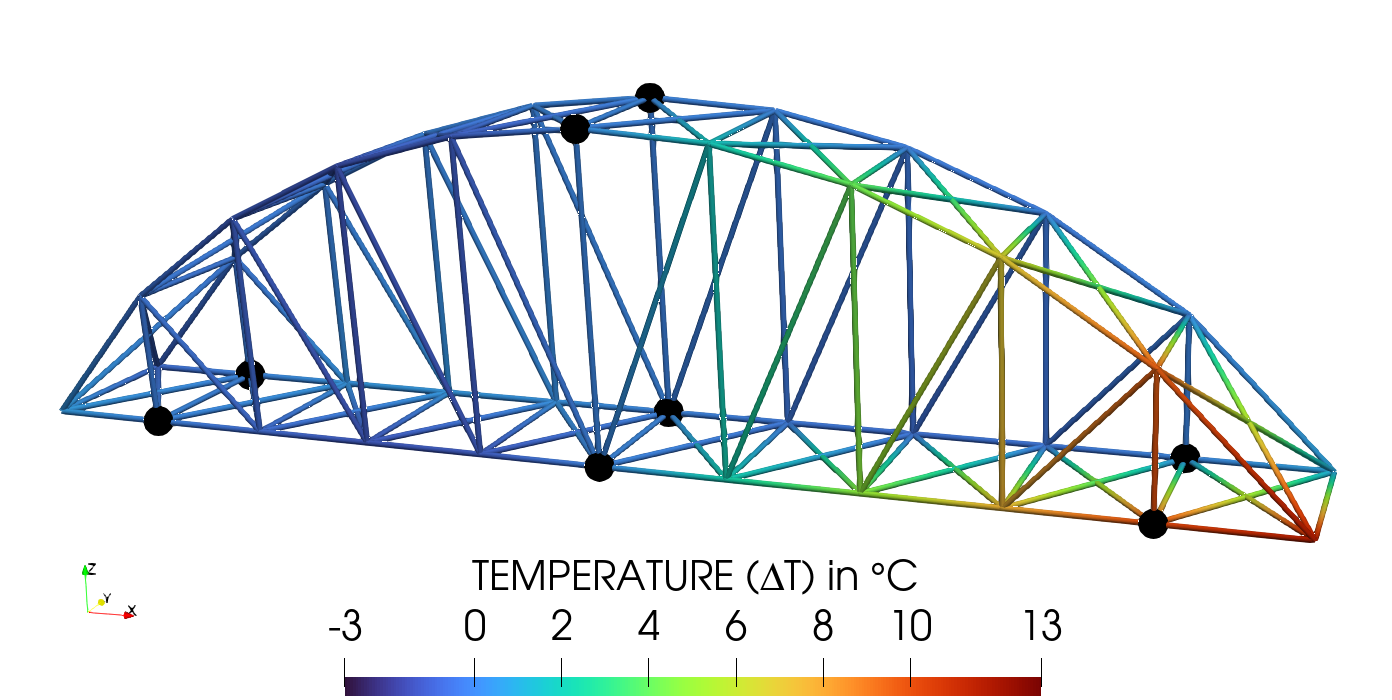}
        \end{minipage}
        \hfill
        \caption{UK: 8 temperature sensors configuration}
        \label{fig:Figure_27a}
    \end{subfigure}
    \begin{subfigure}[t]{\textwidth}
        \centering
        \begin{minipage}[t]{0.62\textwidth}
            \centering
            \includegraphics[trim=0 0 0 20, clip, width=\textwidth]{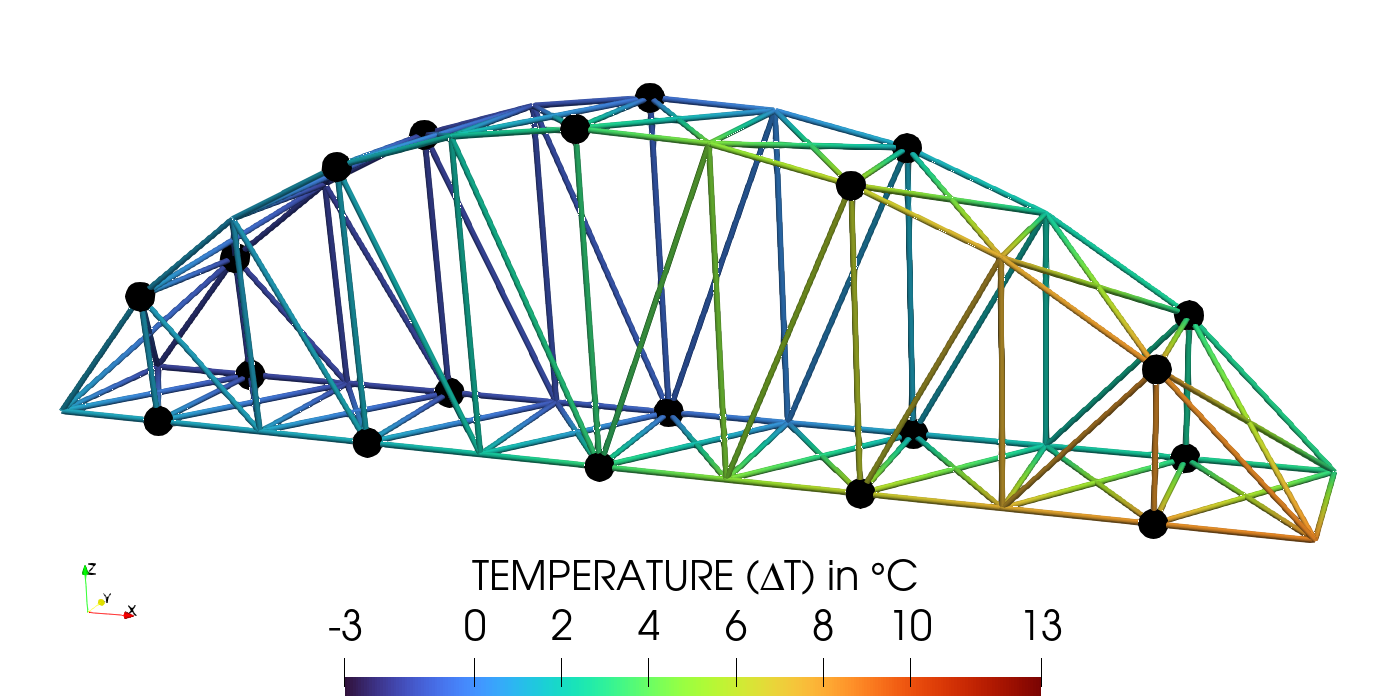}
        \end{minipage}
        \hfill
        \caption{UK: 20 temperature sensors configuration}
        \label{fig:Figure_27b}
    \end{subfigure}
    \caption{Bridge example, 8 and 20 temperature sensors configurations: Universal kriging (UK) interpolation with Gaussian variogram and 'regional linear' trend.}
    \label{fig:Figure_27}
\end{figure}

\medskip
For the Bridge example with 8 and 20 temperature sensors configurations, Figures \ref{fig:Figure_26} and \ref{fig:Figure_27} illustrate the ordinary kriging and universal kriging interpolated temperature distributions, respectively.
For both, ordinary and universal kriging, the Gaussian variogram was specified to quantify the spatial correlation; and a regional linear trend was considered for universal kriging.
The temperature sensor data was extracted from the same prescribed temperature distribution specified in Section \ref{sec:bridge} (shown in Figure \ref{fig:Figure_8} (top)).

For the 8 sensors configuration, both the kriging methods produced better temperature distributions (refer to Figures \ref{fig:Figure_26a} and \ref{fig:Figure_27a}) than kNN interpolation (shown in Figure \ref{fig:Figure_22a}). Among the kriging techniques, universal kriging generated marginally better thermal fields than ordinary kriging. The peak temperature was higher in the case of universal kriging (maximum $\Delta T$ of $9.99$ $\degree C$) compared to the ordinary kriging case (maximum $\Delta T$ of $11.64$ $\degree C$). However, universal kriging results suffered from larger negative peak temperatures (minimum $\Delta T$ of $-1.21$ $\degree C$) compared to the small negative temperatures seen in the ordinary kriging case (minimum $\Delta T$ of $-0.56$ $\degree C$).

In contrast to the results of the 8 sensors configuration, both kriging methods produced subpar results (refer to Figures \ref{fig:Figure_26b} and \ref{fig:Figure_27b}) than the kNN interpolation (shown in Figure \ref{fig:Figure_22b}) for the 20 sensors configuration.
However, in this case also, universal kriging generated slightly better results than ordinary kriging. The peak temperature was also higher in the case of universal kriging (maximum $\Delta T$ of $8.54$ $\degree C$) compared to the ordinary kriging case (maximum $\Delta T$ of $6.96$ $\degree C$). Universal kriging results also suffered from larger negative peak temperatures (minimum $\Delta T$ of $-1.67$ $\degree C$) compared to the small negative temperatures seen in the ordinary kriging case (minimum $\Delta T$ of $-0.31$ $\degree C$).
In general, it can be seen that consideration of a trend in the mean of the underlying spatial process, i.e., universal kriging, results in better results.

\medskip

\begin{figure}[!b]
    \centering
    \begin{subfigure}[t]{\textwidth}
        \centering
        \includegraphics[trim=0 0 0 0, clip, width=0.62\textwidth]{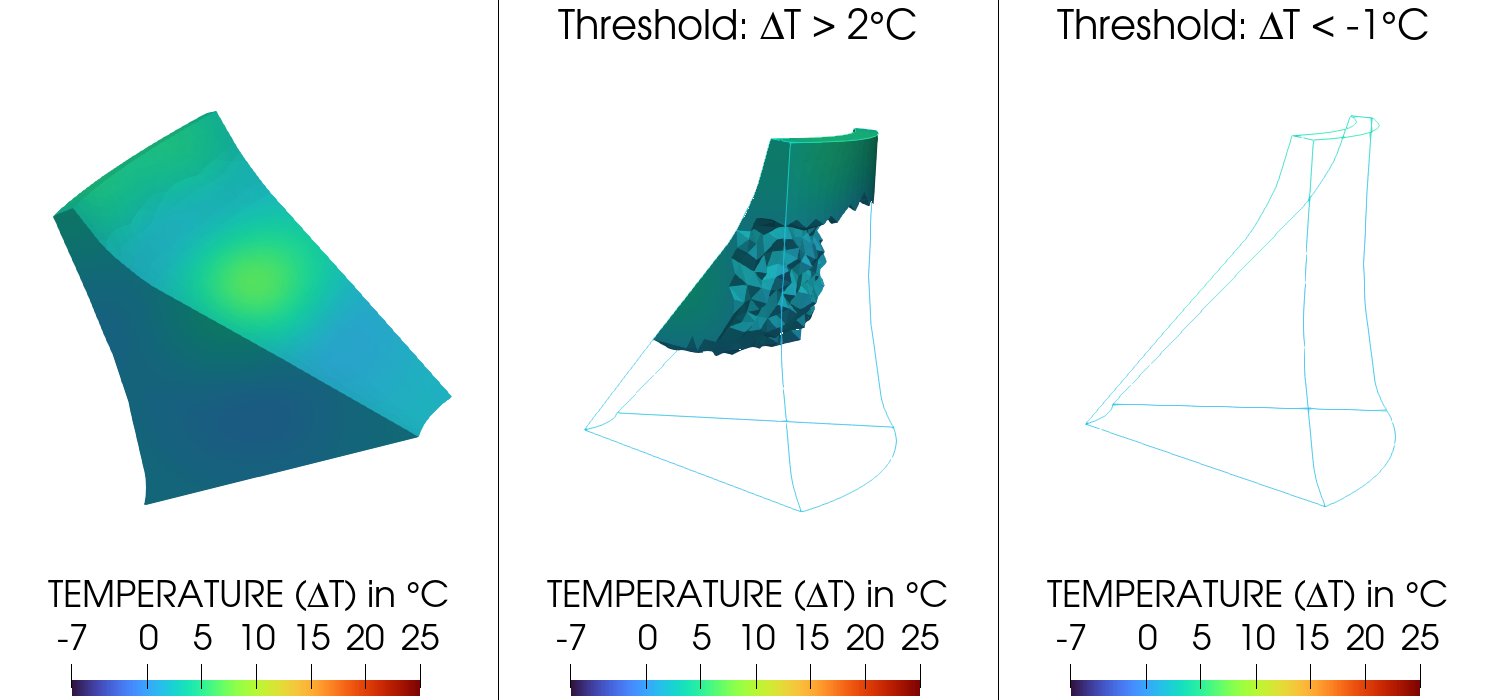}
        \caption{OK: 27 temperature sensors configuration}
        \label{fig:Figure_28a}
        \vspace{1em}
    \end{subfigure}
    \begin{subfigure}[t]{\textwidth}
        \centering
        \includegraphics[trim=0 0 0 -20, clip, width=0.62\textwidth]{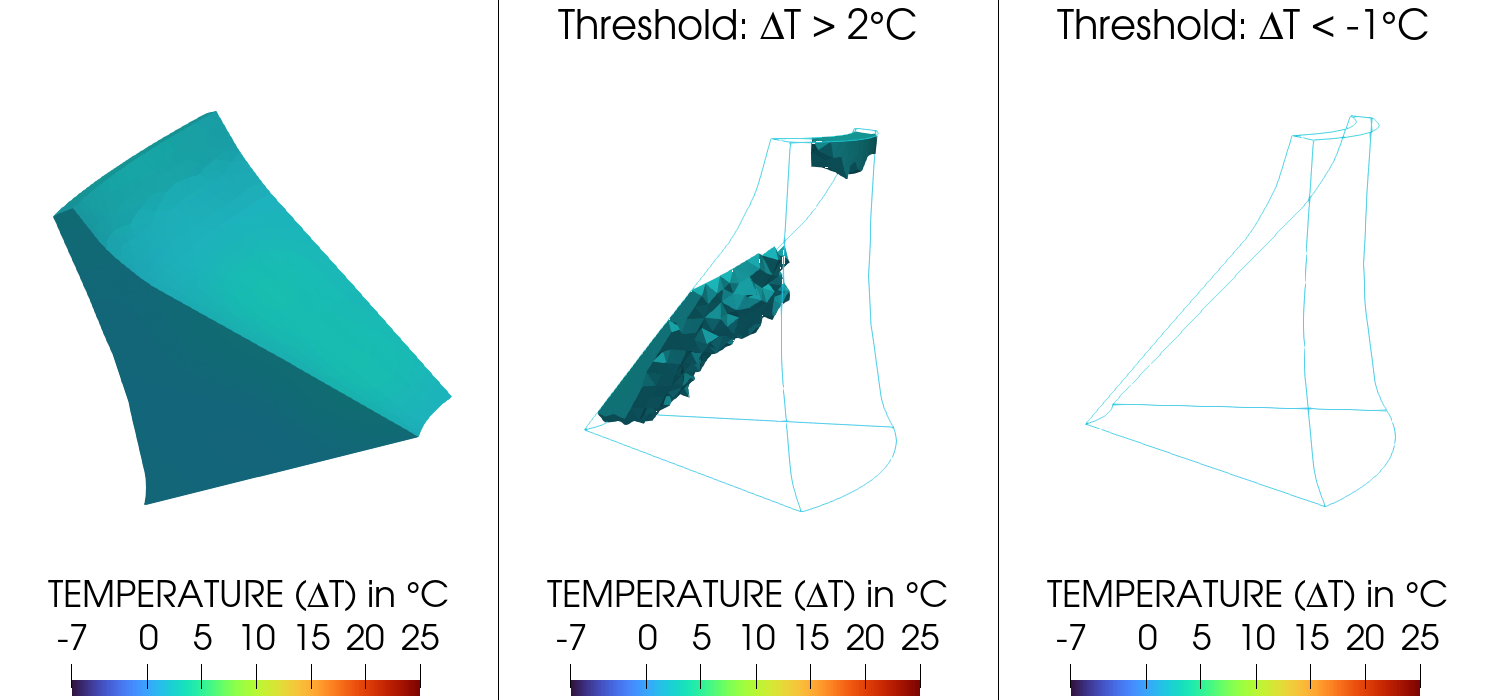}
        \caption{OK: 36 temperature sensors configuration}
        \label{fig:Figure_28b}
        \vspace{1em}
    \end{subfigure}
    \begin{subfigure}[t]{\textwidth}
        \centering
        \includegraphics[trim=0 0 0 -20, clip, width=0.62\textwidth]{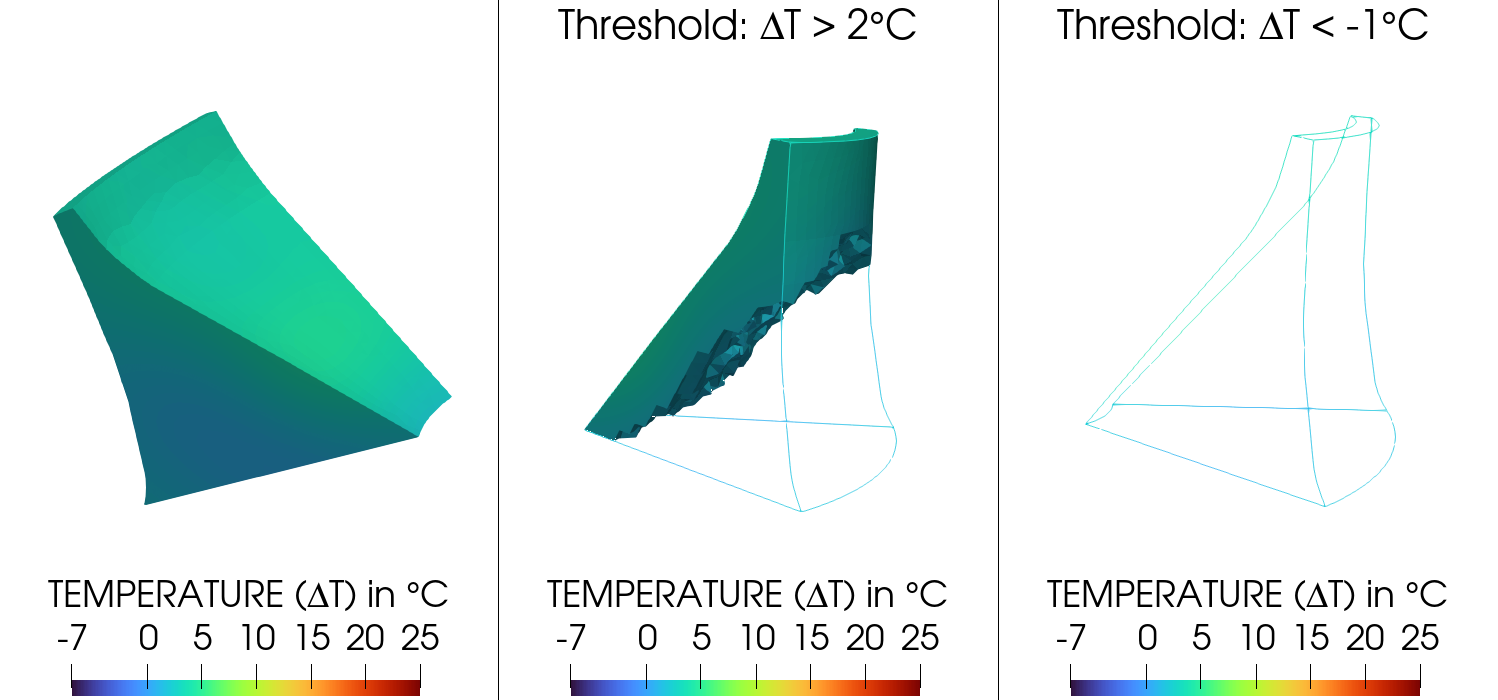}
        \caption{OK: 59 temperature sensors configuration}
        \label{fig:Figure_28c}
    \end{subfigure}
    \caption{Dam example with 27, 36, and 59 temperature sensors configurations: Ordinary kriging (OK) interpolation with Gaussian variogram: Interpolated temperature distributions (left) along with the $\Delta T > 2$ $\degree C$ thresholds (middle), and the $\Delta T < -1$ $\degree C$ thresholds (right).}
    \label{fig:Figure_28}
\end{figure}

\begin{figure}[!b]
    \centering
    \begin{subfigure}[t]{\textwidth}
        \centering
        \includegraphics[trim=0 0 0 0, clip, width=0.62\textwidth]{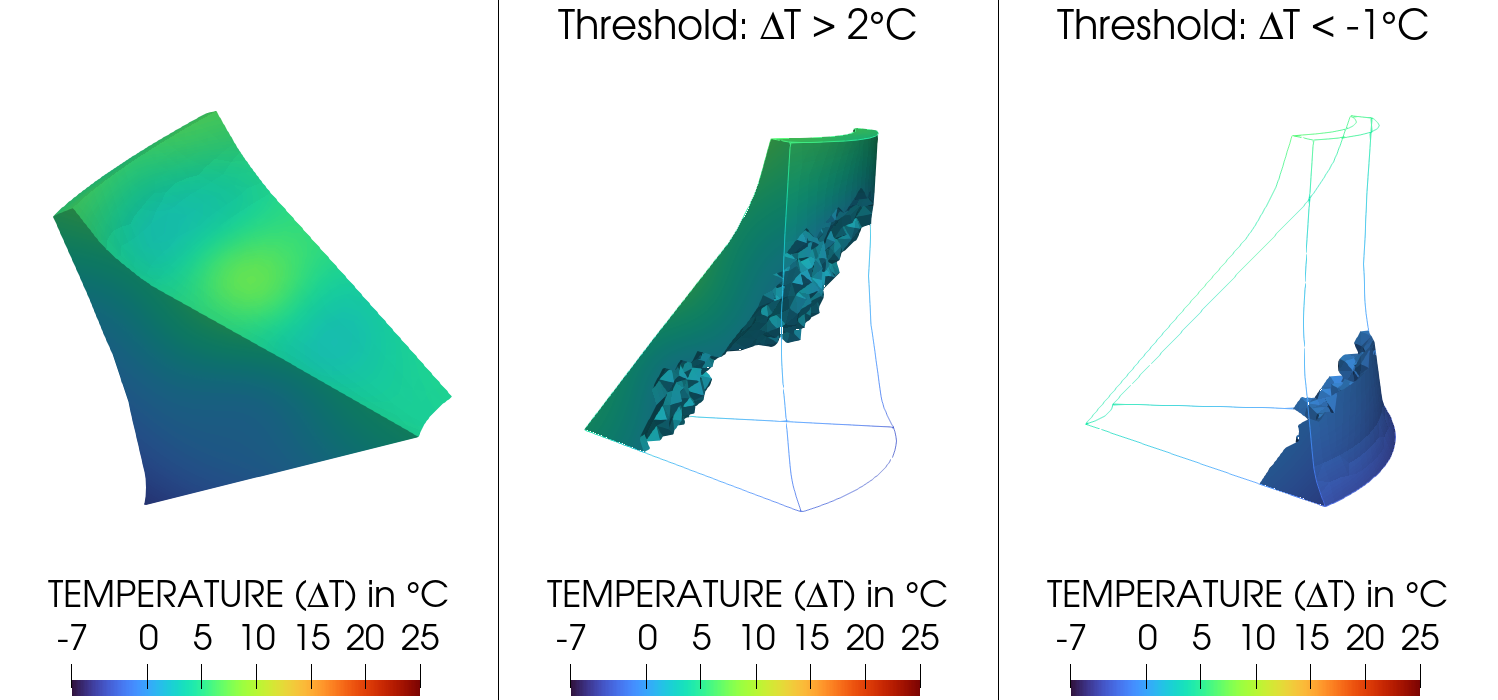}
        \caption{UK: 27 temperature sensors configuration}
        \label{fig:Figure_29a}
        \vspace{1em}
    \end{subfigure}
    \begin{subfigure}[t]{\textwidth}
        \centering
        \includegraphics[trim=0 0 0 0, clip, width=0.62\textwidth]{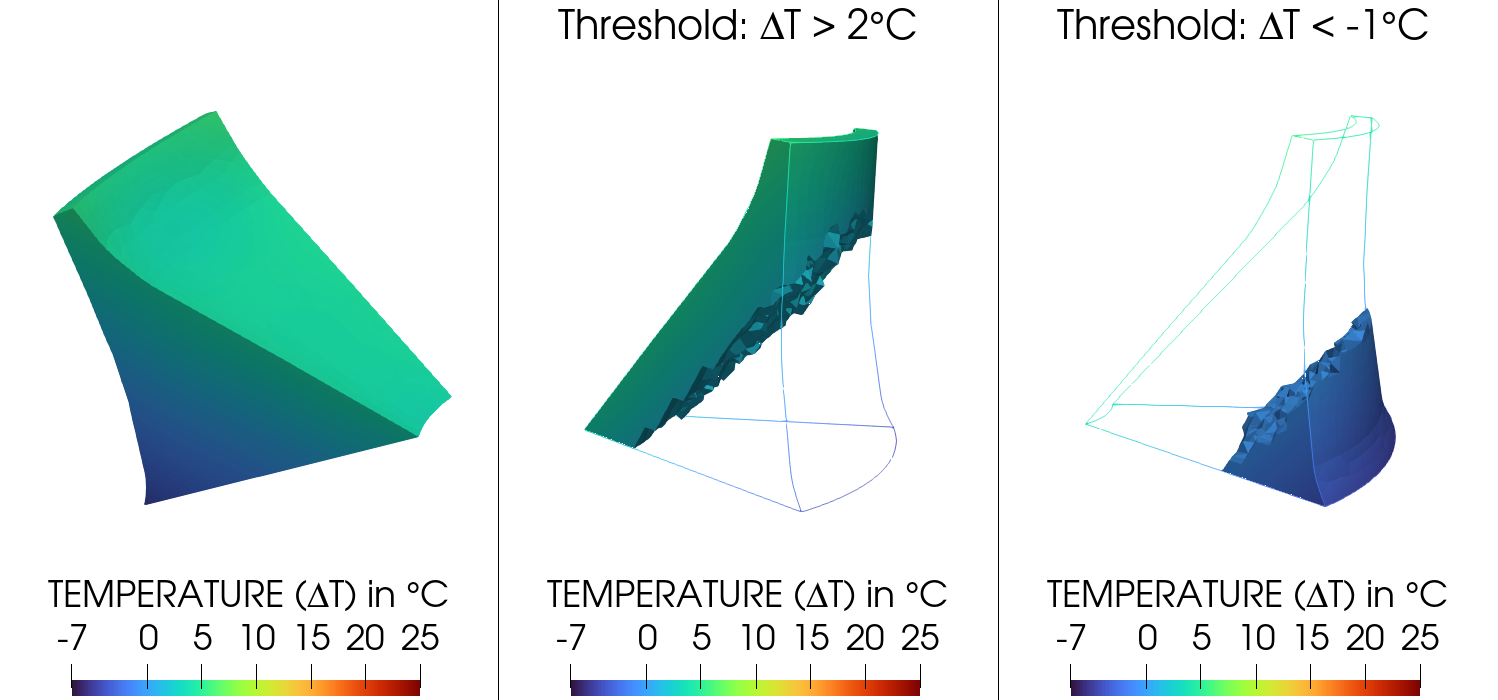}
        \caption{UK: 36 temperature sensors configuration}
        \label{fig:Figure_29b}
        \vspace{1em}
    \end{subfigure}
    \begin{subfigure}[t]{\textwidth}
        \centering
        \includegraphics[trim=0 0 0 0, clip, width=0.62\textwidth]{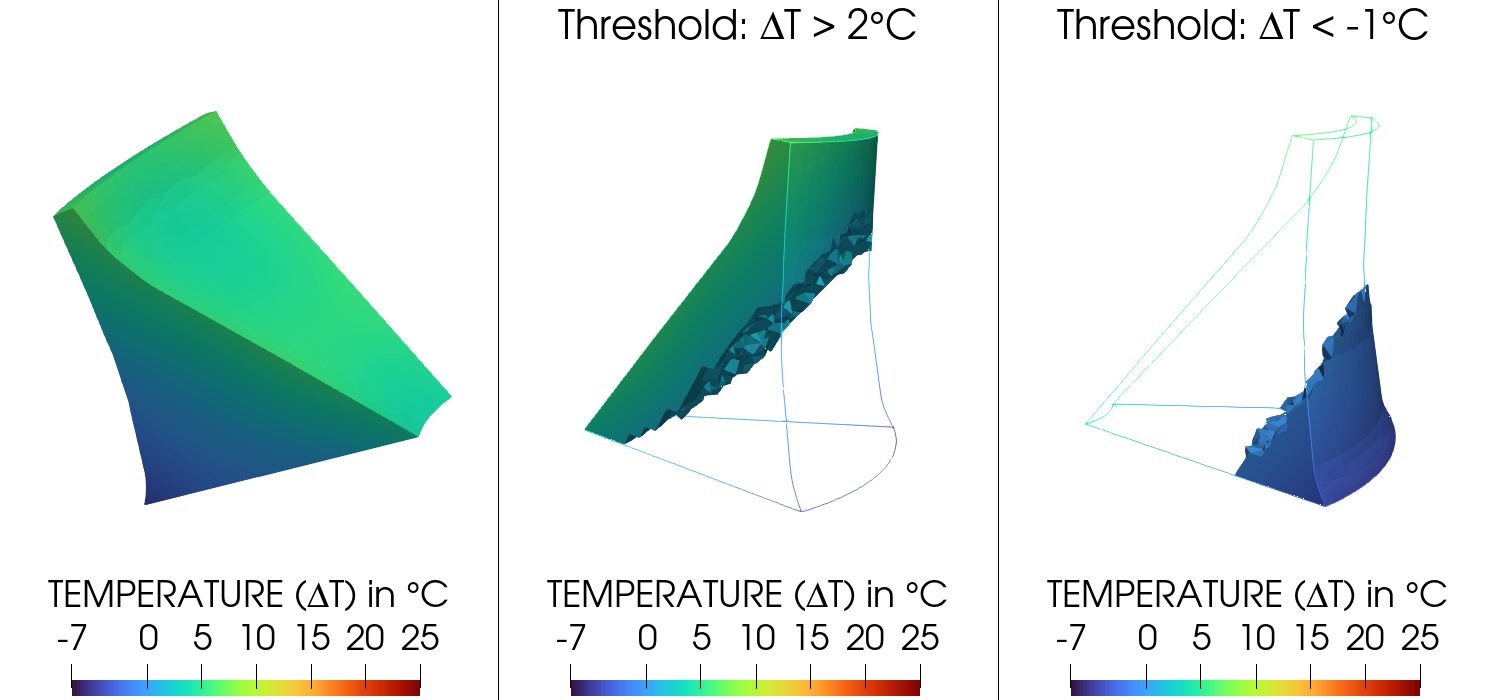}
        \caption{UK: 59 temperature sensors configuration}
        \label{fig:Figure_29c}
    \end{subfigure}
    \caption{Dam example with 27, 36, and 59 temperature sensors configurations: Universal kriging (UK) interpolation with Gaussian variogram and 'regional linear' trend: Interpolated temperature distributions (left) along with the $\Delta T > 2$ $\degree C$ thresholds (middle), and the $\Delta T < -1$ $\degree C$ thresholds (right).}
    \label{fig:Figure_29}
\end{figure}

Figures \ref{fig:Figure_28} and \ref{fig:Figure_29} respectively illustrate the ordinary kriging and universal kriging interpolated temperature distributions (left) and threshold distributions of $\Delta T > 2$ $\degree C$ (middle) and $\Delta T < -1$ $\degree C$ (right) for the 27, 36, and 59 temperature sensors configurations for the (simplified) Hoover dam example. The temperature sensor data was extracted from the same prescribed temperature distribution specified in Section \ref{sec:Hoover dam} (shown in Figure \ref{fig:Figure_14} (right)). For all three sensor configurations, the Gaussian variogram was specified for both ordinary and universal kriging, and a regional linear trend was considered for universal kriging.

For the 27 and 59 temperature sensors configurations (Figures \ref{fig:Figure_28a} and \ref{fig:Figure_28c} respectively), the ordinary kriging results are similar to those obtained by \acrshort{kNN} interpolation, although the temperature variation is much smoother in the ordinary kriging case. For the 36 temperature sensors configuration, ordinary kriging generates very poor results. 

From Figures \ref{fig:Figure_28a}-\ref{fig:Figure_28c}, it is evident that ordinary kriging is probably not suitable for interpolating temperatures in this setup. The downstream heated surface region is inadequately identified with only a fraction of the prescribed $\Delta T$ (maximum $\Delta T$ of $4.2$ $\degree C$ is observed). From the $\Delta T < -1$ $\degree C$ threshold distributions, it can be seen that negative $\Delta T$'s, if any, are all above $-1$ $\degree C$, i.e., close to 0 as in the (prescribed) target temperature distribution.

In comparison, looking at Figures \ref{fig:Figure_29a}-\ref{fig:Figure_29c}, universal kriging has better interpolated the thermal field than ordinary kriging. This can be attributed to the consideration of a trend in the mean of the underlying spatial process in universal kriging. Due to the presence of thermal gradients, the trend in the mean should not be ignored by considering a stationary condition.
However, in the case of universal kriging, the peak temperatures are slightly lower (maximum $\Delta T$ of around $7$ $\degree C$) than the target distribution. Furthermore, the introduction of the regional linear trend produces negative $\Delta T$'s (as low as $-4.5$ $\degree C$) towards the opposite end of the downstream surface, i.e., close to the lower upstream region. 

In terms of localization of the heated region, judged in terms of distribution on the downstream surface and the depth of the high-temperature zone, it is observed that universal kriging generates a smoother downstream surface temperature profile but the smooth transition also generates a larger depth (almost twice compared to that of the prescribed target distribution) tangentially inwards in the dam. 
In line with previous results, interpolation methods also perform better with higher numbers and better distribution of temperature sensors.

\medskip
\FloatBarrier
\subsubsection{Thermal Field Reconstruction Error Comparison}
\label{sec:comparison_error}

In this work, since the target temperature distribution in the structure is prescribed by the authors, the 'identification' error can additionally be analyzed for the overall structure. This is different from the sensor error used in the cost function calculation (Eqn.\eqref{eq:cost_function}). Here, the difference in identified temperature at each node of the \acrshort{FE} model is compared with the prescribed distribution, instead of only at certain limited sensor locations in the cost function calculation.
As explained earlier, the proposed optimization-based approach for temperature field reconstruction is based on the use of displacement/strain sensors while the interpolation techniques are based on the use of temperature sensors. When comparing the different approaches, the sensor configuration in terms of location and number are kept the same, and only the type of sensor (displacement or temperature) is changed.

\begin{table}[!b]
\caption{Plate With a Hole Example: Comparison of RMSE between the temperature distributions obtained by the different approaches and the prescribed target distribution (shown in Figure \ref{fig:Figure_3} with initial RMSE = 1.561738 i.e., the RMSE at optimization start when $\mathbf{\Delta T} = \mathbf{0}$ $\degree C)$.}
\label{tab:Table_1}
\centering
\begin{tblr}{
  width = \linewidth,
  colspec = {Q[260]Q[242]Q[185]Q[125]Q[125]},
  cells = {c},
  cell{1}{1} = {c=3}{0.687\linewidth},
  cell{2}{1} = {r=3}{},
  cell{3}{2} = {r=2}{},
  cell{5}{1} = {r=3}{},
  cell{5}{2} = {c=2}{0.427\linewidth},
  cell{6}{2} = {c=2}{0.427\linewidth},
  cell{7}{2} = {c=2}{0.427\linewidth},
  hline{1-2,5,8} = {-}{},
  hline{3,6-7} = {2-5}{},
  hline{4} = {3-5}{},
} \hline
\textbf{Sensor Configuration / Number of Sensors }   &    &               & \textbf{6 Sensors} & \textbf{14 Sensors} \\
{\textbf{Proposed }\\\textbf{optimization-based }\\\textbf{approach}\\\text{(using displacement}\\\text{sensors)}} & {\textbf{Without Vertex }\\\textbf{Morphing filtering}} & BB Step       & 1.364302   & 1.212493    \\
   & {\textbf{With Vertex }\\\textbf{Morphing filtering }}   & BB Step       & 1.281781   & 1.016281    \\
    &  & Constant Step & 1.281807   & 1.016821    \\ \hline
{\textbf{Interpolation }\\\textbf{techniques}\\\text{(using temperature}\\\text{sensors)}}                        & \textbf{kNN}                                            &               & 1.561738   & 1.518440    \\
   & \textbf{Ordinary kriging}  &  & 1.561738   & 1.650725    \\
       & \textbf{Universal kriging}  &   & 1.561737   & 1.650287   \\ \hline 
\end{tblr}
\end{table}

Table \ref{tab:Table_1} shows the root-mean-squared errors (\acrshort{RMSE}) for the temperature distributions obtained by the proposed optimization-based approach (Figures \ref{fig:Figure_5} and \ref{fig:Figure_6}), \acrshort{kNN} interpolation (Figure \ref{fig:Figure_21}), ordinary kriging (Figure \ref{fig:Figure_24}), and universal kriging (Figure \ref{fig:Figure_25}) compared to the prescribed target distribution for the 6 and 14 sensors configurations (Figure \ref{fig:Figure_3} (right)). The initial \acrshort{RMSE} = $1.561738$ at the beginning of the optimization when $\Delta T = 0$ $\degree C$ everywhere.

Firstly, as presented in Section \ref{sec:Plate With a Hole}, the results obtained from the use of constant stepsize and Barzilai-Borwein (\acrshort{BB}) stepsize for the optimization are almost the same. This is confirmed from the practically identical \acrshort{RMSE} values shown in Table \ref{tab:Table_1}, thus prompting the confident use of \acrshort{BB} step for optimization in this work.

Secondly, it can be observed that in cases using Vertex Morphing filtering, the \acrshort{RMSE} is smaller than when no filtering is used. The RMSE reduces by approximately 6\% and 16\% for the 6 and 14 sensors configurations, respectively, upon incorporation of Vertex Morphing. 

Thirdly, for the 6 temperature sensors configuration used in interpolation methods, since none of the sensors detect a variation in $\Delta T$, the error remains the same as the initial values with any minute difference attributed to the addition of small numerical noise in the measured data. 

Lastly, it can be seen from Table \ref{tab:Table_1} that the proposed optimization-based approach using displacement sensors for temperature field identification performs better than the compared interpolation techniques using temperature sensors in the same configurations. Up to 17.9\% reduction in RMSE (proposed approach with VM and BB step versus all the interpolation techniques considered) for the 6 sensors configuration and 38.4\% reduction in RMSE (proposed approach with VM and BB step versus both the kriging interpolation techniques considered) for the 14 sensors configuration depicts this improvement.
This is also confirmed visually by comparing Figure \ref{fig:Figure_5e} with Figures \ref{fig:Figure_21a},\ref{fig:Figure_24a}, and \ref{fig:Figure_25a} for the 6 sensors configuration and comparing Figure \ref{fig:Figure_6e} with Figures \ref{fig:Figure_24b} and \ref{fig:Figure_25b} for the 14 sensors configuration.

\begin{table}[!t]
\caption{Bridge Example: Comparison of RMSE between the temperature distributions obtained by the different approaches and the prescribed target distribution (shown in Figure \ref{fig:Figure_8} with initial RMSE = 4.743416 i.e., the RMSE at optimization start when $\mathbf{\Delta T} = \mathbf{0}$ $\degree C)$.}
\label{tab:Table_2}
\centering
\begin{tblr}{
  width = \linewidth,
  colspec = {Q[275]Q[258]Q[133]Q[133]},
  cells = {c},
  cell{1}{1} = {c=2}{0.533\linewidth},
  cell{2}{1} = {r=2}{},
  cell{4}{1} = {r=3}{},
  hline{1-2,4,7} = {-}{},
  hline{3,5-6} = {2-4}{},
} \hline
\textbf{Sensor Configuration / Number of Sensors }     &           & \textbf{8 Sensors} & \textbf{20 Sensors} \\
{\textbf{Proposed }\\\textbf{optimization-based }\\\textbf{approach }\\\text{(using displacement}\\\text{sensors)}} & {\textbf{Without Vertex }\\\textbf{Morphing filtering}}    & 1.563510 & 0.199478    \\
    & {\textbf{With Vertex }\\\textbf{Morphing filtering}}       & 1.428890    & 0.150710
    \\ \hline {\textbf{Interpolation }\\\textbf{techniques }\\\text{(using temperature}\\\text{sensors)}}      & \textbf{kNN}            & 3.360823
    & 1.622332
    \\  & \textbf{Ordinary kriging}         & 2.866011   & 2.530294
   \\   & \textbf{Universal kriging}        & 2.367076    & 2.233087
    \\ \hline
\end{tblr}
\end{table}

\medskip

Table \ref{tab:Table_2} shows the root-mean-squared errors (\acrshort{RMSE}) for the temperature distributions obtained by the proposed optimization-based approach (Figures \ref{fig:Figure_11} and \ref{fig:Figure_12}), \acrshort{kNN} interpolation (Figure \ref{fig:Figure_21}), ordinary kriging (Figure \ref{fig:Figure_26}), and universal kriging (Figure \ref{fig:Figure_27}) compared to the prescribed target distribution for the 8 and 20 sensors configurations (Figure \ref{fig:Figure_8}). The initial \acrshort{RMSE} = $4.743416$ at the beginning of the optimization when $\Delta T = 0$ $\degree C$ everywhere.

Similar to the Plate With a Hole example, it can be observed here that the introduction of Vertex Morphing improves (reduces) the RMSE by approximately 8.6\% and 24.4\% for the 8 and 20 sensors configurations, respectively.

For the Bridge example also, it can be observed from Table \ref{tab:Table_2} that the proposed optimization-based approach using displacement sensors for temperature field identification strongly outperforms the compared interpolation techniques using temperature sensors in the same configurations. This improvement is communicated by an up to 57.5\% reduction in RMSE (proposed approach with VM versus kNN interpolation) for the 8 sensors configuration and 94\% reduction in RMSE (proposed approach with VM versus ordinary kriging) for the 20 sensors configuration.
This is also confirmed visually by comparing Figure \ref{fig:Figure_11b} with Figure \ref{fig:Figure_22a} for the 8 sensors configuration and comparing Figure \ref{fig:Figure_12b} with Figure \ref{fig:Figure_26b} for the 20 sensors configuration.

\medskip

Table \ref{tab:Table_3} shows the root-mean-squared errors for the temperature distributions obtained by the proposed optimization-based approach (Figures \ref{fig:Figure_17b}, \ref{fig:Figure_18b}, and \ref{fig:Figure_19b}), kNN interpolation (Figure \ref{fig:Figure_23}), ordinary kriging  (Figure \ref{fig:Figure_28}), and universal kriging  (Figure \ref{fig:Figure_29}) compared to the prescribed target distribution for the 27, 36, and 59 sensors configurations  (Figure \ref{fig:Figure_14}). The initial \acrshort{RMSE} = $3.280578$ at the beginning of the optimization when $\Delta T = 0$ $\degree C$ everywhere.

Similar to the Plate With a Hole example, it can be observed that for cases using Vertex Morphing filtering, the \acrshort{RMSE} is smaller than when no filtering is used. The RMSE reduces by approximately 5.3\%, 17.7\%, and 18\% for the 27, 36, and 59 sensors configurations, respectively, when Vertex Morphing is introduced. 

For the dam example also, it is noted from Table \ref{tab:Table_3} that the proposed optimization-based approach using displacement sensors for temperature field identification outperforms the compared interpolation techniques using temperature sensors in the same configurations. Up to 25\% reduction in RMSE (proposed approach with VM versus kNN interpolation) for the 27 sensors configuration, 40.9\% reduction in RMSE (proposed approach with VM versus ordinary kriging) for the 36 sensors configuration, and 38.8\% reduction in RMSE (proposed approach with VM versus ordinary kriging) for the 59 sensors configuration conveys this improvement.
This is also confirmed visually by comparing Figure \ref{fig:Figure_17b} with Figure \ref{fig:Figure_23a} for the 27 sensors configuration, comparing Figure \ref{fig:Figure_18b} with Figure \ref{fig:Figure_28b} for the 36 sensors configuration and comparing Figure \ref{fig:Figure_19b} with Figure \ref{fig:Figure_28c} for the 59 sensors configuration.

\begin{table}[!bh]
\caption{Dam Example: Comparison of RMSE between the temperature distributions obtained by the different approaches and the prescribed target distribution (shown in Figure \ref{fig:Figure_14} with initial RMSE = 3.280578 i.e., the RMSE at optimization start when $\mathbf{\Delta T} = \mathbf{0}$ $\degree C)$.}
\label{tab:Table_3}
\centering
\begin{tblr}{
  width = \linewidth,
  colspec = {Q[275]Q[258]Q[133]Q[133]Q[133]},
  cells = {c},
  cell{1}{1} = {c=2}{0.533\linewidth},
  cell{2}{1} = {r=2}{},
  cell{4}{1} = {r=3}{},
  hline{1-2,4,7} = {-}{},
  hline{3,5-6} = {2-5}{},
} \hline
\textbf{Sensor Configuration / Number of Sensors }     &            & \textbf{27 Sensors} & \textbf{36 Sensors} & \textbf{59 Sensors} \\
{\textbf{Proposed }\\\textbf{optimization-based }\\\textbf{approach }\\\text{(using displacement}\\\text{sensors)}} & {\textbf{Without Vertex }\\\textbf{Morphing filtering}} & 2.108417    & 1.996023    & 1.777697    \\
    & {\textbf{With Vertex }\\\textbf{Morphing filtering}}    & 1.997094    & 1.642007    & 1.457633    \\ \hline {\textbf{Interpolation }\\\textbf{techniques }\\\text{(using temperature}\\\text{sensors)}}      & \textbf{kNN}         & 2.666437    & 2.471165    & 2.298118    \\  & \textbf{Ordinary kriging}      & 2.482326    & 2.778487    & 2.380026    \\   & \textbf{Universal kriging}     & 2.189399    & 2.453941    & 2.249183    \\ \hline
\end{tblr}
\end{table}

\FloatBarrier
\section{Conclusion}
\label{sec:conclusion}

In this study, an adjoint-based high-fidelity optimization-driven approach to reconstruct the thermal field in structures using displacement/strain measurements was introduced. This physics-included approach mainly focuses on obtaining accurate temperature distributions of structures for use in structural health monitoring and damage detection \& localization. Since the success of this method hinges on deformation sensor measurements, the number and location of sensors, referred to as 'configuration' in this work, should be such that thermal-induced deformations at the region of interest are detected. This removes the explicit requirement for sensors to be located in the region of interest. The methodology was tested using three numerical examples of a Plate With a Hole, a Bridge, and a simplified Hoover Dam FE model. In all cases, a target temperature distribution of the structure was prescribed, together with the external loads. The simulated result was then used to extract the displacements at the sensor locations, thus generating the measurement data. The algorithm started with an initial temperature distribution. It used the measurement data as a reference to tune the temperature field of the structure while aiming to minimize the error between the measured displacement and the computed displacement at the sensor locations accumulated using the cost function. The steepest descent optimization algorithm with the Barzilai-Borwein stepsize calculation method was employed for all three examples. Vertex morphing filtering was employed to mitigate the system's ill-conditioning by filtering the gradients and producing smoother design updates and, thereby, smoother temperature distribution over the structure. Convergence was specified to be achieved when the cost function reaches a 5-magnitude reduction from the initial cost function value.

The Plate with a Hole example was examined for 6 and 14 displacement sensors configurations. Both the configurations satisfactorily reconstructed the overall temperature distribution and localized the heated region, with 14 sensors configuration producing better resolution. The 6 sensors configuration case revealed that the heated area was identified even when none of the displacement sensors were located in the heated region. However, the accuracy in magnitude and localization was low, indicating that more sensors and better placement were required. Without Vertex Morphing, the temperature distributions appeared to have local effects with sharp temperatures and large thermal gradients. Vertex morphing convolutes the gradients over an area, thereby creating smooth temperature changes. 
For comparison, the steepest descent with small constant stepsize cases for both the sensor configurations were also run, and the results were found to be virtually identical to those obtained using the \acrshort{BB} stepsize method. This alleviated concerns regarding the spiking behavior in the cost function convergence plots observed when using the BB step.

The Bridge example was analyzed for 8 and 20 displacement sensors configurations. Similar behavior as the Plate With a Hole example was observed in this case, with both sensor configurations successfully identifying the temperature distribution in the structure and localizing the vicinity of the heated region. Consistent with the Plate With a Hole example, better accuracy and resolution in the temperature distribution were obtained when more sensors are used. In this case, the 20 sensors configuration. Due to the low complexity of the Bridge model compared to the Plate With a Hole and Dam examples, remarkably good thermal fields were obtained even without using Vertex Morphing regularization.

The Hoover Dam example was investigated for 27, 36, and 59 displacement sensors configurations. In line with the Plate With a Hole and the Bridge examples, in this case also, all three sensor configurations sufficiently identified the temperature distribution in the structure and localized the neighborhood of the heated region. The 27 sensors configuration lacked sensors in the lower downstream region of the dam, which resulted in poor temperature identification in this region. An immediate improvement was observed when 9 sensors were added to this sensor deficit region. Consequently, the 36 sensors configuration could identify the general thermal field of the structure and localize the heating over the entire downstream surface as prescribed in the target distribution. The 59 sensors configuration further improved the identification and localization in a sense that both the entire downstream surface as well as the depth of heat penetration, i.e., tangentially inwards to the dam, were more accurately identified. An increase in the number of sensors and the use of Vertex Morphing filtering increased the number of iterations required to converge, with the introduction of Vertex Morphing causing a drastic increase. 
Due to very low sensitivities, especially in regions near corners and edges, some artificial artifacts appearing as large positive and large negative $\Delta T$'s were observed. Vertex morphing smoothing operation tried to reduce the peaks, but some areas still experienced this artificial anomaly. A larger number of sensors, better placement of sensors, testing more load cases to maximize the sensitivity information at all areas of the structure, and regularization are some techniques that can reduce the sensor blind spots and produce a more accurate temperature identification. 

The cost function convergence plots for the Bridge and Hoover Dam examples revealed a trend between the optimization iterations and the reduction in the cost function value. It was observed that a $3-4$ magnitudes reduction in the cost function was achieved fairly quickly within 20\% of the total iterations to converge. The thermal fields at this intermediate progress already contained a sufficient amount of information regarding the rough magnitude and vicinity of the heated and unheated regions. The remaining optimization iterations were utilized in fine-tuning the temperature distribution and dropping the final 1-2 magnitudes of the cost function. Thus, depending on the required precision of the thermal field, the convergence criteria could be changed to balance accuracy and computational effort.

In all cases, the exact prescribed target distribution was not obtained, the peak temperatures were mostly slightly lower than the target peak, and the identified region was moderately larger (i.e., diffused) than the prescribed target distribution. Nevertheless, the distributions were very close to the target distributions. 

For a comparative study, the same sensor configurations were used with temperatures extracted from the target distribution at the measurement locations as input, i.e., temperature sensors. The thermal fields were reconstructed using k-nearest neighbor and kriging interpolation techniques. The spatial interpolation techniques produced inferior thermal field results compared to those obtained from the proposed approach. As expected, in the absence of temperature sensors from the region of interest, the temperature in that region was highly inaccurate. 
In general, the kriging interpolated thermal fields were smoother than the ones produced by \acrshort{kNN} interpolation. Universal kriging performed better than ordinary kriging due to considering a trend in the mean of the spatial process, analogous to the thermal gradient in the structure. However, universal kriging results suffered from negative $\Delta T$'s and, poor localization in the tangential direction to the downstream surface in the Dam example.

It is noteworthy here that although the peak $\Delta T$ identified in the Plate With a Hole example with 6 sensors configuration using the proposed approach was only around a third of the prescribed value, it is still dramatically better than the interpolation techniques, which do not detect any variation due to the scarcity and placement of sensors. 
On the other hand, in the Plate With a Hole example with 14 sensors configuration, even though the interpolation techniques identified peak $\Delta T$ in some areas identical to that in the prescribed distribution, the identified location of the 'heated-zone' was rightward offset and centered around the sensor instead of the region of interest. However, the results obtained from the proposed approach generated a superior localization of the area of interest with some compromise in the peak $\Delta T$. 

The root-mean-squared errors (\acrshort{RMSE}) between the synthetic prescribed temperature distribution and the identified/interpolated temperature distributions were also analyzed. The RMSE values for the three examples and different sensor configurations were consistent with the behavior observed from visual inspection. In a nutshell, the RMSEs from the proposed optimization-based approach with and without Vertex Morphing were lower than the RMSEs from the interpolation techniques. RMSE reductions of up to 38.4\%, 94\%, and 40\% were seen in the Plate With a Hole, the Bridge, and the Hoover Dam examples respectively. Within the proposed approach, the cases with Vertex Morphing registered lower RMSE than the cases not using Vertex Morphing. Improvements as high as 16\%, 24.4\%, and 18\% in the RMSEs were observed in the Plate With a Hole, the Bridge, and the Hoover Dam examples, respectively, when Vertex Morphing was used.

\medskip

Although the proposed approach successfully reconstructs the temperature distribution in the structures, there are several open questions that shall be addressed in future works:
\begin{itemize}
    \item Consideration of noise in the sensor measurements (\cite{FAiraudo_HAntil_RLoehner_URakhimov_2024a}),
    \item Uncertainty in the parameters such as material parameters, sensor locations, etc (\cite{HAntil_SDolgov_AOnwunta_2022b,airaudo2024use}),
    \item Elimination or reduction of sensor 'blind-spots' using better sensor configuration (\cite{warnakulasuriya2024optimal}),
    \item Consideration of other ambient parameters like wind, rain, cloud cover, humidity, angle of irradiation, etc,
    \item Consideration of time lag in the effects observed within the structure due to current ambient conditions, especially relevant in massive concrete structures,
    \item Combining heterogeneous (displacement/strain and temperature sensors) sensor information for a realistic and feedback-driven structure status.
\end{itemize}

\FloatBarrier

\section*{Acknowledgments}

We would like to acknowledge the funding by the Institute of Advanced Studies, Technical University of Munich, under the Hans Fischer Senior Fellowship. The second, fourth, and final authors are also partially supported by NSF grant DMS-2408877, Air Force Office of Scientific Research (AFOSR) under Award NO: FA9550-22-1-0248, and Office of Naval Research (ONR) under Award NO: N00014-24-1-2147. We would also like to acknowledge the funding by the Deutsche Forschungsgemeinschaft (DFG, German Research Foundation) under German's Excellence Strategy - EXC 2163/1 - Sustainable and Energy Efficient Aviation - Project-ID 390881007.

\bibliographystyle{elsarticle-harv} 
\bibliography{bibliography}

\begin{thebibliography}{71}
\expandafter\ifx\csname natexlab\endcsname\relax\def\natexlab#1{#1}\fi
\providecommand{\url}[1]{\texttt{#1}}
\providecommand{\href}[2]{#2}
\providecommand{\path}[1]{#1}
\providecommand{\DOIprefix}{DOI:}
\providecommand{\ArXivprefix}{arXiv:}
\providecommand{\URLprefix}{URL: }
\providecommand{\Pubmedprefix}{pmid:}
\providecommand{\doi}[1]{\href{http://dx.doi.org/#1}{\path{#1}}}
\providecommand{\Pubmed}[1]{\href{pmid:#1}{\path{#1}}}
\providecommand{\bibinfo}[2]{#2}
\ifx\xfnm\relax \def\xfnm[#1]{\unskip,\space#1}\fi
\bibitem[{Abid et~al.(2016)Abid, Tay{\c{s}}i and
  {\"O}zak{\c{c}}a}]{abid2016experimental}
\bibinfo{author}{Abid, S.R.}, \bibinfo{author}{Tay{\c{s}}i, N.},
  \bibinfo{author}{{\"O}zak{\c{c}}a, M.}, \bibinfo{year}{2016}.
\newblock \bibinfo{title}{Experimental analysis of temperature gradients in
  concrete box-girders}.
\newblock \bibinfo{journal}{Construction and Building Materials}
  \bibinfo{volume}{106}, \bibinfo{pages}{523--532}.
\newblock \DOIprefix\doi{https://doi.org/10.1016/j.conbuildmat.2015.12.144}.
\bibitem[{Abid et~al.(2022)Abid, Xue, Liu, Tay{\c{s}}i, Liu, {\"O}zak{\c{c}}a
  and Briseghella}]{abid2022temperatures}
\bibinfo{author}{Abid, S.R.}, \bibinfo{author}{Xue, J.}, \bibinfo{author}{Liu,
  J.}, \bibinfo{author}{Tay{\c{s}}i, N.}, \bibinfo{author}{Liu, Y.},
  \bibinfo{author}{{\"O}zak{\c{c}}a, M.}, \bibinfo{author}{Briseghella, B.},
  \bibinfo{year}{2022}.
\newblock \bibinfo{title}{Temperatures and gradients in concrete bridges:
  Experimental, finite element analysis and design}, in:
  \bibinfo{booktitle}{Structures}, \bibinfo{organization}{Elsevier}. pp.
  \bibinfo{pages}{960--976}.
\newblock \DOIprefix\doi{https://doi.org/10.1016/j.istruc.2022.01.070}.
\bibitem[{Airaudo et~al.(2024a)Airaudo, Antil, L\"ohner and
  Rakhimov}]{FAiraudo_HAntil_RLoehner_URakhimov_2024a}
\bibinfo{author}{Airaudo, F.}, \bibinfo{author}{Antil, H.},
  \bibinfo{author}{L\"ohner, R.}, \bibinfo{author}{Rakhimov, U.},
  \bibinfo{year}{2024}a.
\newblock \bibinfo{title}{On the use of risk measures in digital twins to
  identify weaknesses in structures}, in: \bibinfo{booktitle}{AIAA SCITECH 2024
  Forum}, \bibinfo{publisher}{American Institute of Aeronautics and
  Astronautics}. p. \bibinfo{pages}{2622}.
\newblock \DOIprefix\doi{http://dx.doi.org/10.2514/6.2024-2622}.
\bibitem[{Airaudo et~al.(2024b)Airaudo, Antil, L{\"o}hner, Warnakulasuriya,
  Antonau and W{\"u}chner}]{airaudo2024use}
\bibinfo{author}{Airaudo, F.}, \bibinfo{author}{Antil, H.},
  \bibinfo{author}{L{\"o}hner, R.}, \bibinfo{author}{Warnakulasuriya, S.},
  \bibinfo{author}{Antonau, I.}, \bibinfo{author}{W{\"u}chner, R.},
  \bibinfo{year}{2024}b.
\newblock \bibinfo{title}{On the use of risk measures in digital twins to
  identify weaknesses in structures}, in: \bibinfo{booktitle}{ECCOMAS Congress
  2024: 9th European Congress on Computational Methods in Applied Sciences and
  Engineering}, \bibinfo{address}{Lisboa, Portugal}.
\bibitem[{Airaudo et~al.(2023)Airaudo, L{\"o}hner, W{\"u}chner and
  Antil}]{airaudo2023adjoint}
\bibinfo{author}{Airaudo, F.N.}, \bibinfo{author}{L{\"o}hner, R.},
  \bibinfo{author}{W{\"u}chner, R.}, \bibinfo{author}{Antil, H.},
  \bibinfo{year}{2023}.
\newblock \bibinfo{title}{Adjoint-based determination of weaknesses in
  structures}.
\newblock \bibinfo{journal}{Computer Methods in Applied Mechanics and
  Engineering} \bibinfo{volume}{417}, \bibinfo{pages}{116471}.
\newblock \DOIprefix\doi{https://doi.org/10.1016/j.cma.2023.116471}.
\bibitem[{{American Association of State Highway and Transportation
  Officials}(2020)}]{aashto2020lrfd}
\bibinfo{author}{{American Association of State Highway and Transportation
  Officials}}, \bibinfo{year}{2020}.
\newblock \bibinfo{title}{LRFD Bridge Design Specifications}.
  \bibinfo{edition}{9th} ed.
\newblock \bibinfo{address}{Washington, DC}.
\newblock ISBN: 1560517387. ISBN (Online): 9781560517382, \URLprefix
  \url{https://aashtojournal.transportation.org/aashto-issues-updated-lrfd-bridge-design-guide/}.
  \bibinfo{note}{({A}ccessed online on 28/10/2024)}.
\bibitem[{Antil(2024)}]{HAntil_2024a}
\bibinfo{author}{Antil, H.}, \bibinfo{year}{2024}.
\newblock \bibinfo{title}{Mathematical opportunities in digital twins
  (math-dt)}.
\newblock \URLprefix \url{https://arxiv.org/abs/2402.10326},
  \href{http://arxiv.org/abs/2402.10326}{{\tt arXiv:2402.10326}}.
\bibitem[{Antil et~al.(2022)Antil, Dolgov and
  Onwunta}]{HAntil_SDolgov_AOnwunta_2022b}
\bibinfo{author}{Antil, H.}, \bibinfo{author}{Dolgov, S.},
  \bibinfo{author}{Onwunta, A.}, \bibinfo{year}{2022}.
\newblock \bibinfo{title}{Ttrisk: Tensor train decomposition algorithm for risk
  averse optimization}.
\newblock \bibinfo{journal}{Numerical Linear Algebra with Applications}
  \bibinfo{volume}{30}.
\newblock \DOIprefix\doi{https://doi.org/10.1002/nla.2481}.
\bibitem[{Antil et~al.(2018)Antil, Kouri, Lacasse and
  Ridzal}]{HAntil_DPKouri_MDLacasse_DRidzal_2018a}
\bibinfo{editor}{Antil, H.}, \bibinfo{editor}{Kouri, D.P.},
  \bibinfo{editor}{Lacasse, M.D.}, \bibinfo{editor}{Ridzal, D.} (Eds.),
  \bibinfo{year}{2018}.
\newblock \bibinfo{title}{Frontiers in {PDE}-constrained optimization}. volume
  \bibinfo{volume}{163} of \textit{\bibinfo{series}{The IMA Volumes in
  Mathematics and its Applications}}.
\newblock \bibinfo{publisher}{Springer, New York}.
\newblock ISBN: 978-1-4939-8635-4; 978-1-4939-8636-1,
  \DOIprefix\doi{https://doi.org/10.1007/978-1-4939-8636-1}.
  \bibinfo{note}{{P}apers based on the workshop held at the Institute for
  Mathematics and its Applications, Minneapolis, MN, June 6--10, 2016}.
\bibitem[{Antonau(2023)}]{antonau2023enhanced}
\bibinfo{author}{Antonau, I.}, \bibinfo{year}{2023}.
\newblock \bibinfo{title}{Enhanced computational design methods for large
  industrial node-based shape optimization problems}.
\newblock Ph.D. thesis. Technische Universit{\"a}t M{\"u}nchen.
\newblock \URLprefix \url{https://mediatum.ub.tum.de/1695301}.
\bibitem[{Antonau et~al.(2024)Antonau, Warnakulasuriya, Airaudo, L{\"o}hner,
  Antil and W{\"u}chner}]{antonau2024formulation}
\bibinfo{author}{Antonau, I.}, \bibinfo{author}{Warnakulasuriya, S.},
  \bibinfo{author}{Airaudo, F.}, \bibinfo{author}{L{\"o}hner, R.},
  \bibinfo{author}{Antil, H.}, \bibinfo{author}{W{\"u}chner, R.},
  \bibinfo{year}{2024}.
\newblock \bibinfo{title}{Formulation and regularization of the optimization
  problem for adjoint-based digital twin construction}, in:
  \bibinfo{booktitle}{ECCOMAS Congress 2024: 9th European Congress on
  Computational Methods in Applied Sciences and Engineering},
  \bibinfo{address}{Lisbon, Portugal}.
\bibitem[{Antonau et~al.(2022)Antonau, Warnakulasuriya, Bletzinger, Bluhm,
  Hojjat and W{\"u}chner}]{antonau2022latest}
\bibinfo{author}{Antonau, I.}, \bibinfo{author}{Warnakulasuriya, S.},
  \bibinfo{author}{Bletzinger, K.U.}, \bibinfo{author}{Bluhm, F.M.},
  \bibinfo{author}{Hojjat, M.}, \bibinfo{author}{W{\"u}chner, R.},
  \bibinfo{year}{2022}.
\newblock \bibinfo{title}{Latest developments in node-based shape optimization
  using vertex morphing parameterization}.
\newblock \bibinfo{journal}{Structural and Multidisciplinary Optimization}
  \bibinfo{volume}{65}, \bibinfo{pages}{198}.
\newblock \DOIprefix\doi{https://doi.org/10.1007/s00158-022-03279-w}.
\bibitem[{B{\k{a}}czkiewicz et~al.(2018)B{\k{a}}czkiewicz, Malaska, Pajunen and
  Heinisuo}]{bkaczkiewicz2018experimental}
\bibinfo{author}{B{\k{a}}czkiewicz, J.}, \bibinfo{author}{Malaska, M.},
  \bibinfo{author}{Pajunen, S.}, \bibinfo{author}{Heinisuo, M.},
  \bibinfo{year}{2018}.
\newblock \bibinfo{title}{Experimental and numerical study on temperature
  distribution of square hollow section joints}.
\newblock \bibinfo{journal}{Journal of Constructional Steel Research}
  \bibinfo{volume}{142}, \bibinfo{pages}{31--43}.
\newblock \DOIprefix\doi{https://doi.org/10.1016/j.jcsr.2017.12.006}.
\bibitem[{Bayraktar et~al.(2022)Bayraktar, Akk{\"o}se, Ta{\c{s}}, Erdi{\c{s}}
  and Kur{\c{s}}un}]{bayraktar2022long}
\bibinfo{author}{Bayraktar, A.}, \bibinfo{author}{Akk{\"o}se, M.},
  \bibinfo{author}{Ta{\c{s}}, Y.}, \bibinfo{author}{Erdi{\c{s}}, A.},
  \bibinfo{author}{Kur{\c{s}}un, A.}, \bibinfo{year}{2022}.
\newblock \bibinfo{title}{Long-term strain behavior of in-service cable-stayed
  bridges under temperature variations}.
\newblock \bibinfo{journal}{Journal of Civil Structural Health Monitoring}
  \bibinfo{volume}{12}, \bibinfo{pages}{833--844}.
\newblock \DOIprefix\doi{https://doi.org/10.1007/s13349-022-00578-0}.
\bibitem[{Bui et~al.(2019)Bui, Sancharoen, Tanapornraweekit, Tangtermsirikul
  and Nanakorn}]{bui2019evaluation}
\bibinfo{author}{Bui, K.A.}, \bibinfo{author}{Sancharoen, P.},
  \bibinfo{author}{Tanapornraweekit, G.}, \bibinfo{author}{Tangtermsirikul,
  S.}, \bibinfo{author}{Nanakorn, P.}, \bibinfo{year}{2019}.
\newblock \bibinfo{title}{An evaluation of thermal effects on behavior of a
  concrete arch dam}.
\newblock \bibinfo{journal}{Songklanakarin Journal of Science and Technology}
  \bibinfo{volume}{41}, \bibinfo{pages}{1059--1068}.
\bibitem[{Castilho et~al.(2018)Castilho, Schclar, Tiago and
  Farinha}]{castilho2018fea}
\bibinfo{author}{Castilho, E.}, \bibinfo{author}{Schclar, N.},
  \bibinfo{author}{Tiago, C.}, \bibinfo{author}{Farinha, M.L.B.},
  \bibinfo{year}{2018}.
\newblock \bibinfo{title}{{FEA} model for the simulation of the hydration
  process and temperature evolution during the concreting of an arch dam}.
\newblock \bibinfo{journal}{Engineering Structures} \bibinfo{volume}{174},
  \bibinfo{pages}{165--177}.
\newblock \DOIprefix\doi{https://doi.org/10.1016/j.engstruct.2018.07.065}.
\bibitem[{Chen and Liu(2019)}]{chen2019temperature}
\bibinfo{author}{Chen, H.}, \bibinfo{author}{Liu, Z.}, \bibinfo{year}{2019}.
\newblock \bibinfo{title}{Temperature control and thermal-induced stress field
  analysis of {GongGuoQiao} {RCC} dam}.
\newblock \bibinfo{journal}{Journal of Thermal Analysis and Calorimetry}
  \bibinfo{volume}{135}, \bibinfo{pages}{2019--2029}.
\newblock \DOIprefix\doi{https://doi.org/10.1007/s10973-018-7450-1}.
\bibitem[{Corbally and Malekjafarian(2022)}]{corbally2022data}
\bibinfo{author}{Corbally, R.}, \bibinfo{author}{Malekjafarian, A.},
  \bibinfo{year}{2022}.
\newblock \bibinfo{title}{A data-driven approach for drive-by damage detection
  in bridges considering the influence of temperature change}.
\newblock \bibinfo{journal}{Engineering Structures} \bibinfo{volume}{253},
  \bibinfo{pages}{113783}.
\newblock \DOIprefix\doi{https://doi.org/10.1016/j.engstruct.2021.113783}.
\bibitem[{Dadvand et~al.(2010)Dadvand, Rossi and O{\~n}ate}]{dadvand2010object}
\bibinfo{author}{Dadvand, P.}, \bibinfo{author}{Rossi, R.},
  \bibinfo{author}{O{\~n}ate, E.}, \bibinfo{year}{2010}.
\newblock \bibinfo{title}{An object-oriented environment for developing finite
  element codes for multi-disciplinary applications}.
\newblock \bibinfo{journal}{Archives of computational methods in engineering}
  \bibinfo{volume}{17}, \bibinfo{pages}{253--297}.
\newblock \DOIprefix\doi{https://doi.org/10.1007/s11831-010-9045-2}.
\bibitem[{Fan et~al.(2022)Fan, Li, Liu and Liu}]{fan2022efficient}
\bibinfo{author}{Fan, J.S.}, \bibinfo{author}{Li, B.L.}, \bibinfo{author}{Liu,
  C.}, \bibinfo{author}{Liu, Y.F.}, \bibinfo{year}{2022}.
\newblock \bibinfo{title}{An efficient model for simulation of temperature
  field of steel-concrete composite beam bridges}.
\newblock \bibinfo{journal}{Structures} \bibinfo{volume}{43},
  \bibinfo{pages}{1868--1880}.
\newblock \DOIprefix\doi{https://doi.org/10.1016/j.istruc.2022.05.079}.
\bibitem[{Ferr\'{a}ndiz et~al.(2024)Ferr\'{a}ndiz, Bucher, Zorrilla,
  Warnakulasuriya, Rossi, Cornejo, jcotela, Roig, Maria, tteschemacher,
  Mas\'{o}, Casas, Núñez, Dadvand, Latorre, de~Pouplana, Gonz\'{a}lez,
  AFranci, Arrufat, riccardotosi, Ghantasala, Sautter, Wilson, dbaumgaertner,
  Chandra, Geiser, Lopez, llu\'{i}s, jgonzalezusua and
  G\'{a}rate}]{vicente_mataix_ferrandiz_2024_6926179}
\bibinfo{author}{Ferr\'{a}ndiz, V.M.}, \bibinfo{author}{Bucher, P.},
  \bibinfo{author}{Zorrilla, R.}, \bibinfo{author}{Warnakulasuriya, S.},
  \bibinfo{author}{Rossi, R.}, \bibinfo{author}{Cornejo, A.},
  \bibinfo{author}{jcotela}, \bibinfo{author}{Roig, C.},
  \bibinfo{author}{Maria, J.}, \bibinfo{author}{tteschemacher},
  \bibinfo{author}{Mas\'{o}, M.}, \bibinfo{author}{Casas, G.},
  \bibinfo{author}{Núñez, M.}, \bibinfo{author}{Dadvand, P.},
  \bibinfo{author}{Latorre, S.}, \bibinfo{author}{de~Pouplana, I.},
  \bibinfo{author}{Gonz\'{a}lez, J.I.}, \bibinfo{author}{AFranci},
  \bibinfo{author}{Arrufat, F.}, \bibinfo{author}{riccardotosi},
  \bibinfo{author}{Ghantasala, A.}, \bibinfo{author}{Sautter, K.B.},
  \bibinfo{author}{Wilson, P.}, \bibinfo{author}{dbaumgaertner},
  \bibinfo{author}{Chandra, B.}, \bibinfo{author}{Geiser, A.},
  \bibinfo{author}{Lopez, I.}, \bibinfo{author}{llu\'{i}s},
  \bibinfo{author}{jgonzalezusua}, \bibinfo{author}{G\'{a}rate, J.},
  \bibinfo{year}{2024}.
\newblock \bibinfo{title}{Kratosmultiphysics/kratos: Release v9.5}.
\newblock \DOIprefix\doi{https://doi.org/10.5281/zenodo.6926179}.
\bibitem[{Fletcher(2005)}]{fletcher2005barzilai}
\bibinfo{author}{Fletcher, R.}, \bibinfo{year}{2005}.
\newblock \bibinfo{title}{On the barzilai-borwein method}, in:
  \bibinfo{booktitle}{Optimization and Control with Applications},
  \bibinfo{organization}{Springer}. pp. \bibinfo{pages}{235--256}.
\newblock \DOIprefix\doi{https://doi.org/10.1007/0-387-24255-4_10}.
\bibitem[{Ghantasala et~al.(2021)Ghantasala, Najian~Asl, Geiser, Brodie,
  Papoutsis and Bletzinger}]{ghantasala2021realization}
\bibinfo{author}{Ghantasala, A.}, \bibinfo{author}{Najian~Asl, R.},
  \bibinfo{author}{Geiser, A.}, \bibinfo{author}{Brodie, A.},
  \bibinfo{author}{Papoutsis, E.}, \bibinfo{author}{Bletzinger, K.U.},
  \bibinfo{year}{2021}.
\newblock \bibinfo{title}{Realization of a framework for simulation-based
  large-scale shape optimization using vertex morphing}.
\newblock \bibinfo{journal}{Journal of Optimization Theory and Applications}
  \bibinfo{volume}{189}, \bibinfo{pages}{164--189}.
\newblock \DOIprefix\doi{https://doi.org/10.1007/s10957-021-01826-x}.
\bibitem[{Glashier et~al.(2024a)Glashier, Kromanis and
  Buchanan}]{glashier2024iterative}
\bibinfo{author}{Glashier, T.}, \bibinfo{author}{Kromanis, R.},
  \bibinfo{author}{Buchanan, C.}, \bibinfo{year}{2024}a.
\newblock \bibinfo{title}{An iterative regression-based thermal response
  prediction methodology for instrumented civil infrastructure}.
\newblock \bibinfo{journal}{Advanced Engineering Informatics}
  \bibinfo{volume}{60}, \bibinfo{pages}{102347}.
\newblock \DOIprefix\doi{https://doi.org/10.1016/j.aei.2023.102347}.
\bibitem[{Glashier et~al.(2024b)Glashier, Kromanis and
  Buchanan}]{glashier2024temperature}
\bibinfo{author}{Glashier, T.}, \bibinfo{author}{Kromanis, R.},
  \bibinfo{author}{Buchanan, C.}, \bibinfo{year}{2024}b.
\newblock \bibinfo{title}{Temperature-based measurement interpretation of the
  {MX3D} {B}ridge}.
\newblock \bibinfo{journal}{Engineering Structures} \bibinfo{volume}{305},
  \bibinfo{pages}{116736}.
\newblock \DOIprefix\doi{https://doi.org/10.1016/j.engstruct.2023.116736}.
\bibitem[{Gu et~al.(2017)Gu, Gul and Wu}]{gu2017damage}
\bibinfo{author}{Gu, J.}, \bibinfo{author}{Gul, M.}, \bibinfo{author}{Wu, X.},
  \bibinfo{year}{2017}.
\newblock \bibinfo{title}{Damage detection under varying temperature using
  artificial neural networks}.
\newblock \bibinfo{journal}{Structural Control and Health Monitoring}
  \bibinfo{volume}{24}, \bibinfo{pages}{e1998}.
\newblock \DOIprefix\doi{https://doi.org/10.1002/stc.1998}.
\bibitem[{Hagedorn et~al.(2019)Hagedorn, Mart{\'\i}-Vargas, Dang, Hale and
  Floyd}]{hagedorn2019temperature}
\bibinfo{author}{Hagedorn, R.}, \bibinfo{author}{Mart{\'\i}-Vargas, J.R.},
  \bibinfo{author}{Dang, C.N.}, \bibinfo{author}{Hale, W.M.},
  \bibinfo{author}{Floyd, R.W.}, \bibinfo{year}{2019}.
\newblock \bibinfo{title}{Temperature gradients in bridge concrete {I}-girders
  under heat wave}.
\newblock \bibinfo{journal}{Journal of Bridge Engineering}
  \bibinfo{volume}{24}.
\newblock \DOIprefix\doi{https://doi.org/10.1061/(asce)be.1943-5592.0001454}.
\bibitem[{Hojjat et~al.(2014)Hojjat, Stavropoulou and
  Bletzinger}]{hojjat2014vertex}
\bibinfo{author}{Hojjat, M.}, \bibinfo{author}{Stavropoulou, E.},
  \bibinfo{author}{Bletzinger, K.U.}, \bibinfo{year}{2014}.
\newblock \bibinfo{title}{The {Vertex} {Morphing} method for node-based shape
  optimization}.
\newblock \bibinfo{journal}{Computer Methods in Applied Mechanics and
  Engineering} \bibinfo{volume}{268}, \bibinfo{pages}{494--513}.
\newblock \DOIprefix\doi{https://doi.org/10.1016/j.cma.2013.10.015}.
\bibitem[{Hou et~al.(2022)Hou, Chai, Cheng, Ning, Yang and
  Zhou}]{hou2022simulation}
\bibinfo{author}{Hou, C.}, \bibinfo{author}{Chai, D.}, \bibinfo{author}{Cheng,
  H.}, \bibinfo{author}{Ning, S.}, \bibinfo{author}{Yang, B.},
  \bibinfo{author}{Zhou, Y.}, \bibinfo{year}{2022}.
\newblock \bibinfo{title}{Simulation feedback of temperature field of
  super-high arch dam during operation and its difference with design
  temperature}.
\newblock \bibinfo{journal}{Water} \bibinfo{volume}{14}, \bibinfo{pages}{4028}.
\newblock \DOIprefix\doi{https://doi.org/10.3390/w14244028}.
\bibitem[{Jin et~al.(2016)Jin, Jang, Sun, Li and Christenson}]{jin2016damage}
\bibinfo{author}{Jin, C.}, \bibinfo{author}{Jang, S.}, \bibinfo{author}{Sun,
  X.}, \bibinfo{author}{Li, J.}, \bibinfo{author}{Christenson, R.},
  \bibinfo{year}{2016}.
\newblock \bibinfo{title}{Damage detection of a highway bridge under severe
  temperature changes using extended kalman filter trained neural network}.
\newblock \bibinfo{journal}{Journal of Civil Structural Health Monitoring}
  \bibinfo{volume}{6}, \bibinfo{pages}{545--560}.
\newblock \DOIprefix\doi{https://doi.org/10.1007/s13349-016-0173-8}.
\bibitem[{Jin et~al.(2015)Jin, Li, Jang, Sun and
  Christenson}]{jin2015structural}
\bibinfo{author}{Jin, C.}, \bibinfo{author}{Li, J.}, \bibinfo{author}{Jang,
  S.}, \bibinfo{author}{Sun, X.}, \bibinfo{author}{Christenson, R.},
  \bibinfo{year}{2015}.
\newblock \bibinfo{title}{Structural damage detection for in-service highway
  bridge under operational and environmental variability}, in:
  \bibinfo{booktitle}{Sensors and {S}mart {S}tructures {T}echnologies for
  {C}ivil, {M}echanical, and {A}erospace {S}ystems 2015},
  \bibinfo{publisher}{SPIE}. pp. \bibinfo{pages}{943--952}.
\newblock \DOIprefix\doi{https://doi.org/10.1117/12.2084384}.
\bibitem[{Jin et~al.(2010)Jin, Chen, Wang and Yang}]{jin2010practical}
\bibinfo{author}{Jin, F.}, \bibinfo{author}{Chen, Z.}, \bibinfo{author}{Wang,
  J.}, \bibinfo{author}{Yang, J.}, \bibinfo{year}{2010}.
\newblock \bibinfo{title}{Practical procedure for predicting non-uniform
  temperature on the exposed face of arch dams}.
\newblock \bibinfo{journal}{Applied Thermal Engineering} \bibinfo{volume}{30},
  \bibinfo{pages}{2146--2156}.
\newblock \DOIprefix\doi{https://doi.org/10.1016/j.applthermaleng.2010.05.027}.
\bibitem[{Kosti{\'c} and G{\"u}l(2017)}]{kostic2017vibration}
\bibinfo{author}{Kosti{\'c}, B.}, \bibinfo{author}{G{\"u}l, M.},
  \bibinfo{year}{2017}.
\newblock \bibinfo{title}{Vibration-based damage detection of bridges under
  varying temperature effects using time-series analysis and artificial neural
  networks}.
\newblock \bibinfo{journal}{Journal of Bridge Engineering}
  \bibinfo{volume}{22}.
\newblock \DOIprefix\doi{https://doi.org/10.1061/(asce)be.1943-5592.0001085}.
\bibitem[{Kromanis and Kripakaran(2013)}]{kromanis2013support}
\bibinfo{author}{Kromanis, R.}, \bibinfo{author}{Kripakaran, P.},
  \bibinfo{year}{2013}.
\newblock \bibinfo{title}{Support vector regression for anomaly detection from
  measurement histories}.
\newblock \bibinfo{journal}{Advanced Engineering Informatics}
  \bibinfo{volume}{27}, \bibinfo{pages}{486--495}.
\newblock \DOIprefix\doi{https://doi.org/10.1016/j.aei.2013.03.002}.
\bibitem[{Kromanis and Kripakaran(2014)}]{kromanis2014predicting}
\bibinfo{author}{Kromanis, R.}, \bibinfo{author}{Kripakaran, P.},
  \bibinfo{year}{2014}.
\newblock \bibinfo{title}{Predicting thermal response of bridges using
  regression models derived from measurement histories}.
\newblock \bibinfo{journal}{Computers \& Structures} \bibinfo{volume}{136},
  \bibinfo{pages}{64--77}.
\newblock \DOIprefix\doi{https://doi.org/10.1016/j.compstruc.2014.01.026}.
\bibitem[{Kulprapha and Warnitchai(2012)}]{kulprapha2012structural}
\bibinfo{author}{Kulprapha, N.}, \bibinfo{author}{Warnitchai, P.},
  \bibinfo{year}{2012}.
\newblock \bibinfo{title}{Structural health monitoring of continuous
  prestressed concrete bridges using ambient thermal responses}.
\newblock \bibinfo{journal}{Engineering Structures} \bibinfo{volume}{40},
  \bibinfo{pages}{20--38}.
\newblock \DOIprefix\doi{https://doi.org/10.1016/j.engstruct.2012.02.001}.
\bibitem[{Lee and Kalkan(2012)}]{lee2012analysis}
\bibinfo{author}{Lee, J.H.}, \bibinfo{author}{Kalkan, I.},
  \bibinfo{year}{2012}.
\newblock \bibinfo{title}{Analysis of thermal environmental effects on precast,
  prestressed concrete bridge girders: temperature differentials and thermal
  deformations}.
\newblock \bibinfo{journal}{Advances in Structural Engineering}
  \bibinfo{volume}{15}, \bibinfo{pages}{447--459}.
\newblock \DOIprefix\doi{https://doi.org/10.1260/1369-4332.15.3.447}.
\bibitem[{Lin et~al.(2021)Lin, Peng, Fan, Xiang, Yang and Yang}]{lin20213d}
\bibinfo{author}{Lin, P.}, \bibinfo{author}{Peng, H.}, \bibinfo{author}{Fan,
  Q.}, \bibinfo{author}{Xiang, Y.}, \bibinfo{author}{Yang, Z.},
  \bibinfo{author}{Yang, N.}, \bibinfo{year}{2021}.
\newblock \bibinfo{title}{A 3{D} thermal field restructuring method for
  concrete dams based on real-time temperature monitoring}.
\newblock \bibinfo{journal}{KSCE Journal of Civil Engineering}
  \bibinfo{volume}{25}, \bibinfo{pages}{1326--1340}.
\newblock \DOIprefix\doi{https://doi.org/10.1007/s12205-021-1084-8}.
\bibitem[{Liu et~al.(2015)Liu, Zhang, Chang, Zhou, Cheng and
  Duan}]{liu2015precise}
\bibinfo{author}{Liu, X.}, \bibinfo{author}{Zhang, C.}, \bibinfo{author}{Chang,
  X.}, \bibinfo{author}{Zhou, W.}, \bibinfo{author}{Cheng, Y.},
  \bibinfo{author}{Duan, Y.}, \bibinfo{year}{2015}.
\newblock \bibinfo{title}{Precise simulation analysis of the thermal field in
  mass concrete with a pipe water cooling system}.
\newblock \bibinfo{journal}{Applied Thermal Engineering} \bibinfo{volume}{78},
  \bibinfo{pages}{449--459}.
\newblock \DOIprefix\doi{https://doi.org/10.1016/j.applthermaleng.2014.12.050}.
\bibitem[{Lohner et~al.(2024a)Lohner, Airaudo, Antil, W{\"u}chner, Meister and
  Warnakulasuriya}]{lohner2024high}
\bibinfo{author}{Lohner, R.}, \bibinfo{author}{Airaudo, F.},
  \bibinfo{author}{Antil, H.}, \bibinfo{author}{W{\"u}chner, R.},
  \bibinfo{author}{Meister, F.}, \bibinfo{author}{Warnakulasuriya, S.},
  \bibinfo{year}{2024}a.
\newblock \bibinfo{title}{High-fidelity digital twins: Detecting and localizing
  weaknesses in structures}, in: \bibinfo{booktitle}{AIAA SCITECH 2024 Forum},
  p. \bibinfo{pages}{2621}.
\newblock \DOIprefix\doi{https://doi.org/10.2514/6.2024-2621}.
\bibitem[{Lohner et~al.(2024b)Lohner, Airaudo, Antil, W{\"u}chner, Meister and
  Warnakulasuriya}]{RLoehner_Fairaudo_HAntil_RWuechner_FMeister_SWarnakulasuriya_2024b}
\bibinfo{author}{Lohner, R.}, \bibinfo{author}{Airaudo, F.},
  \bibinfo{author}{Antil, H.}, \bibinfo{author}{W{\"u}chner, R.},
  \bibinfo{author}{Meister, F.}, \bibinfo{author}{Warnakulasuriya, S.},
  \bibinfo{year}{2024}b.
\newblock \bibinfo{title}{High-fidelity digital twins: Detecting and localizing
  weaknesses in structures}.
\newblock \bibinfo{journal}{International Journal for Numerical Methods in
  Engineering (IJNME)} \bibinfo{volume}{125}, \bibinfo{pages}{e7568}.
\newblock \DOIprefix\doi{https://doi.org/10.1002/nme.7568}.
\bibitem[{Ma et~al.(2023)Ma, Wu, Qin, Zhao and Yang}]{ma2023statistical}
\bibinfo{author}{Ma, W.}, \bibinfo{author}{Wu, B.}, \bibinfo{author}{Qin, D.},
  \bibinfo{author}{Zhao, B.}, \bibinfo{author}{Yang, X.}, \bibinfo{year}{2023}.
\newblock \bibinfo{title}{Statistical analyses of the non-uniform longitudinal
  temperature distribution in steel box girder bridge}.
\newblock \bibinfo{journal}{Buildings} \bibinfo{volume}{13},
  \bibinfo{pages}{1316}.
\newblock \DOIprefix\doi{https://doi.org/10.3390/buildings13051316}.
\bibitem[{Murphy et~al.(2024)Murphy, Yurchak and M{\"u}ller}]{pykrige}
\bibinfo{author}{Murphy, B.}, \bibinfo{author}{Yurchak, R.},
  \bibinfo{author}{M{\"u}ller, S.}, \bibinfo{year}{2024}.
\newblock \bibinfo{title}{Geostat-framework/pykrige: v1.7.2}.
\newblock \DOIprefix\doi{https://doi.org/10.5281/zenodo.11360184}.
\bibitem[{Nguyen et~al.(2014)Nguyen, Mahowald, Golinval and
  Maas}]{nguyen2014damage}
\bibinfo{author}{Nguyen, V.H.}, \bibinfo{author}{Mahowald, J.},
  \bibinfo{author}{Golinval, J.C.}, \bibinfo{author}{Maas, S.},
  \bibinfo{year}{2014}.
\newblock \bibinfo{title}{Damage detection in civil engineering structure
  considering temperature effect}, in: \bibinfo{booktitle}{Dynamics of {C}ivil
  {S}tructures, {V}olume 4: {P}roceedings of the 32nd {IMAC}, {A C}onference
  and {Exposition} on Structural Dynamics, 2014},
  \bibinfo{organization}{Springer}. pp. \bibinfo{pages}{187--196}.
\newblock \DOIprefix\doi{https://doi.org/10.1007/978-3-319-04546-7_22}.
\bibitem[{Nguyen et~al.(2016)Nguyen, Schommer, Maas and
  Z{\"u}rbes}]{nguyen2016static}
\bibinfo{author}{Nguyen, V.H.}, \bibinfo{author}{Schommer, S.},
  \bibinfo{author}{Maas, S.}, \bibinfo{author}{Z{\"u}rbes, A.},
  \bibinfo{year}{2016}.
\newblock \bibinfo{title}{Static load testing with temperature compensation for
  structural health monitoring of bridges}.
\newblock \bibinfo{journal}{Engineering Structures} \bibinfo{volume}{127},
  \bibinfo{pages}{700--718}.
\newblock \DOIprefix\doi{https://doi.org/10.1016/j.engstruct.2016.09.018}.
\bibitem[{Pan et~al.(2022)Pan, Liu, Wang, Jin and Chi}]{pan2022novel}
\bibinfo{author}{Pan, J.}, \bibinfo{author}{Liu, W.}, \bibinfo{author}{Wang,
  J.}, \bibinfo{author}{Jin, F.}, \bibinfo{author}{Chi, F.},
  \bibinfo{year}{2022}.
\newblock \bibinfo{title}{A novel reconstruction method of temperature field
  for thermomechanical stress analysis of arch dams}.
\newblock \bibinfo{journal}{Measurement} \bibinfo{volume}{188},
  \bibinfo{pages}{110585}.
\newblock \DOIprefix\doi{https://doi.org/10.1016/j.measurement.2021.110585}.
\bibitem[{Peeters et~al.(2001)Peeters, Maeck and
  De~Roeck}]{peeters2001vibration}
\bibinfo{author}{Peeters, B.}, \bibinfo{author}{Maeck, J.},
  \bibinfo{author}{De~Roeck, G.}, \bibinfo{year}{2001}.
\newblock \bibinfo{title}{Vibration-based damage detection in civil
  engineering: excitation sources and temperature effects}.
\newblock \bibinfo{journal}{Smart Materials and Structures}
  \bibinfo{volume}{10}, \bibinfo{pages}{518--527}.
\newblock \DOIprefix\doi{https://doi.org/10.1088/0964-1726/10/3/314}.
\bibitem[{Peng et~al.(2020)Peng, Lin, Xiang, Chen, Zhou, Yang and
  Qiao}]{peng2020positioning}
\bibinfo{author}{Peng, H.}, \bibinfo{author}{Lin, P.}, \bibinfo{author}{Xiang,
  Y.}, \bibinfo{author}{Chen, W.}, \bibinfo{author}{Zhou, S.},
  \bibinfo{author}{Yang, N.}, \bibinfo{author}{Qiao, Y.}, \bibinfo{year}{2020}.
\newblock \bibinfo{title}{A positioning method of temperature sensors for
  monitoring dam global thermal field}.
\newblock \bibinfo{journal}{Frontiers in Materials} \bibinfo{volume}{7},
  \bibinfo{pages}{587738}.
\newblock \DOIprefix\doi{https://doi.org/10.3389/fmats.2020.587738}.
\bibitem[{Potgieter(1983)}]{potgieter1983response}
\bibinfo{author}{Potgieter, I.C.}, \bibinfo{year}{1983}.
\newblock \bibinfo{title}{Response of highway bridges to nonlinear temperature
  distributions}.
\newblock \bibinfo{publisher}{University of Illinois at Urbana-Champaign}.
\bibitem[{Reynders et~al.(2014)Reynders, Wursten and
  De~Roeck}]{reynders2014output}
\bibinfo{author}{Reynders, E.}, \bibinfo{author}{Wursten, G.},
  \bibinfo{author}{De~Roeck, G.}, \bibinfo{year}{2014}.
\newblock \bibinfo{title}{Output-only structural health monitoring in changing
  environmental conditions by means of nonlinear system identification}.
\newblock \bibinfo{journal}{Structural Health Monitoring} \bibinfo{volume}{13},
  \bibinfo{pages}{82--93}.
\newblock \DOIprefix\doi{https://doi.org/10.1177/1475921713502836}.
\bibitem[{Roberts-Wollman et~al.(2002)Roberts-Wollman, Breen and
  Cawrse}]{roberts2002measurements}
\bibinfo{author}{Roberts-Wollman, C.L.}, \bibinfo{author}{Breen, J.E.},
  \bibinfo{author}{Cawrse, J.}, \bibinfo{year}{2002}.
\newblock \bibinfo{title}{Measurements of thermal gradients and their effects
  on segmental concrete bridge}.
\newblock \bibinfo{journal}{Journal of Bridge Engineering} \bibinfo{volume}{7},
  \bibinfo{pages}{166--174}.
\newblock
  \DOIprefix\doi{https://doi.org/10.1061/(asce)1084-0702(2002)7:3(166)}.
\bibitem[{Santill{\'a}n et~al.(2015)Santill{\'a}n, Salete and
  Toledo}]{santillan2015new}
\bibinfo{author}{Santill{\'a}n, D.}, \bibinfo{author}{Salete, E.},
  \bibinfo{author}{Toledo, M.}, \bibinfo{year}{2015}.
\newblock \bibinfo{title}{A new 1{D} analytical model for computing the thermal
  field of concrete dams due to the environmental actions}.
\newblock \bibinfo{journal}{Applied Thermal Engineering} \bibinfo{volume}{85},
  \bibinfo{pages}{160--171}.
\newblock \DOIprefix\doi{https://doi.org/10.1016/j.applthermaleng.2015.04.023}.
\bibitem[{Sharma and Sen(2021)}]{sharma2021bridge}
\bibinfo{author}{Sharma, S.}, \bibinfo{author}{Sen, S.}, \bibinfo{year}{2021}.
\newblock \bibinfo{title}{Bridge damage detection in presence of varying
  temperature using two-step neural network approach}.
\newblock \bibinfo{journal}{Journal of Bridge Engineering}
  \bibinfo{volume}{26}.
\newblock \DOIprefix\doi{https://doi.org/10.1061/(asce)be.1943-5592.0001708}.
\bibitem[{Sheng et~al.(2022)Sheng, Zhou, Huang, Cai and
  Shi}]{sheng2022prediction}
\bibinfo{author}{Sheng, X.}, \bibinfo{author}{Zhou, T.},
  \bibinfo{author}{Huang, S.}, \bibinfo{author}{Cai, C.}, \bibinfo{author}{Shi,
  T.}, \bibinfo{year}{2022}.
\newblock \bibinfo{title}{Prediction of vertical temperature gradient on
  concrete box-girder considering different locations in china}.
\newblock \bibinfo{journal}{Case Studies in Construction Materials}
  \bibinfo{volume}{16}, \bibinfo{pages}{e01026}.
\newblock \DOIprefix\doi{https://doi.org/10.1016/j.cscm.2022.e01026}.
\bibitem[{Soo Lon~Wah and Chen(2020)}]{soo2020new}
\bibinfo{author}{Soo Lon~Wah, W.}, \bibinfo{author}{Chen, Y.T.},
  \bibinfo{year}{2020}.
\newblock \bibinfo{title}{A new approach toward damage localization and
  quantification of structures under changing temperature condition}.
\newblock \bibinfo{journal}{Journal of Low Frequency Noise, Vibration and
  Active Control} \bibinfo{volume}{39}, \bibinfo{pages}{572--587}.
\newblock \DOIprefix\doi{https://doi.org/10.1177/1461348418793079}.
\bibitem[{Sun et~al.(2019)Sun, Zhang and Nagarajaiah}]{sun2019bridge}
\bibinfo{author}{Sun, L.}, \bibinfo{author}{Zhang, W.},
  \bibinfo{author}{Nagarajaiah, S.}, \bibinfo{year}{2019}.
\newblock \bibinfo{title}{Bridge real-time damage identification method using
  inclination and strain measurements in the presence of temperature
  variation}.
\newblock \bibinfo{journal}{Journal of Bridge Engineering}
  \bibinfo{volume}{24}, \bibinfo{pages}{04018111}.
\newblock \DOIprefix\doi{https://doi.org/10.1061/(asce)be.1943-5592.0001325}.
\bibitem[{Tong et~al.(2001)Tong, Tham, Au and Lee}]{tong2001numerical}
\bibinfo{author}{Tong, M.}, \bibinfo{author}{Tham, L.}, \bibinfo{author}{Au,
  F.}, \bibinfo{author}{Lee, P.}, \bibinfo{year}{2001}.
\newblock \bibinfo{title}{Numerical modelling for temperature distribution in
  steel bridges}.
\newblock \bibinfo{journal}{Computers \& Structures} \bibinfo{volume}{79},
  \bibinfo{pages}{583--593}.
\newblock \DOIprefix\doi{https://doi.org/10.1016/s0045-7949(00)00161-9}.
\bibitem[{{Waka Kotahi NZ Transport Agency}(2013)}]{waka2013bridge}
\bibinfo{author}{{Waka Kotahi NZ Transport Agency}}, \bibinfo{year}{2013}.
\newblock \bibinfo{title}{Bridge Manual (SP/M/022)}.
\newblock ISBN: 978-0-478-37161-1. ISBN (Online): 978-0-478-37162-8, \URLprefix
  \url{https://www.nzta.govt.nz/resources/bridge-manual/}.
  \bibinfo{note}{{T}hird edition, {A}mendment 4, Effective from May 2022.
  ({A}ccessed online on 28/10/2024)}.
\bibitem[{Wang et~al.(2023)Wang, Chai, Liu and Gu}]{wang2023causal}
\bibinfo{author}{Wang, S.}, \bibinfo{author}{Chai, B.}, \bibinfo{author}{Liu,
  Y.}, \bibinfo{author}{Gu, H.}, \bibinfo{year}{2023}.
\newblock \bibinfo{title}{A causal prediction model for the measured
  temperature field of high arch dams with dual simulation of lag influencing
  mechanism}.
\newblock \bibinfo{journal}{Structures} \bibinfo{volume}{58}.
\newblock \DOIprefix\doi{https://doi.org/10.1016/j.istruc.2023.105568}.
\bibitem[{Wang et~al.(2020)Wang, Gao and Liu}]{wang2020damage}
\bibinfo{author}{Wang, X.}, \bibinfo{author}{Gao, Q.}, \bibinfo{author}{Liu,
  Y.}, \bibinfo{year}{2020}.
\newblock \bibinfo{title}{Damage detection of bridges under environmental
  temperature changes using a hybrid method}.
\newblock \bibinfo{journal}{Sensors} \bibinfo{volume}{20},
  \bibinfo{pages}{3999}.
\newblock \DOIprefix\doi{https://doi.org/10.3390/s20143999}.
\bibitem[{Warnakulasuriya et~al.(2024)Warnakulasuriya, Antonau, Airaudo,
  L{\"o}hner, Antil and W{\"u}chner}]{warnakulasuriya2024optimal}
\bibinfo{author}{Warnakulasuriya, S.}, \bibinfo{author}{Antonau, I.},
  \bibinfo{author}{Airaudo, F.}, \bibinfo{author}{L{\"o}hner, R.},
  \bibinfo{author}{Antil, H.}, \bibinfo{author}{W{\"u}chner, R.},
  \bibinfo{year}{2024}.
\newblock \bibinfo{title}{Optimal sensor placement methodology for
  adjoint-based weakening localization in structures}, in:
  \bibinfo{booktitle}{ECCOMAS Congress 2024: 9th European Congress on
  Computational Methods in Applied Sciences and Engineering},
  \bibinfo{address}{Lisbon, Portugal}.
\bibitem[{Yang et~al.(2012)Yang, Chang and Liu}]{yang2012fem}
\bibinfo{author}{Yang, X.P.}, \bibinfo{author}{Chang, X.L.},
  \bibinfo{author}{Liu, X.H.}, \bibinfo{year}{2012}.
\newblock \bibinfo{title}{{FEM} simulation of temperature and thermal stress of
  xiaowan arch dam}.
\newblock \bibinfo{journal}{Applied Mechanics and Materials}
  \bibinfo{volume}{212}, \bibinfo{pages}{887--890}.
\newblock
  \DOIprefix\doi{https://doi.org/10.4028/www.scientific.net/amm.212-213.887}.
\bibitem[{Zhang et~al.(2019)Zhang, G{\"u}l and
  Kosti{\'c}}]{zhang2019eliminating}
\bibinfo{author}{Zhang, H.}, \bibinfo{author}{G{\"u}l, M.},
  \bibinfo{author}{Kosti{\'c}, B.}, \bibinfo{year}{2019}.
\newblock \bibinfo{title}{Eliminating temperature effects in damage detection
  for civil infrastructure using time series analysis and autoassociative
  neural networks}.
\newblock \bibinfo{journal}{Journal of Aerospace Engineering}
  \bibinfo{volume}{32}.
\newblock \DOIprefix\doi{https://doi.org/10.1061/(asce)as.1943-5525.0000987}.
\bibitem[{Zhang et~al.(2024)Zhang, Zhang, Wang and Chen}]{zhang2024temperature}
\bibinfo{author}{Zhang, W.m.}, \bibinfo{author}{Zhang, Z.h.},
  \bibinfo{author}{Wang, Z.w.}, \bibinfo{author}{Chen, B.},
  \bibinfo{year}{2024}.
\newblock \bibinfo{title}{Temperature analysis and prediction for road-rail
  steel truss cable-stayed bridges based on the structural health monitoring}.
\newblock \bibinfo{journal}{Engineering Structures} \bibinfo{volume}{315},
  \bibinfo{pages}{118476}.
\newblock \DOIprefix\doi{https://doi.org/10.1016/j.engstruct.2024.118476}.
\bibitem[{Zheng et~al.(2019)Zheng, Du and Peng}]{zheng2019simulation}
\bibinfo{author}{Zheng, D.}, \bibinfo{author}{Du, Z.K.}, \bibinfo{author}{Peng,
  S.C.}, \bibinfo{year}{2019}.
\newblock \bibinfo{title}{Simulation of two-dimensional temperature field via
  kriging method based on limited conditioning points in an arc dam}, in:
  \bibinfo{booktitle}{IOP Conference Series: Earth and Environmental Science},
  \bibinfo{organization}{IOP Publishing}. p. \bibinfo{pages}{032058}.
\newblock \DOIprefix\doi{https://doi.org/10.1088/1755-1315/304/3/032058}.
\bibitem[{Zhong et~al.(2017)Zhong, Hou and Qiang}]{zhong2017improved}
\bibinfo{author}{Zhong, R.}, \bibinfo{author}{Hou, G.p.},
  \bibinfo{author}{Qiang, S.}, \bibinfo{year}{2017}.
\newblock \bibinfo{title}{An improved composite element method for the
  simulation of temperature field in massive concrete with embedded cooling
  pipe}.
\newblock \bibinfo{journal}{Applied Thermal Engineering} \bibinfo{volume}{124},
  \bibinfo{pages}{1409--1417}.
\newblock \DOIprefix\doi{https://doi.org/10.1016/j.applthermaleng.2017.06.124}.
\bibitem[{Zhou and Yi(2013)}]{zhou2013thermal}
\bibinfo{author}{Zhou, G.D.}, \bibinfo{author}{Yi, T.H.}, \bibinfo{year}{2013}.
\newblock \bibinfo{title}{Thermal load in large-scale bridges: a
  state-of-the-art review}.
\newblock \bibinfo{journal}{International Journal of Distributed Sensor
  Networks} \bibinfo{volume}{9}, \bibinfo{pages}{217983}.
\newblock \DOIprefix\doi{https://doi.org/10.1155/2013/217983}.
\bibitem[{Zhou and Yi(2014)}]{zhou2014summary}
\bibinfo{author}{Zhou, G.D.}, \bibinfo{author}{Yi, T.H.}, \bibinfo{year}{2014}.
\newblock \bibinfo{title}{A summary review of correlations between temperatures
  and vibration properties of long-span bridges}.
\newblock \bibinfo{journal}{Mathematical Problems in Engineering}
  \bibinfo{volume}{2014}, \bibinfo{pages}{1--19}.
\newblock \DOIprefix\doi{https://doi.org/10.1155/2014/638209}.
\bibitem[{Zhou et~al.(2019)Zhou, Pan, Liang, Zhao, Zhou and
  Wang}]{zhou2019temperature}
\bibinfo{author}{Zhou, H.}, \bibinfo{author}{Pan, Z.}, \bibinfo{author}{Liang,
  Z.}, \bibinfo{author}{Zhao, C.}, \bibinfo{author}{Zhou, Y.},
  \bibinfo{author}{Wang, F.}, \bibinfo{year}{2019}.
\newblock \bibinfo{title}{Temperature field reconstruction of concrete dams
  based on distributed optical fiber monitoring data}.
\newblock \bibinfo{journal}{KSCE Journal of Civil Engineering}
  \bibinfo{volume}{23}, \bibinfo{pages}{1911--1922}.
\newblock \DOIprefix\doi{https://doi.org/10.1007/s12205-019-0787-6}.
\bibitem[{Zhouc et~al.(2016)Zhouc, Zhou, Zhao and
  Liang}]{zhouc2016optimization}
\bibinfo{author}{Zhouc, H.}, \bibinfo{author}{Zhou, Y.}, \bibinfo{author}{Zhao,
  C.}, \bibinfo{author}{Liang, Z.}, \bibinfo{year}{2016}.
\newblock \bibinfo{title}{Optimization of the temperature control scheme for
  roller compacted concrete dams based on finite element and sensitivity
  analysis methods}.
\newblock \bibinfo{journal}{Stavebn{\'\i} obzor-Civil Engineering Journal}
  \bibinfo{volume}{25}.
\newblock \DOIprefix\doi{https://doi.org/10.14311/cej.2016.03.0014}.
\bibitem[{Zuo et~al.(2015)Zuo, Hu, Li and Liu}]{zuo2015extended}
\bibinfo{author}{Zuo, Z.}, \bibinfo{author}{Hu, Y.}, \bibinfo{author}{Li, Q.},
  \bibinfo{author}{Liu, G.}, \bibinfo{year}{2015}.
\newblock \bibinfo{title}{An extended finite element method for pipe-embedded
  plane thermal analysis}.
\newblock \bibinfo{journal}{Finite Elements in Analysis and Design}
  \bibinfo{volume}{102}, \bibinfo{pages}{52--64}.
\newblock \DOIprefix\doi{https://doi.org/10.1016/j.finel.2015.05.002}.

\end{thebibliography}

\end{document}